\documentstyle[amscd,12pt]{amsart}
\numberwithin{equation}{subsection}
\theoremstyle{plain}
      \newtheorem{thm}{Theorem}[subsection]
      \newtheorem{lem}[thm]{Lemma}
      \newtheorem{sublem}[thm]{Sublemma}
      \newtheorem{prop}[thm]{Proposition}
      \newtheorem{cor}[thm]{Corollary}
      \newtheorem{conj}[thm]{Conjecture}
      \newtheorem{cond}[thm]{Condition}
      \newtheorem{defn}[thm]{Definition}
      \newtheorem{rem}[thm]{Remark}

\begin{document}
\title[$p$-adic \'etale Tate twists]
          {}
\author[K. Sato]
          {}
\date{1 October, 2006}
\thanks{ \endgraf
    Appendix A was written by Kei Hagihara}
\maketitle
\begin{center}
{\Large \bf
$p$-adic \'etale Tate twists and arithmetic duality}
\end{center}
\bigskip
\begin{center}
{\Large \bf
(Twists de Tate $p$-adiques \'etale et}
\end{center}
\begin{center}
{\Large \bf
dualit\'e arithm\'etique)}
\end{center}

\vspace{1cm}
\begin{center}
{\large \bf Kanetomo Sato}
\end{center}
\vspace{1cm}
\begin{center}
Graduate School of Mathematics
\end{center}
\begin{center}
Nagoya University
\end{center}
\begin{center}
Furo-cho, Chikusa-ku
\end{center}
\begin{center}
Nagoya 464-8602, JAPAN
\end{center}
\vspace{1cm}
\begin{center}
e-mail: kanetomo@@math.nagoya-u.ac.jp
\end{center}
\begin{center}
Tel: +81-52-789-2549 \quad Fax: +81-52-789-2829
\end{center}

\newpage
%
%
%
%
%
%
%
%
%
%
%
%
\def\AA{{\mathbb A}}
\def\abs{{\mathrm{abs}}}    
\def\ad{{\mathrm{ad}}}      
\def\an{{\mathrm{an}}}      
\def\BB{{\cal B}}           %
\def\Br{{\mathrm{Br}}}      
\def\C{{\Bbb C}}            %
\def\cd{{\mathrm{cd}}}      
\def\ch{{\mathrm{ch}}}      
\def\CH{{\mathrm{CH}}}      
\def\Chain{{\mathrm{Ch}}}   
\def\cl{{\mathrm{cl}}}      
\def\codim{{\mathrm{codim}}}
\def\Coker{{\mathrm{Coker}}}
\def\cone{{\mathrm{Cone}}}  
\def\Cor{{\mathrm{cor}}}    
\def\cores{{\mathrm{cores}}}
\def\cotor{{\mathrm{cotors}}}
\def\cris{{\mathrm{cris}}}  
\def\D{{\cal D}}
\def\DD{{\Bbb D}}
\def\dlog{d{\mathrm{log}}}  
\def\ep{\epsilon}
\def\et{{\mathrm{\acute{e}t}}}    
\def\Ext{{\mathrm{Ext}}}    
\def\F{{\mathrm{F}}}        
\def\FF{{\cal F}}           
\def\Fp{{\Bbb F}_p}         
\def\Frac{{\mathrm{Frac}}}  
\def\GG{{\cal G}}           
\def\Gal{{\mathrm{Gal}}}    
\def\gp{{\mathrm{gp}}}
\def\gr{{\mathrm{gr}}}      
\def\gys{{\mathrm{Gys}}}    
\def\h{{\mathrm{h}}}        
\def\H{{\mathrm{H}}}        
\def\Hom{{\mathrm{Hom}}}    
\def\sHom{{\cal H}om}       
\def\i{\imath}
\def\id{{\mathrm{id}}}     
\def\Image{{\mathrm{Im}}}   
\def\ind{{\mathrm{ind}}}    
\def\K{{\cal K}}            
\def\KK{{\Bbb K}}           %
\def\ker{{\mathrm{Ker}}}    
\def\L{{\Bbb L}}            %
\def\lin{{\mathrm{lin}}}    
\def\loc{{\mathrm{loc}}}    
\def\mod{\;{\mathrm{mod}}\;}    
\def\O{{\cal O}}            
\def\p{{\frak p}}           
\def\P{{\mathbb P}}
\def\pcotor{p\text{-}\cotor}
\def\Pic{{\mathrm{Pic}}}    
\def\pr{{\mathrm{pr}}}      
\def\ptor{p\text{-}\tor}    
\def\Q{{\Bbb Q}}            
\def\R{{\Bbb R}}            
\def\res{{\mathrm{res}}}    
\def\sh{{\mathrm{sh}}}      
\def\sing{{\mathrm{sing}}}  
\def\sp{{\mathrm{sp}}}      
\def\Spec{{\mathrm{Spec}}}  
\def\syn{{\mathrm{syn}}}    
\def\Symb{{\mathrm{Symb}}}  
\def\T{{\frak T}}
\def\tor{{\mathrm{tors}}}   
\def\tot{{\mathrm{t}}}      
\def\tr{{\mathrm{tr}}}      
\def\tsym{{\mathrm{tsym}}}  
\def\ur{{\mathrm{ur}}}      
\def\val{{\mathrm{val}}}    
\def\Z{{\Bbb Z}}            
\def\zar{{\mathrm{Zar}}}    
\def\ZZ{{\cal Z}}
\def\Lam{\Lambda}
\def\lam{\lambda}
\def\vare{\varepsilon}
\def\vG{\Gamma}
\def\vvG{\varGamma}
\def\ra{\rightarrow}
\def\lra{\longrightarrow}
\def\Lra{\Longrightarrow}
\def\hra{\hookrightarrow}
\def\lmt{\longmapsto}
\def\sm{\setminus}
\def\wt#1{\widetilde{#1}}
\def\wh#1{\widehat{#1}}
\def\spt{\sptilde}
\def\ol#1{\overline{#1}}
\def\ul#1{\underline{#1}}
\def\us#1#2{\underset{#1}{#2}}
\def\os#1#2{\overset{#1}{#2}}
\def\lim#1{\us{#1}{\varinjlim}}
\def\Gm{{\Bbb G}_{\hspace{-1pt}\mathrm{m}}}
\def\zp{\Z_p}
\def\qzp{\Q_p/\Z_p}
\def\qz{\Q/\Z}
\def\mwitt#1#2#3{W_{\hspace{-2pt}#2}{\hspace{1pt}}\omega_{#1}^{#3}}
\def\witt#1#2#3{W_{\hspace{-2pt}#2}{\hspace{1pt}}\Omega_{#1}^{#3}}
\def\mlogwitt#1#2#3{W_{\hspace{-2pt}#2}{\hspace{1pt}}\omega_{{#1},{\log}}^{#3}}
\def\logwitt#1#2#3{W_{\hspace{-2pt}#2}{\hspace{1pt}}\Omega_{{#1},{\log}}^{#3}}
\vspace{1cm}
\begin{quote}
{\sc Abstract}:
In this paper, we define,
for arithmetic schemes with semistable reduction,
$p$-adic objects playing the roles of Tate twists in \'etale topology,
     and establish their fundamental properties.
\end{quote}
\medskip
\vspace{1cm}
\begin{quote}
{\sc R\'esum\'e}:
Dans ce papier, nous definissions,
pour les sch\'emas arithm\'etiques \`a r\'eduction semistable,
des objets $p$-adiques jouant les r\^oles de twists \`a la Tate en
topologie \'etale,
et nous \'etablissons leurs propri\'et\'es fondamentales.
\end{quote}
\medskip
\vspace{1cm}
\begin{quote}
\par
Beilinson-Lichtenbaum axioms for motivic complexes\par
$p$-adic \'etale Tate twists\par
twisted duality for arithmetic schemes\par
\'etale sheaves of $p$-adic vanishing cycles\par
\end{quote}

\vspace{1cm}
\medskip
{\small
\tableofcontents
}
\newpage
\par
%
\vspace{4mm}

\section{\bf Introduction}\label{sect1}
\medskip
Let $k$ be a finite field of characteristic $p>0$, and
  let $X$ be a proper smooth variety over $\Spec(k)$
   of dimension $d$.
For a positive integer $m$ prime to $p$,
   we have the \'etale sheaf $\mu_m$ on $X$ consisting of
   $m$-th roots of unity.
The sheaves $\Z/m\Z(n):=\mu_m^{\otimes n}$ ($n \ge 0$),
  so called Tate twists, satisfy
   Poincar\'e duality of the following form:
There is a non-degenerate pairing of finite groups for any $i \in \Z$
$$
\begin{CD}
   \H^i_{\et}(X,\Z/m\Z(n)) \times
     \H^{2d+1-i}_{\et}(X,\Z/m\Z(d-n)) @>>> \Z/m\Z.
\end{CD}
$$
On the other hand, we have the \'etale subsheaf
   $\logwitt X r n$ ($n \ge 0, r \ge 1$) of the logarithmic part of
    the Hodge-Witt sheaf $\witt X r n$ (\cite{B}, \cite{il}).
When we put $\Z/p^r\Z(n):=\logwitt X r n[-n]$, we have
  an analogous duality fact due to Milne
    \cite{milne:duality}, \cite{milne:mot}.

In this paper, for a regular scheme $X$ which is flat of finite type over 
$\Spec(\Z)$
     and a prime number $p$,
    we construct an object $\T_r(n)_X$
    playing the role of `$\Z/p^r\Z(n)$'
   in $D^b(X_{\et},\Z/p^r\Z)$, the derived category of
      bounded complexes of \'etale $\Z/p^r\Z$-sheaves on $X$.
The fundamental idea is due to Schneider \cite{sch}, that is,
  we will grue $\mu_{p^r}^{\otimes n}$
   on $X[1/p]$ and a logarithmic Hodge-Witt sheaf on the fibers
   of characteristic $p$ to define $\T_r(n)_X$
   (cf.\ Lemma \ref{lem0-1} below).
We will further prove a duality result analogous to the above
   Poincar\'e duality.
The object $\T_r(n)_X$ is a $p$-adic analogue of
  the Beilinson-Deligne complex $\R(n)_{\D}$ on
      the complex manifold $(X\otimes_{\Z}\C)^{\an}$,
     while $\mu_{p^r}^{\otimes n}$ on $X[1/p]$
       corresponds to $(2 \pi \sqrt{-1})^n \cdot \R$
         on $(X\otimes_{\Z}\C)^{\an}$.
\subsection{Existence result}\label{sect1.1}
We fix the setting as follows.
Let $p$ be a rational prime number.
Let $A$ be a Dedekind ring whose fraction field has characteristic zero
     and which has a residue field of characteristic $p$.
We assume that
\begin{quote}
{\it every residue field of $A$ of characteristic $p$ is perfect.}
\end{quote}
Let $X$ be a noetherian regular scheme of pure-dimension
    which is flat of finite type over $B:=\Spec(A)$
      and satisfies the following condition:
\begin{quote}
{\it $X$ is a smooth or semistable family around any fiber of
      $X/B$ of characteristic $p$.}
\end{quote}
Let $j$ be the open immersion $X[1/p] \hra X$.
The first main result of this paper is the following:
\begin{thm}\label{thm0-2}
For each $n \geq 0$ and $r \ge 1$,
    there exists an object $\T_r(n)_X\in D^b(X_{\et},\Z/p^r\Z)$,
      which we call a $p$-adic \'etale Tate twist,
    satisfying the following properties$:$
\smallskip
\begin{quote}
{\bf T1 (Trivialization, cf.\ \ref{def:Z/pZ}).}
There is an isomorphism
       $t : j^*\T_r(n)_X \simeq \mu_{p^r}^{\otimes n}$.
\par
\medskip
\noindent
{\bf T2 (Acyclicity, cf.\ \ref{def:Z/pZ}).}
$\T_r(n)_X$ is concentrated in $[0,n]$, i.e.,
   the $q$-th cohomology sheaf is zero unless $0 \leq q \leq n$.
\par
\medskip
\noindent
{\bf T3 (Purity, cf.\ \ref{thm:purity}).}
For a locally closed regular subscheme $i : Z \hra X$
          of characteristic $p$ and of codimension
          $c \,(\ge 1)$, there is a Gysin isomorphism
$$
\begin{CD}
     \logwitt Z r {n-c}[-n-c] @>{\simeq}>> \tau_{\leq n+c} R i^!\T_r(n)_X
      \quad \hbox{ in } \; D^b(Z_{\et},\Z/p^r\Z).
\end{CD}
$$
\par
\medskip
\noindent
{\bf T4 (Compatibility, cf.\ \ref{thm:comp}).}
Let $i_y: y \hra X$ and $i_x: x \hra X$ be points on $X$
    with $\ch(x)=p$, $x \in \ol{ \{ y \} }$ and $\codim_X(x)=\codim_X(y)+1$.
Put $c:=\codim_X(x)$. Then the connecting homomorphism
$$
\begin{CD}
    R^{n+c-1}i_{y*}(Ri_y^!\T_r(n)_X) @>>> R^{n+c}i_{x*}(Ri_x^!\T_r(n)_X)
\end{CD}
$$
in localization theory $($cf.\ \eqref{note.3.3} below$)$
agrees with the $($sheafified\,$)$ boundary map of
    Galois cohomology groups due to Kato $($cf.\ $\S\ref{sect1}.8$ below$)$
$$
\begin{CD}
    \left. \begin{cases}
     R^{n-c+1}i_{y*}\mu_{p^r}^{\otimes n-c+1} \quad& (\ch(y)=0)\\
     i_{y*}\logwitt y r {n-c+1}  \quad& (\ch(y)=p)
    \end{cases}
    \right\}
      @>>>
         i_{x*}\logwitt x r {n-c}
\end{CD}
$$
up to a sign depending only on $(\ch(y),c)$, via Gysin isomorphisms.
Here the Gysin map for $i_y$ with $\ch(y)=0$
   is defined by the isomorphism $t$ in {\bf T1} and Deligne's cycle class
     in $R^{2c-2}i_y^!\mu_{p^r}^{\otimes c-1}$.
\par
\medskip
\noindent
{\bf T5 (Product structure, cf.\ \ref{prop:prod}).}
{\it There is a unique morphism
$$
\begin{CD}
\T_r(m)_X \otimes^{\L} \T_r(n)_X
       @>>> \T_r(m+n)_X \quad \hbox { in } \; D^-(X_{\et},\Z/p^r\Z)
\end{CD}
$$
that extends the natural isomorphism
         $\mu_{p^r}^{\otimes m} \otimes \mu_{p^r}^{\otimes n}
           \simeq \mu_{p^r}^{\otimes m+n}$ on $X[1/p]$.}
\end{quote}
\end{thm}
\noindent
If $X$ is smooth over $B$, the object
    $\T_r(n)_X$ is already considered by Schneider \cite{sch}, \S7
(see also \ref{rem:jss} below).
The properties {\bf T1}--{\bf T3} and {\bf T5}
   are $\Z/p^r\Z$-coefficient variants of the Beilinson-Lichtenbaum axioms
     on the conjectural \'etale motivic complex $\vvG(n)^{\et}_X$
      \cite{Be}, \cite{lich1}, \cite{lich2}.
More precisely, {\bf T1} (resp.\ {\bf T2}) corresponds to
    the axiom of Kummer theory for $\vvG(n)^{\et}_V$
    (resp.\ the acyclicity axiom for $\vvG(n)^{\et}_X$),
    and {\bf T3} is suggested by the purity axiom and the axiom of Kummer theory
     for $\vvG(n-c)^{\et}_Z$.
Although {\bf T4} is not among the Beilinson-Lichtenbaum axioms,
     it is a natural property to be satisfied.
We deal with this rather technical property for two reasons.
One is that
     the pair $(\T_r(n)_X,t)$ ($t$ is that in {\bf T1})
       is characterized by the properties {\bf T2}, {\bf T3} and {\bf T4}
      (see \ref{cor0-3} below).
The other is that
    we need {\bf T4} to prove
    the property {\bf T7} in the following functoriality result.
\begin{thm}\label{thm0-11}
Let $X$ be as in $\ref{thm0-2}$, and
let $Z$ be another scheme which is flat of finite type over $B$
    and for which the objects $\T_r(n)_Z$ $(n \ge 0, r \ge 1)$ are defined.
Let $f : Z \to X$ be a morphism of schemes and let
    $\psi : Z[1/p] \to X[1/p]$ be the induced morphism. Then$:$
\smallskip
\begin{quote}
{\bf T6 (Contravariant functoriality, cf.\ \ref{prop:funct}).}
There is a unique morphism
$$
\begin{CD}
f^*\T_r(n)_X
       @>>> \T_r(n)_Z \quad \hbox { in } \; D^b(Z_{\et},\Z/p^r\Z)
\end{CD}
$$
that extends the natural isomorphism
    $\psi^*\mu_{p^r}^{\otimes n} \simeq \mu_{p^r}^{\otimes n}$
     on $Z[1/p]$.
\par
\medskip
\noindent
{\bf T7 (Covariant funtoriality, cf.\ \ref{prop:trace}).}
Assume that $f$ is proper, and put
    $c:=\dim(X)-\dim(Z)$. Then there is a unique morphism
$$
\begin{CD}
   Rf_*\T_r(n-c)_Z[-2c]
       @>>> \T_r(n)_X \quad \hbox { in } \; D^b(X_{\et},\Z/p^r\Z)
\end{CD}
$$
that extends the trace morphism
   $R\psi_*\mu_{p^r}^{\otimes n-c}[-2c] \to \mu_{p^r}^{\otimes n}$
     on $X[1/p]$.
\end{quote}
\smallskip
\noindent
Furthermore, these morphisms satisfy
   a projection formula $($cf.\ $\ref{cor:proj}$ below$)$.
\end{thm}
\noindent
We will explain how we find $\T_r(n)_X$ in \S\ref{sect1.3} below.
\subsection{Arithmetic duality}\label{sect1.2}
We explain the second main result of this paper,
    the arithmetic duality for $p$-adic \'etale Tate twists.
We assume that $A$ is an algebraic integer ring,
    and that $X$ is proper over $B$.
Put $V:=X[1/p]$ and $d:=\dim(X)$.
For a scheme $Z$ which is separated of finite type over $B$,
   let $\H^{*}_{c}(Z,\bullet)$ be the \'etale cohomology with compact support
      (cf.\ \S\ref{sect10.2} below).
There is a well-known pairing
$$
\begin{CD}
\H^q_c(V,\mu_{p^r}^{\otimes n})
     \times \H^{2d+1-q}_{\et}(V,\mu_{p^r}^{\otimes d-n})
       @>>> \Z/p^r\Z,
\end{CD}
$$
and it is a non-degenerate pairing of finite groups
   by the Artin-Verdier duality
     (\cite{arver}, \cite{mazur}, \cite{milne:adual}, \cite{de}, \cite{Sp})
    and the relative Poincar\'e duality for regular schemes
       (\cite{sga4}, XVIII, \cite{Th}, \cite{fujiwara}).
We extend this duality to a twisted duality for $X$ with coefficients in
    the $p$-adic \'etale Tate twists.
A key ingredient is a global trace map
$$
\begin{CD}
\H^{2d+1}_{c}(X,\T_r(d)_X) @>>> \Z/p^r\Z,
\end{CD}
$$
which is obtained from the trace morphism in {\bf T7}
    for the structural morphism $X \to B$ and
    the classical global class field theory.
See \S\ref{sect10.2} below for details.
The product structure {\bf T5} and the global trace map give rise to a pairing
\begin{equation}\label{pair:global}
\begin{CD}
\H^q_{c}(X,\T_r(n)_X) \times
     \H^{2d+1-q}_{\et}(X,\T_r(d-n)_X)
       @>>> \Z/p^r\Z.
\end{CD}
\end{equation}
The second main result of this paper is the following:
\addtocounter{thm}{1}
\begin{thm}[{{\bf \ref{thm:main2}}}]\label{thm0-4}
The pairing \eqref{pair:global}
     is a non-degenerate pairing of finite groups
       for any $q$ and $n$ with $0 \leq n \leq d$.
\end{thm}
\noindent
A crucial point of this duality result is the non-degeneracy of a pairing
\addtocounter{equation}{1}
\begin{equation}\notag
\begin{CD}
\H^q_{\et}(X_{\Sigma},\T_r(n)_{X_{\Sigma}}) \times
     \H^{2d+1-q}_{Y}(X,\T_r(d-n)_X)
       @>>> \Z/p^r\Z,
\end{CD}
\end{equation}
which is an extension of a duality result of Niziol \cite{niziol}
    for crystalline local systems.
Here $\Sigma$ denotes the set of the closed points on $B$
     of characteristic $p$,
      $X_{\Sigma}$ denotes $\coprod_{s \in \Sigma}~X \times_B B_s$
        with $B_s$ the henselization of $B$ at $s$,
    and $Y$ denotes $X \times _B \Sigma$.
To calculate this pairing,
   we will provide an explicit formula (cf.\ \ref{lem:trace} below)
    for a pairing of \'etale sheaves of $p$-adic vanishing cycles.
\par
\medskip
We state a consequence of Theorem \ref{thm0-4}.
For an abelian group $M$, let $M_{\ptor}$ be
     the subgroup of $p$-primary torsion elements
       and let $M_{\pcotor}$ be the quotient of $M_{\ptor}$
         by its maximal $p$-divisible subgroup.
The following corollary is originally
   due to Cassels and Tate (\cite{Ca}, \cite{tate:duality}, 3.2)
     in the case that the structural morphism $X \ra B$ has a section,
    and due to Saito \cite{saito2} in the general case.
\begin{cor}
Assume
  $d=2$ and either $p \geq 3$ or $A$ has no real places.
Then $\Br(X)_{\pcotor}$ is finite and
    carries a non-degenerate skew-symmetric bilinear form with values in $\qzp$.
In particular, if $p \geq 3$ then it is alternating and
   the order of $\Br(X)_{p\text{-}\cotor}$
     is a square number.
\end{cor}
\noindent
Indeed, by a Bockstein triangle (cf.\ \S\ref{sect4.3})
    and a standard limit argument,
   Theorem \ref{thm0-4} yields a non-degenerate pairing
     of cofinitely and finitely generated $\zp$-modules
$$
\begin{CD}
\H^q(X,\T_{\qzp}(1)) \times \H^{5-q}(X,\T_{\zp}(1))
        @>>> \qzp,
\end{CD}
$$
where
$\H^*(X,\T_{\qzp}(1)):=
    \varinjlim_{r \geq 1}\; \H^*_{\et}(X,\T_r(1)_X)$ and
   $\H^*(X,\T_{\zp}(1)):=
    \varprojlim_{r \geq 1}\; \H^*_{\et}(X,\T_r(1)_X)$.
By the Kummer theory for $\Gm$
    (cf.\ \ref{prop:modp} below), one can easily check that
$$
\Br(X)_{\pcotor} \simeq
   \H^2_{\et}(X,\T_{\infty}(1)_X)_{\pcotor}
   \simeq \H^3_{\et}(X,\T_{\zp}(1)_X)_{\ptor}.
$$
Hence the corollary follows from the same argument as
     for \cite{U}, 1.5 (cf.\ \cite{Ta2})
    and the fact that the bigraded algebra
    $\bigoplus_{q,n\geq 0}~\H^q_{\et}(X,\T_{r}(n)_X)$
     with respect to the cup product
    is anti-commutative in $q$.
\subsection{Construction of $\T_r(n)_X$}\label{sect1.3}
We explain how to find $\T_r(n)_X$ satisfying the properties
     in Theorem \ref{thm0-2}.
Let $Y \subset X$ be the divisor on $X$ defined by the radical ideal of
    $(p) \subset \O_X$ and
    let $V$ be the complement $X \setminus Y=X[1/p]$.
Let $\iota$ and $j$ be as follows:
$$
\begin{CD}
     V @>{j}>> X @<{\iota}<< Y.
\end{CD}
$$
We start with necessary conditions for $\T_r(n)_X$ to exist.
\begin{lem}\label{lem0-1}
Assume that there exists an object $\T_r(n)_X \in D^b(X_{\et},\Z/p^r\Z)$
     satisfying {\bf T1}--{\bf T4}.
For a point $x \in X$, let $i_x$ be the natural map $x \hra X$.
Then$:$
\begin{enumerate}
\item[(1)]
There is an exact sequence of sheaves on $X_{\et}$
\addtocounter{equation}{1}
\begin{equation}\label{intro.1.1}
\begin{CD}
R^nj_{*}\mu_{p^r}^{\otimes n}
     @>>> \bigoplus_{y \in Y^0} i_{y*}\logwitt y r {n-1}
     @>>> \bigoplus_{x \in Y^1} i_{x*}\logwitt x r {n-2},
\end{CD}
\end{equation}
where each arrow arises from the boundary maps of Galois cohomology groups.
\item[(2)]
There is a distinguished triangle
     in $D^b(X_{\et},\Z/p^r\Z)$ of the form
\begin{equation}\label{intro.1.2}
\hspace{-30pt}
\begin{CD}
\iota_*\nu_{Y,r}^{n-1}[-n-1] @>{g}>>
     \T_r(n)_X @>{t'}>>
       \tau_{\leq n}Rj_{*}\mu_{p^r}^{\otimes n}
            @>{\sigma_{X,r}(n)}>>
            \iota_*\nu_{Y,r}^{n-1}[-n].
\end{CD}
\hspace{-30pt}
\end{equation}
Here $t'$ is induced by $t$ in {\bf T1}
      and the acyclicity property in {\bf T2}, and
    $\tau_{\leq n}$ denotes the truncation at degree $n$.
The object $\nu_{Y,r}^{n-1}$
    is an \'etale sheaf on $Y$ defined as the kernel of
       the second arrow in \eqref{intro.1.1} $($restricted onto $Y)$,
     and the arrow $\sigma_{X,r}(n)$ is
    induced by the exact sequence \eqref{intro.1.1}.
\end{enumerate}
\end{lem}
\noindent
The sheaf $\nu_{Y,r}^{n-1}$ agrees with $\logwitt Y r {n-1}$
     if $Y$ is smooth.
See \S\ref{sect2.2} below for fundamental properties of $\nu_{Y,r}^{n-1}$.
Because this lemma is quite simple,
     we include a proof here.
\begin{pf}
There is a localization distinguished triangle
    (cf.\ \eqref{DT:local} below)
\begin{equation}\label{intro.1.3}
\hspace{-30pt}
\begin{CD}
\T_r(n)_X @>{j^*}>> Rj_*j^*\T_r(n)_X @>{\delta_{U,Z}^{\loc}}>>
R\iota_*R\iota^!\T_r(n)_X[1]
           @>{\iota_*}>>\T_r(n)_X[1].
\end{CD}
\hspace{-30pt}
\end{equation}
By {\bf T1}, we have $j^*\T_r(n)_X \simeq j^*\mu_{p^r}^{\otimes n}$
    via $t$.
On the other hand, one can easily check
$$
\tau_{\leq n}(Ri_*Ri^!\T_r(n)_X[1]) \simeq \iota_*\nu_{Y,r}^{n-1}[-n]
$$
    by {\bf T3} and {\bf T4} (cf.\ \eqref{note.3.4} below).
Because the map
    $R^nj_*\mu_{p^r}^{\otimes n} \ra \iota_*\nu_{Y,r}^{n-1}$
      of cohomology sheaves
      induced by $\delta_{U,Z}^{\loc}$
    is compatible with Kato's boundary maps up to a sign
     (again by {\bf T4}),
    the sequence \eqref{intro.1.1} must be a complex
      and we obtain the morphism $\sigma_{X,r}(n)$.
Finally by {\bf T2},
    we obtain the triangle \eqref{intro.1.2}
      by truncating and shifting the triangle \eqref{intro.1.3} suitably.
The exactness of \eqref{intro.1.1} also follows from {\bf T2}.
Thus we obtain the lemma.
\end{pf}
We will prove the exactness of the sequence \eqref{intro.1.1},
   independently of this lemma,
        in \ref{lem:kato} and \ref{thm:vcyc} below.
By this exactness, we are provided with the morphism
   $\sigma_{X,r}(n)$ in \eqref{intro.1.2},
    and it turns out that any object $\T_r(n)_X \in D^b(X_{\et},\Z/p^r\Z)$
      fitting into a distinguished triangle of the form \eqref{intro.1.2}
        is concentrated in $[0,n]$.
Because $D^b(X_{\et},\Z/p^r\Z)$ is a triangulated category,
   there is at least one such $\T_r(n)_X$.
Moreover, an elementary homological algebra argument
   (cf.\ \ref{lem:CD2} (3) below)
    shows that a triple $(\T_r(n)_X,t',g)$ fitting into \eqref{intro.1.2}
     is unique up to a unique isomorphism
        (and that $g$ is determined by $(\T_r(n)_X,t')$).
Thus there is a unique pair $(\T_r(n)_X,t')$ fitting into \eqref{intro.1.2}.
Our task is to prove that this pair satisfies the listed properties,
   which will be carried out in \S\S\ref{sect4}--\ref{sect7} below.
As a consequence of Theorem \ref{thm0-2} and
     Lemma \ref{lem0-1}, we obtain
\addtocounter{thm}{3}
\begin{thm}\label{cor0-3}
The pair $(\T_r(n)_X,t)$ in $\ref{thm0-2}$ is
    the only pair that satisfies {\bf T2}--{\bf T4},
    up to a unique isomorphism
      in $D^b(X_{\et},\Z/p^r\Z)$.
\end{thm}
\subsection{Comparison with known objects}\label{sect1.1'}
We mention relations between $\T_r(n)_X$ and other
     cohomology theories (or coefficients).
Assume that $A$ is local with residue field $k$
      and $p \geq n+2$.
If $X$ is smooth over $B$, then $\iota^*\T_r(n)_X$
      is isomorphic to $R\epsilon_*S_r(n)$,
     where $S_r(n)$ denotes the syntomic complex of Fontaine-Messing
     \cite{fm} on the crystalline site
     $(X_r/W_r)_{\cris}$ with $X_r := X \otimes_A A/p^rA$ and $W_r:=W_r(k)$,
           and $\epsilon$ denotes the natural continuous map
              $(X_r/W_r)_{\cris} \ra (X_r)_{\et}$ of sites.
This fact follows from
    a result of Kurihara \cite{kurihara} (cf.\ \cite{kk:vcyc})
    and Lemma \ref{lem0-1} (see also \ref{lem:CD2} (3) below).
The isomorphism $t$ in {\bf T1}
      corresponds to the Fontaine-Messing morphism
        $R\epsilon_*S_r(n) \ra
         \tau_{\leq n}\iota^*Rj_*\mu_{p^r}^{\otimes n}$.
On the other hand, $\iota^*\T_r(n)_X$ is {\it not}
     a log syntomic complex of Kato and Tsuji (\cite{kk:ss}, \cite{ts})
      unless $n > \dim(X)$,
       because the latter object is isomorphic to
        $\tau_{\leq n}\iota^*Rj_*\mu_{p^r}^{\otimes n}$
           by a result of Tsuji \cite{tsuji2}.
Therefore $\T_r(n)_X$ is a new object particularly
      on semistable families.
\par
We turn to the setting in \S\ref{sect1.1},
     and mention what can be hoped for $\T_r(n)_X$ in comparison with
      the \'etale sheafification $\Z(n)_X^{\et}$ and
      the Zariski sheafification $\Z(n)_X^{\zar}$ of Bloch's cycle complex
       (\cite{B1}, \cite{Le1}).
By works of Levine (\cite{Le1}, \cite{levine:mot}),
     these two objects are strong candidates for
       the motivic complexes $\vvG(n)_X^{\et}$ and $\vvG(n)_X^{\zar}$,
         respectively.
So Theorem \ref{cor0-3} leads us to the following:
\begin{conj}\label{conj0-3}
\begin{enumerate}
\item[(1)]
There is an isomorphism in $D^b(X_{\et},\Z/p^r\Z)$
$$
\begin{CD}
      \Z(n)_X^{\et} \otimes^{\L} \Z/p^r\Z @>{\simeq}>> \T_r(n)_X.
\end{CD}
$$
\item[(2)]
Let $\vare$ be the natural continuous map $X_{\et} \ra X_{\zar}$ of sites.
Then the isomorphism in $(1)$ induces an isomorphism
   in $D^b(X_{\zar},\Z/p^r\Z)$
$$
\begin{CD}
\Z(n)_X^{\zar} \otimes^{\L} \Z/p^r\Z @>{\simeq}>>
       \tau_{\leq n} R\vare_* \T_r(n)_X.
\end{CD}
$$
\end{enumerate}
\end{conj}
\noindent
The case $n=0$ is obvious, because $\T_r(0)_X=\Z/p^r\Z$ (by definition).
The case $n=1$ holds by
   the Kummer theory for $\Gm$ (cf.\ \ref{prop:modp} below)
    and the isomorphisms
\begin{align*}
&\Z(1)_X^{\et} \simeq \Gm[-1], \quad
\Z(1)_X^{\zar} \simeq \vare_*\Gm[-1] \quad
   \hbox{(Levine, \cite{levine:mot}, 11.2)},\\
&R^1\vare_*\Gm =0 \quad \hbox{(Hilbert's theorem 90)}.
\end{align*}
As for $n \geq 2$,
   by results of Geisser (\cite{geisser}, 1.2 (2), (4), 1.3),
   Conjecture \ref{conj0-3} holds if $X/B$ is smooth,
   under the Bloch-Kato conjecture on Galois symbol maps \cite{bk}, \S5.
A key step in his proof is to show that
$\Z(n)_X^{\zar}\otimes^{\L} \Z/p^r\Z$ is concentrated in degrees $\le n$.
We have nothing to say about this problem for the general case
    in this paper.
\subsection{Guide for the readers}
This paper is organized as follows.
In \S\ref{sect2}, we will review some preliminary facts from
    homological algebra and results in \cite{sato:ss},
      which will be used frequently in this paper.
In \S\ref{sect3}, which is the technical heart of this paper,
we will provide preliminary results on
  \'etale sheaves of $p$-adic vanishing cycles
   (cf.\ Theorem \ref{thm:vcyc}, Corollary \ref{cor:vcyc})
   using the Bloch-Kato-Hyodo theorem (Theorem \ref{thm:hyodo}).
In \S\ref{sect4}, we will define $p$-adic \'etale Tate twists
in a slightly more general situation and prove fundamental properties including
     the product structure {\bf T5},
     the contravariant functoriality {\bf T6},
     the purity property {\bf T3} and the Kummer theory for $\Gm$.
In \S\S\ref{sect5}--\ref{sect6},
     we are concerned with the compatibility property {\bf T4}.
Using this property,
      we will prove the covariant functoriality {\bf T7}
        and a projection formula in \S\ref{sect7}.
In \S\S\ref{sect8}--\ref{sect10},
    we will study pairings of $p$-adic vanishing cycles
    and prove Theorem \ref{thm0-4}.
The appendix A due to Kei Hagihara
    includes a proof of a semi-purity of the \'etale sheaves of
     $p$-adic vanishing cycles
     (cf.\ Theorem \ref{lem:fin} below),
      which plays an important role in this paper.
He applies his semi-purity result
    to the coniveau filtration on \'etale cohomology groups
      of varieties over $p$-adic fields
       (cf.\ Theorems \ref{thm:hmain1}, \ref{thm:hmain2} and
        Corollary \ref{coro} below).
\medskip
\section*{\bf Notation and conventions}
\stepcounter{subsection}
\medskip
\noindent
1.6.
For an abelian group $M$ and a positive integer $n$,
      ${}_nM$ and $M/n$ denote the kernel and the cokernel of the map
         $M \os{\times n}{\lra} M$, respectively.
For a field $k$, $\ol k$ denotes a fixed separable closure, and
     $G_k$ denotes the absolute Galois group $\Gal(\ol k/k)$.
For a discrete $G_k$-module $M$, $\H^*(k,M)$ denote
     the Galois cohomology groups $\H^*_{\Gal}(G_k,M)$,
       which are the same as the \'etale cohomology groups of $\Spec(k)$
         with coefficients in the \'etale sheaf associated with $M$.
\par
\stepcounter{subsection}
\medskip
\smallskip
\noindent
1.7.
Unless indicated otherwise, all cohomology groups of schemes are
     taken over the \'etale topology.
We fix some general notation for a scheme $X$.
For a point $x \in X$, $\kappa(x)$ denotes its residue field
     and $\ol x$ denotes $\Spec(\ol{\kappa(x)})$.
If $X$ has pure dimension, then for a non-negative integer $q$,
      $X^q$ denotes the set of all points on $X$ of codimension $q$.
For an \'etale sheaf $\Lam$ of commutative rings on $X$,
      we write $D(X_{\et},\Lam)$
        for the derived category of \'etale $\Lam$-modules on $X$
          (cf.\ \cite{H:RD}, I, \cite{pervers}, \S1).
We write $D^+(X_{\et},\Lam)$ for the full subcategory
     of $D(X_{\et},\Lam)$ consisting of objects
         coming from complexes of \'etale $\Lam$-modules
           bounded below.
For $x \in X$ and the natural map $i_x:x \hra X$, we define the functor
     $Ri_x^!:D^+(X_{\et},\Lam) \ra D^+(x_{\et},\Lam)$
       as $\xi^* R\varphi^!$,
     where $\varphi$ denotes the closed immersion $\ol {\{ x \}} \hra X$
         and $\xi$ denotes the natural map $x \hra \ol {\{ x \}}$.
If $\xi$ is of finite type, then
   $Ri_x^!$ is right adjoint to $Ri_{x*}$,
   but otherwise it is not.
If $x$ is a generic point of $X$, $Ri_x^!$ agrees with $i_x^*$.
For $\FF \in D^+(X_{\et},\Lam)$,
    we often write $\H^*_x(X,\FF)$ for $\H^*_x(\Spec(\O_{X,x}),\FF)$.
\par
\stepcounter{subsection}
\medskip
\smallskip
\noindent
1.8.
We fix some notation of arithmetic objects defined for a scheme $X$.
For a positive integer $m$ invertible on $X$,
    $\mu_m$ denotes the \'etale sheaf of $m$-th roots of unity.
If $X$ is a smooth variety
     over a perfect field of positive characteristic $p>0$,
    then for integers $r \geq 1$ and $q \ge 0$,
      $\logwitt X r q$ denotes the \'etale subsheaf of
    the logarithmic part of the Hodge-Witt sheaf $\witt X r q$
     (\cite{B}, \cite{il}).
For $q < 0$, we define $\logwitt X r q$ as the zero sheaf.
For a noetherian excellent scheme $X$
     (all schemes in this paper are of this kind),
    we will use the following notation.
Let $y$ and $x$ be points
on $X$ such that
       $x$ has codimension $1$ in the closure $\ol {\{ y \}} \subset X$.
Let $p$ be a prime number, and
   let $i$ and $n$ be non-negative integers.
In \cite{kk:hasse}, \S1, Kato defined the boundary maps
\begin{eqnarray*}
   \H^{i+1}(y,\mu_{p^r}^{\otimes n+1})
        &  \lra \;
\H^i(x,\mu_{p^r}^{\otimes n}) \hspace{20pt} \; \quad
          & (\hbox{if } \, \ch(x)\not = p), \\
      \H^{0}(y,\logwitt y r {n+1})
        & \lra \;
   \H^0(x,\logwitt x r n) \quad
          & (\hbox{if } \, \ch(y)=\ch(x)=p), \\
   \H^{n+1}(y,\mu_{p^r}^{\otimes n+1})
        & \lra \;
       \H^0(x,\logwitt x r n)
          \quad & (\hbox{if } \, \ch(y)=0 \, \hbox{ and } \, \ch(x)=p).
\end{eqnarray*}
We write $\partial_{y,x}^{\val}$ for these maps.
See \eqref{isom:bk} for the construction of the last map.
\par
\stepcounter{subsection}
\medskip
\smallskip
\noindent
1.9.
Let $X$ be a scheme,
let $i:Z \hra X$ be a closed immersion,
     and let $j:U \hra X$ be the open complement $X \sm Z$.
Let $m$ be a non-negative integer.
For ${\cal K} \in D^+(X_{\et},\Z/m\Z)$,
     we define the morphism
$$
\begin{CD}
\delta_{U,Z}^{\loc}({\cal K}):
      Rj_*j^*{\cal K} @>>> Ri_*Ri^!{\cal K}[1]
          \quad \hbox{ in } \; D^+(X_{\et},\Z/m\Z)
\end{CD}
$$
as the connecting morphism associated with the semi-splitting
      short exact sequence of complexes
$0 \ra i_*i^!I^{\bullet} \ra I^{\bullet} \ra j_*j^*I^{\bullet} \ra 0$
       (\cite{sga4.5}, Cat\'egories D\'eriv\'ees, I.1.2.4),
         where $I^{\bullet}$ is an injective resolution of ${\cal K}$.
The morphism $\delta_{X,U}^{\loc}({\cal K})$ is functorial in ${\cal K}$,
    and
\begin{equation}\label{sign}
\begin{CD}
\delta_{U,Z}^{\loc}({\cal K})[q] =
    (-1)^q \cdot \delta_{U,Z}^{\loc}({\cal K}[q])
\end{CD}
\end{equation}
for an integer $q$.
Note also that the triangle
\begin{equation}\label{DT:local}
\begin{CD}
{\cal K} @>{j^*}>> Rj_*j^*{\cal K} @>{\delta_{U,Z}^{\loc}({\cal K})}>>
       Ri_*Ri^!{\cal K}[1]
            @>{i_*}>> {\cal K}[1]
\end{CD}
\end{equation}
is distinguished in $D^+(X_{\et},\Z/m\Z)$,
where the arrow $i_*$ (resp.\ $j^*$) denotes
      the adjunction map $Ri_*Ri^! \ra \id$
         (resp.\ $\id \ra Rj_*j^*$).
We generalize the above connecting morphism to
           the following situation.
Let $y$ and $x$ be points on $X$ such that
      $x$ has codimension $1$
       in the closure $\ol {\{ y \}} \subset X$.
Put $T:=\ol{\{ y \}} \subset X$ and $S:=\Spec(\O_{T,x})$, and
     let $i_T$ (resp.\ $i_y$, $i_x$, $\psi$)
      be the natural map $T \hra X$
      (resp.\ $y \hra X$, $x \hra X$, $S \hra T$).
Then for ${\cal K} \in D^+(X_{\et},\Z/m\Z)$,
    we define
\begin{equation}\label{note.3.3}
\begin{CD}
\delta_{y,x}^{\loc}({\cal K}):=
     Ri_{T*}R\psi_*\{ \delta_{y,x}^{\loc}(\psi^*Ri_T^!{\cal K}) \}:
      Ri_{y*}Ri_y^!{\cal K} @>>> Ri_{x*}Ri_x^!{\cal K}[1]
\end{CD}
\end{equation}
($Ri_y^!$ and $Ri_x^!$ were defined in \S1.6),
which is a morphism in $D^+(X_{\et},\Z/m\Z)$.
These connecting morphisms for all points on $X$ give rise to
    a local-global spectral sequence of sheaves on $X_{\et}$
$$
\begin{CD}
E_1^{u,v}= \bigoplus{}_{x \in X^u}~R^{u+v}i_{x*}(Ri_x^!{\cal K})
     \Lra {\cal H}^{u+v}({\cal K}).
\end{CD}
$$
For a closed immersion $i:Z \hra X$, there is
   a localized variant
\begin{equation}\label{note.3.4}
\begin{CD}
E_1^{u,v}= \bigoplus{}_{x \in X^u \cap Z}~R^{u+v}i_{x*}(Ri_x^!{\cal K})
    \Lra i_*R^{u+v}i^!{\cal K}.
\end{CD}
\end{equation}
\par
\stepcounter{subsection}
\medskip
\smallskip
\noindent
1.10.
Let $k$ be a field, and let $X$ be a pure-dimensional scheme
     which is of finite type over $\Spec(k)$.
We call $X$ a {\it normal crossing scheme} over $\Spec(k)$,
   if it is everywhere \'etale locally isomorphic to
$$
\Spec(k[T_0,T_1,\dotsb,T_N]/(T_0T_1\dotsb T_a))
$$
for some integer $a$ with $0 \leq a \leq N=\dim(X)$.
This condition is equivalent to the assumption that
$X$ is everywhere \'etale locally embedded into
       a smooth variety over $\Spec(k)$
          as a normal crossing divisor.
\par
\stepcounter{subsection}
\medskip
\smallskip
\noindent
1.11.
Let $A$
be a discrete valuation ring, and
    let $K$ (resp.\ $k$) be the fraction field (resp.\ residue field) of $A$.
Let $X$ be a pure-dimensional scheme
   which is flat of finite type over $\Spec(A)$.
We call $X$ a {\it regular semistable family} over $\Spec(A)$,
      if it is regular and
everywhere \'etale locally isomorphic to
$$
\Spec(A[T_0,T_1\dotsc,T_N]/(T_0T_1\dotsb T_a-\pi))
$$
for some integer $a$ with $0 \leq a \leq N=\dim(X/A)$,
where $\pi$ denotes a prime element of $A$.
This condition is equivalent to the assumption that
     $X$ is regular, $X \otimes_A K$ is smooth over $\Spec(K)$,
   and $X \otimes_A k$ is reduced and a normal crossing divisor on $X$.
If $X$ is a regular semistable family over $\Spec(A)$, then
    the closed fiber $X \otimes_A k$ is a normal crossing scheme over
      $\Spec(k)$.
\medskip
\section{\bf Preliminaries}\label{sect2}
\medskip
In this section we review some fundamental facts
    on homological algebra and
    results of the author in \cite{sato:ss},
    which will be used frequently in this paper.
\subsection{Elementary facts from homological algebra}\label{sect2.1}
Let ${\cal A}$ be an abelian category with enough injective objects,
    and let $D({\cal A})$ be the derived
category of complexes of objects of ${\cal A}$.
\begin{lem}\label{lem:CD}
Let $m$ and $q$ be integers.
Let $\K$ be an object
        of $D({\cal A})$ concentrated in degrees $ \leq m$
    and let $\K '$ be an object of $D({\cal A})$
      concentrated in degrees $ \geq 0$.
Then we have
$$
\begin{CD}
\Hom_{D({\cal A})}(\K, \K'[-q])
     = \begin{cases}
            \Hom_{{\cal A}}({\cal H}^m(\K),{\cal H}^0(\K'))
                   \qquad & (\hbox{if } \, q=m), \\
           0 \qquad & (\hbox{if } \, q>m), \\
       \end{cases}
\end{CD}
$$
where for $n \in \Z$ and ${\cal L} \in D({\cal A})$,
    ${\cal H}^{n}({\cal L})$ denotes
      the $n$-th cohomology object of $\cal L$.
\end{lem}
\begin{pf}
By the assumption that ${\cal A}$ has enough injectives,
   the left hand side is written as the group of morphisms
    in the homotopy category of complexes of objects of ${\cal A}$
     (\cite{sga4.5}, Cat\'egories D\'eriv\'ees, II, 2.3 (4)).
The assertion follows from this fact.
\end{pf}
\begin{lem}\label{lem:CD2}
Let ${\cal N}_1 \os{f}{\ra} {\cal N}_2 \os{g}{\ra} {\cal N}_3
       \os{h}{\ra} {\cal N}_1 [1]$
     be a distinguished triangle in $D({\cal A})$.
\begin{enumerate}
\item[(1)]
Let $i:\K \ra {\cal N}_2$ be a morphism
   with $g \circ i=0$ and
suppose that
$
\Hom_{D({\cal A})}(\K, {\cal N}_3[-1])=0.
$
Then there exists a unique morphism $i':\K \ra {\cal N}_1$
that $i$ factors through.
\item[(2)]
Let $p:{\cal N}_2 \ra \K$ be a morphism
          with $p \circ f=0$ and
            suppose that
$
\Hom_{D({\cal A})}({\cal N}_1[1], \K)=0.
$
Then there exists a unique morphism $p':{\cal N}_3 \ra \K$
       that $p$ factors through.
\item[(3)]
Suppose that $\Hom_{D({\cal A})}({\cal N}_{2},{\cal N}_1)=0$
     $($resp.\ $\Hom_{D({\cal A})}({\cal N}_{3},{\cal N}_2)=0)$.
Then relatively to
      a fixed triple $({\cal N}_1,{\cal N}_3,h)$,
      the other triple $({\cal N}_2,f,g)$ is unique
         up to a unique isomorphism, and
     $f$ $($resp.\ $g)$ is determined
       by the pair $({\cal N}_2,g)$
     $($resp.\ $({\cal N}_2,f))$.
\end{enumerate}
\end{lem}
\begin{pf}
These claims follow from the same arguments as in \cite{pervers}, 1.1.9.
The details are straight-forward and left to the reader.
\end{pf}
\subsection{Logarithmic Hodge-Witt sheaves}\label{sect2.2}
Throughout this subsection, $n$ denotes a non-negative integer
     and $r$ denotes a positive integer.
Let $k$ be a perfect field of positive characteristic $p$.
Let $X$ be a pure-dimensional scheme of finite type over $\Spec(k)$.
For a point $x \in X$, let $i_x$ be the canonical map $x \hra X$.
We define the \'etale sheaves $\nu_{X,r}^n$ and $\lam_{X,r}^n$ on $X$ as
\begin{align*}
\nu_{X,r}^n~
    & := \ker \left( \bigoplus{}_{x \in X^0}~ i_{x*} \logwitt x r n
             \os{\partial^{\val}}{\lra}
   \bigoplus{}_{x \in X^1}~ i_{x*} \logwitt x r {n-1}
   \right),\\
\lam_{X,r}^n~
    & := \Image \left( (\Gm{}_{,X})^{\otimes n}
   \os{\dlog}{\lra} \bigoplus{}_{x \in X^0}~ i_{x*} \logwitt x r n
   \right),
\end{align*}
where $\partial^{\val}$ denotes the sum of
   $\partial^{\val}_{y,x}$'s with
          $y \in X^0$ and $x \in X^1$ (cf.\ \S\ref{sect1}.8).
By definition, $\lam_{X,r}^n$ is a subsheaf of $\nu_{X,r}^n$.
If $X$ is smooth, then both
    $\nu_{X,r}^n$ and $\lam_{X,r}^n$ agree with
      the sheaf $\logwitt {X} r n$.
See also Remark \ref{rem:hyodo} (4) below.

We define the Gysin morphism for logarithmic Hodge-Witt sheaves as follows.
We define the complex of sheaves $C^{\bullet}_r(X,n)$ on $X_{\et}$ to be
$$
\begin{CD}
       \bigoplus_{x \in X^0} i_{x*} \logwitt x r n
           @>{(-1)^{n-1}\cdot\partial}>>
             \bigoplus_{x \in X^1} i_{x*} \logwitt x r {n-1}
            @>{(-1)^{n-1}\cdot\partial}>> \dotsb \\
          @>{(-1)^{n-1}\cdot\partial}>>
                \bigoplus_{x \in X^q} i_{x*} \logwitt x r {n-q}
            @>{(-1)^{n-1}\cdot\partial}>> \dotsb.
\end{CD}
$$
Here the first term is placed in degree $0$ and
     $\partial$ denotes the sum of sheafified variants of
      $\partial^{\val}_{y,x}$'s with
           $y \in X^q$ and $x \in X^{q+1}$ $($cf.\ $\S\ref{sect1}.8)$.
The fact $\partial \circ \partial=0$ is due to Kato
            $($\cite{kk:hasse}, $1.7)$.
If $X$ is a normal crossing scheme,
    this complex is quasi-isomorphic to the sheaf $\nu_{X,r}^n$ by
      \cite{sato:ss}, $2.2.5$ $(1)$.
\begin{defn}
[cf.\ loc.\ cit., 2.4.1]
\label{rem:nugysin}
Let $X$ be a normal crossing scheme over $\Spec(k)$
   and let $i: Z \hra X$ be a closed immersion of pure codimension
     $c \geq 0$.
We define the Gysin morphism
\begin{equation}\notag
\begin{CD}
   \gys_{i}^n  :  \nu_{Z,r}^{n-c}[-c] @>>>
         Ri^! \nu_{X,r}^n   \quad \hbox{ in } \; D^b(Z_{\et},\Z/p^r\Z)
\end{CD}
\end{equation}
as the adjoint morphism of the composite morphism
    in $D^b(X_{\et},\Z/p^r\Z)$
$$
\begin{CD}
i_*\nu_{Z,r}^{n-c}[-c] @>>>
i_*C^{\bullet}_r(Z,n-c)[-c] @>>> C^{\bullet}_r(X,n)
     @<{\simeq}<< \nu_{X,r}^n,
\end{CD}
$$
where the second arrow is the natural inclusion of complexes.
See also Remark $\ref{rem:newtrace}$ below.
\end{defn}
\begin{thm}[Purity, loc.\ cit., 2.4.2]
\label{thm:wpurity}
For $i: Z \hra X$ as in Definition $\ref{rem:nugysin}$,
$\gys_{i}^n$ induces an isomorphism
$\tau_{\le c } (\gys_{i}^n):\nu_{Z,r}^{n-c}[-c] \os{\simeq}{\lra}
         \tau_{\leq c}Ri^! \nu_{X,r}^n$.
\end{thm}
\noindent
We next state the duality result in loc.\ cit.
For integers $m,n \ge0$, there is a natural biadditive
      pairing of sheaves
\addtocounter{equation}{2}
\begin{equation}\label{def:prod}
\begin{CD}
   \nu_{X,r}^n \times  \lam_{X,r}^m @>>> \nu_{X,r}^{m+n}
\end{CD}
\end{equation}
induced by the corresponding pairing on the generic points of $X$
     (cf.\ loc.\ cit., 3.1.1).
\addtocounter{thm}{1}
\begin{thm}[Duality, loc.\ cit., 1.2.2]
\label{thm:milne}
Let $k$ be a finite field, and let
$X$ be a normal crossing scheme of dimension $N$
   which is proper over $\Spec(k)$.
Then$:$
\begin{enumerate}
\item[(1)]
There is a trace map
$\tr_{X} : \H^{N+1}(X,\nu_{X,r}^N)
         \ra \Z/p^r\Z$ such that for an arbitrary closed point $x \in X$
           the composite map
$$
\begin{CD}
    \H^1(x,\Z/p^r\Z) @>{\gys_{i_x}^{N}}>>
      \H^{N+1}(X,\nu_{X,r}^N) @>{\tr_X}>> \Z/p^r\Z
\end{CD}
$$
coincides with the trace map of $x$, i.e.,
the map that sends a continuous character of $G_{\kappa(x)}$
     to its value at the Frobenius substitution.
Furthermore $\tr_X$ is bijective if $X$ is connected.
\item[(2)]
For integers $q$ and $n$
      with $0 \leq n \leq N$, the natural pairing
\addtocounter{equation}{1}
\begin{equation}\label{logHW:pair}
\hspace{-70pt}
\begin{CD}
\H^q(X,\nu_{X,r}^{n}) \times \H^{N+1-q}(X,\lam_{X,r}^{N-n})
   @>{\eqref{def:prod}}>> \H^{N+1}(X,\nu_{X,r}^N)
         @>{\tr_{X}}>> \Z/p^r\Z
\end{CD}
\hspace{-70pt}
\end{equation}
is a non-degenerate pairing of finite $\Z/p^r\Z$-modules.
\end{enumerate}
\end{thm}
\noindent
We will give the definition of $\tr_X$ in Remark \ref{rem:newtrace} (4) below.
\addtocounter{thm}{1}
\begin{rem}\label{rem:newtrace}
We summarize the properties of the Gysin morphisms
    and the trace morphisms,
      which will be used in this paper.
\begin{enumerate}
\item[(1)]
The Gysin morphisms defined in Definition $\ref{rem:nugysin}$
   satisfy the transitivity property.
\item[(2)]
For $i : Z \hra X$ as in Definition $\ref{rem:nugysin}$,
   $\gys_{i}^n$ agrees with the Gysin morphism
    considered in \cite{sato:ss}, $2.4.1$,
     up to a sign of $(-1)^c$.
In particular if $X$ and $Z$ are smooth,
  then $\gys_{i}^n$ agrees with the Gysin morphism
      $\logwitt Z r {n-c}[-c] \to
         Ri^! \logwitt X r n$ of Gros
    $($\cite{gros:purity}, {\rm II.1)}
     up to the sign
     $(-1)^c$ by \cite{sato:ss}, $2.3.1$.
This fact will be used in Lemma $\ref{lem:trace2}$ below.
\item[(3)]
Let $X$ and $Z$ be normal crossing schemes over $\Spec(k)$
    of dimension $N$ and $d$, respectively, and
   let $f : Z \to X$ be a separated morphism of schemes.
We define the morphism
$$
\begin{CD}
     \tr_f : Rf_! \nu_{Z,r}^{d} [d] @>>> \nu_{X,r}^{N}[N]
       \quad \hbox{ in } \; D^b(X_{\et},\Z/p^r\Z)
\end{CD}
$$
by applying the same arguments as for \cite{jss}, Theorem $2.9$
      to the complexes $C^{\bullet}_r(Z,d)[d]$ and $C^{\bullet}_r(X,N)[N]$.
Then $\tr_f$ agrees with that in loc.\ cit., Theorem $2.9$
    up to the sign $(-1)^{N-d}$.
In particular if $X$ and $Z$ are smooth and $f$ is proper,
   then $\tr_f$ agrees with the Gysin morphism
      $Rf_*\logwitt Z r {d}[d] \to
        \logwitt X r N [N]$ due to Gros
    $($\cite{gros:purity}, {\rm II.1)}
     up to the sign $(-1)^{N-d}$.
\item[(4)]
We define the trace map $\tr_X$ in Theorem $\ref{thm:milne}$ $(1)$
  as the map induced by $\tr_f$ for $f:X \to \Spec(k)$
    and the trace map of $\Spec(k)$.
The map $\tr_X$ agrees with $(-1)^N$-times of
     the trace morphism constructed in loc.\ cit., {\rm \S3.4}.
\end{enumerate}
\end{rem}
\medskip

\section{\bf Boundary maps on the sheaves of $p$-adic vanishing
cycles}\label{sect3}
\medskip
This section is devoted to technical preparations on
     the \'etale sheaves of $p$-adic vanishing cycles.
The main results of this section are Theorem \ref{thm:vcyc}
    and Corollary \ref{cor:vcyc} below.
Throughout this section,
      $n$ and $r$ denote integers with $n \geq 0$ and $r \geq 1$.
\subsection{Milnor $K$-groups and boundary maps}\label{sect3.1}
We prepare some notation on Milnor $K$-groups.
Let $R$ be a commutative ring with unity.
We define the $0$-th Milnor $K$-group $K^M_0(R)$ as $\Z$.
For $n \geq 1$,
      we define the $n$-th Milnor $K$-group $K^M_n(R)$
      as $(R^{\times})^{\otimes n}/J$,
        where $J$ denotes the subgroup of $(R^{\times})^{\otimes n}$
           generated by elements of the form
             $x_1 \otimes \dotsb \otimes x_n$
               with $x_i + x_j=0$ or $1$ for some
                  $1 \leq i < j \leq n$.
An element $x_1 \otimes \dotsb \otimes x_n ~\mod ~J$ will be
       denoted by $\{x_1,\dotsc,x_n \}$.
Now let $L$ be a field endowed with a discrete valuation $v$.
Let $O_v$ be the valuation ring with respect to $v$,
       and let $F_v$ be its residue field.
Fix a prime element $\pi_v$ of $O_v$.
We define the homomorphism
\begin{align*}
\partial_{\pi_v}^M : &~
      K^M_n(L) \lra K^M_{n-1}(F_v) \\
\hspace{-30pt}
({\mathrm{resp.}}~ \sp_{\pi_v} : &~
      K^M_n(L) \lra K^M_{n}(F_v))
\end{align*}
by the assignment
\begin{eqnarray*}
\hspace{30pt}
\{ \pi_v, x_1, \dotsc, x_{n-1} \}
     & \mapsto \; \{ \ol {x_1}, \dotsc, \ol {x_{n-1}} \}
       &(\hbox{resp. } 0)\\
\{x_1, \dotsc, x_{n} \} & \mapsto \; 0 \qquad \qquad \;\;\,\quad
     &(\hbox{resp. } \{ \ol {x_1} , \dotsc, \ol {x_{n}} \})
\end{eqnarray*}
with each $x_i \in O_v^{\times}$ (cf.\ \cite{bt}, I.4.3).
Here for $x \in O_v^{\times}$,
       $\ol x$ denotes its residue class in $F_v^{\times}$.
The map $\partial_{\pi_v}^M$ is called the
    {\it boundary map of Milnor $K$-groups},
      and depends only on the valuation ideal $\p_v \subset O_v$.
We will denote $\partial_{\pi_v}^M$ by $\partial_{\p_v}^M$.
On the other hand,
      {\it the specialization map} $\sp_{\pi_v}$ depends on
       the choice of $\pi_v$, and
      its restriction to
          $\ker(\partial_{\p_v}^M) \subset K^M_n(L)$
      depends only on $\p_v$.
Indeed, $\ker(\partial_{\pi_v}^M)$
              is generated by the image of $(O_v^{\times})^{\otimes n}$
               and symbols of the form
                 $\{1+a,x_1,\dotsc,x_{n-1} \}$
             with $a \in \p_v$ and each $x_i \in L^{\times}$.
\subsection{Boundary map in a geometric setting}\label{sect3.2}
Let $p$ be a prime number.
Let $K$ be a henselian discrete valuation field
     of characteristic $0$ whose residue field $k$ has characteristic $p$.
Let $O_K$ be the integer ring of $K$.
Let $X$ be a regular semistable family
         over $\Spec(O_K)$ of pure dimension (cf.\ \S\ref{sect1}.11),
            or more generally,
    a scheme over $\Spec(O_K)$
      satisfying the following condition:
\begin{cond}\label{cond0}
There exist a discrete valuation subring
      $O' \subset O_K$ with $O_K/O'$ finite
        and a pure-dimensional regular semistable family
           $X'$ over $\Spec(O')$
             with $X \simeq X'\otimes_{O'}O_K$.
\end{cond}
\noindent
Later in \S\ref{sect3.4} and \S\ref{sect3.5} below,
   the extension $O_K/O'$ will be assumed to be unramified or tamely ramified.
Let $Y$ be the reduced divisor on $X$
       defined by a prime element $\pi \in O_K$,
        and let $\iota$ and $j$ be as follows:
$$
\begin{CD}
       X_K @>{j}>> X @<{\iota}<< Y.
\end{CD}
$$
In this section,
    we are concerned with the \'etale sheaf
$$
\begin{CD}
M_{r}^n
    :=\iota^*R^nj_{*}\mu_{p^r}^{\otimes n}
\end{CD}
$$
   on $Y$ and the composite map of \'etale sheaves
\stepcounter{equation}
\begin{equation}\label{def:sigma}
\begin{CD}
\hspace{-30pt}
\partial_{X,r}^n:
M_r^n
    @>>> {\bigoplus}_{y \in Y^0}~ i_{y*}i_y^* M_{r}^n
    @>{\partial^{\val}}>> {\bigoplus}_{y \in Y^0}~ i_{y*} \logwitt {y} {r} 
{n-1}.
\hspace{-30pt}
\end{CD}
\end{equation}
Here for a point $y \in Y$,
    $i_y$ denotes the canonical map $y \hra Y$.
For each $y \in Y^0$ the second arrow $\partial^{\val}$
       is defined as follows:
\begin{equation}\label{isom:bk}
\hspace{-50pt}
\begin{CD}
(i_y^* M_{r}^n)_{\ol y} \simeq K^M_n(\O_{X,\ol y}^{\sh}[1/p])/p^r
     @>{\partial_{\p _y}^M}>> K^M_{n-1}(\kappa(\ol y))/p^r
       @>{\dlog}>> \logwitt {\ol y} {r} {n-1},
\end{CD}
\hspace{-50pt}
\end{equation}
where $\p_y$ denotes the maximal ideal of
    the discrete valuation ring $\O_{X,\ol y}^{\sh}$.
The first isomorphism is due to Bloch-Kato \cite{bk}, (5.12),
   and the arrow $\partial_{\p _y}^M$ denotes
    the boundary map of Milnor $K$-groups.
We first show the following fundamental fact:
\addtocounter{thm}{2}
\begin{lem}\label{lem:kato}
The image of $\partial_{X,r}^n$
     is contained in $\nu_{Y,r}^{n-1}$.
See $\S\ref{sect2.2}$ for the definition of $\nu_{Y,r}^{n-1}$.
\end{lem}
\begin{pf}
For $x \in X_K$, let $\i_x$ be the natural map $x \hra X$.
Consider a diagram on $Y_{\et}$
$$
\begin{CD}
    @. M_r^n
         @>{a}>> \bigoplus_{x \in (X_K)^0}~
       \iota^*R^n\i_{x*}\mu_{p^r}^{\otimes n}
       @>{\partial_2}>>
        \bigoplus_{x \in (X_K)^1}~
       \iota^*R^{n-1}\i_{x*} \mu_{p^r}^{\otimes n-1}\\
     @. @. @V{\partial_1}VV  @VV{\partial_3}V\\
0 @>>> \nu_{Y,r}^{n-1}
         @>>>
           \bigoplus_{y \in Y^0}~i_{y*}\logwitt y r {n-1}
              @>{\partial_4}>>
             \bigoplus_{y \in Y^1}~ i_{y*}\logwitt y r {n-2}.
\end{CD}
$$
Here $a$ denotes the canonical adjunction map
    and each $\partial_i$ ($i=1,\dotsc,4$) is the sum of
      sheafified variants of boundary maps in \S\ref{sect1}.8.
The right square is anti-commutative
       by a result of Kato \cite{kk:hasse}, 1.7.
The upper row is a complex by the smoothness of $X_K$.
The lower row is exact by the definition of $\nu_{Y,r}^{n-1}$.
Hence we have
    $\Image(\partial_{X,r}^n)=
    \Image(\partial_1 \circ a) \subset \nu_{Y,r}^{n-1}$.
\end{pf}
\noindent
By this lemma, $\partial_{X,r}^n$ induces a map
\addtocounter{equation}{1}
\begin{equation}\label{sigma}
\begin{CD}
    \sigma_{X,r}^n: M_r^n @>>> \nu_{Y,r}^{n-1},
\end{CD}
\end{equation}
which is a geometric version of the boundary map
      of Milnor $K$-groups (modulo $p^r$).
\subsection{Bloch-Kato-Hyodo theorem}\label{sect3.3}
We give a brief review of the Bloch-Kato-Hyodo theorem
        on the structure of $M_r^n$,
   which will be useful in this and later sections.
    See also Remark \ref{rem:hyodo} below.
We define the \'etale sheaf $\K^M_{n,X_K/Y}$ on $Y$ as
     $(\iota^*j_*\O^{\times}_{X_K})^{\otimes n}/J$,
       where $J$ denotes the subsheaf generated by
          local sections of the form
             $x_1 \otimes \dotsb \otimes x_n$
             ($x_i \in \iota^*j_*\O^{\times}_{X_K}$) with
             $x_i + x_j=0$ or $1$ for some
                  $1 \leq i < j \leq n$.
There is a natural map due to Bloch and Kato \cite{bk}, (1.2)
\begin{equation}\label{symbol}
\begin{CD}
       \K ^M_{n,X_K/Y} @>>> M_r^ n,
\end{CD}
\end{equation}
which is a geometric version of Tate's Galois symbol map.
We define the filtrations $U^{\bullet}$ and
      $V^{\bullet}$ on $M_r^n$ using this map, as follows.
\stepcounter{thm}
\begin{defn}\label{def:vcyc}
\begin{enumerate}
\item[(1)]
Let $\pi$ be a prime element of $O_K$.
Let $U_{X_K}^0$ be the full sheaf $\iota^*j_*\O^{\times}_{X_K}$.
For $q \geq 1$, let $U_{X_K}^q$
         be the \'etale subsheaf of $\iota^*j_*\O^{\times}_{X_K}$
            generated by local sections of the form
              $1+\pi^q \cdot a$ with $a \in \iota^*\O_X$.
We define the subsheaf $U^q\K ^M_{n,X_K/Y}$ $(q \geq 0)$ of
       $\K ^M_{n,X_K/Y}$
         as the part generated by
           $U_{X_K}^q \otimes
       \{ \iota^*j_*\O^{\times}_{X_K} \}^{ \otimes n-1}$.
\item[(2)]
We define the subsheaf
       $U^qM_r^n$ $(q \geq 0)$ of $M_r^n$ as the image of
         $U^q\K ^M_{n,X_K/Y}$ under \eqref{symbol}.
We define the subsheaf $V^{q}M_r^n$ $(q \geq 0)$ of $M_r^n$
        as the part generated by $U^{q+1}M_r^n$
           and the image of
            $U^q\K ^M_{n-1,X_K/Y}\otimes \langle \pi \rangle$
        under \eqref{symbol}.
\end{enumerate}
\end{defn}
\begin{rem}
\begin{enumerate}
\item[(1)]
$U^{\bullet}\K ^M_{n,X_K/Y}$
    and $U^{\bullet}M_r^n$ are independent of the choice of $\pi \in O_K$
      by definition.
\item[(2)]
$V^0M_r^n$ and $V^{\bullet}M_1^n$
    are independent of the choice of $\pi\in O_K$
     by Theorem $\ref{thm:hyodo}$ below.
\end{enumerate}
\end{rem}
\noindent
To describe the graded pieces
     $\gr^q_{U/V}M_r^n:=U^qM_r^n/V^qM_r^n$ and
        $\gr^q_{V/U}M_r^n:=V^qM_r^n/U^{q+1}M_r^n$
         (especially in the case where $Y$ is not smooth),
          we introduce some
            notation from log geometry
             in \'etale topology
See \cite{kk:log} for the general framework of log schemes
          in the Zariski topology.
See also e.g., \cite{kf:log}, \S2 and \S3
          for the corresponding framework in the \'etale topology.
For a regular scheme $Z$ and
       a normal crossing divisor $D$ on $Z$,
         we define the \'etale sheaf ${\cal L}_Z(D)$
          of pointed sets on $Z$ as
$$
{\cal L}_Z(D):=\{f \in \O_Z ;~f \hbox{ is invertible outside of } D \}
       \subset \O_Z.
$$
We regard this sheaf as a sheaf of monoids by the multiplication
         of functions.
The natural inclusion ${\cal L}_Z(D) \hra \O_Z$ gives a log structure on $Z$,
and the associated sheaf ${\cal L}_Z(D)^{\gp}$ of abelian groups
        is \'etale locally generated by
          $\O^{\times}_{Z}$ and primes of $\O_Z$
               defining irreducible components of $D$.
Now we return to our situation.
Put $B:=\Spec(O_K)$ and $s:=\Spec(k)$,
      and let ${\cal L}_s$
         be the inverse image of
              ${\cal L}_B:={\cal L}_B(s)$ onto $s_{\et}$
                  in the sense of log structures.
We define the log structure ${\cal L}_Y$ on $Y_{\et}$ as follows.
By \ref{cond0}, there exist
      a discrete valuation subring $O' \subset O_K$
        and a regular semistable
          family over $B':=\Spec(O')$
          such that
       $O_K/O'$ is finite totally ramified and
            such that $X' \otimes_{O'} O_K \simeq X$.
Note that $Y$ is a normal crossing divisor on $X'$.
We fix such a pair $(O',X')$
      and define the log structure ${\cal L}_X$ on $X_{\et}$
          as that obtained from ${\cal L}_{X'}(Y)$
           by base-change (in the category of log schemes):
$$
\begin{CD}
(X,{\cal L}_{X}) :=
       (X',{\cal L}_{X'}(Y))
         \times_{(B',{\cal L}_{B'}(s))}
            (B,{\cal L}_B).
\end{CD}
$$
Finally we define ${\cal L}_Y$
       as the inverse image of ${\cal L}_{X}$
           onto $Y_{\et}$ in the sense of log structures.
Let us recall the following fundamental facts:
\begin{itemize}
\item
The log scheme $(X',{\cal L}_{X'}(Y))$
      (resp.\ $(X,{\cal L}_X)$, $(Y,{\cal L}_Y)$) is smooth over
      the log scheme $(B',{\cal L}_{B'}(s))$
        (resp.\ $(B,{\cal L}_B)$, $(s,{\cal L}_s)$)
         with respect to the natural map induced by
           the structure map $X' \ra B'$.
\item
The relative differential modules
      $\omega_{(Y,{\cal L}_Y)/(s,{\cal L}_s)}^*$ on $Y_{\et}$
         are locally free $\O_Y$-modules of finite rank and
        coincide with the modified differential modules
           $\omega_Y^{*}$ defined in \cite{hyodo}.
\item
There is a natural surjective
       homomorphism
\addtocounter{equation}{2}
\begin{equation}\label{log:surj}
\begin{CD}
     \iota^*j_*\O^{\times}_{X_K} \simeq \iota^*({\cal L}_X^{\gp})
      @>>> {\cal L}_Y^{\gp}
\end{CD}
\end{equation}
of sheaves of abelian groups on $Y_{\et}$
     (see \cite{ts}, (3.2.1) for the first isomorphism).
\end{itemize}
Let us recall further some facts relating
       log structures and differential modules.
\begin{itemize}
\item
By the definition of $\omega_{(Y,{\cal L}_Y)/(s,{\cal L}_s)}^1
       =\omega_{Y}^1$,
    there is a natural map taking the logarithmic differentials of local 
sections
      of ${\cal L}_Y^{\gp}$:
\begin{equation}\label{dlog:log2}
\begin{CD}
\dlog:  {\cal L}_Y^{\gp} @>>> \omega_{Y}^1.
\end{CD}
\end{equation}
\item
There is an analogous map
   for each $n \geq 0$ and $r>0$
\begin{equation}\label{dlog:log}
\begin{CD}
\dlog:  ({\cal L}_Y^{\gp} )^{\otimes n} @>>>
          \bigoplus_{y \in Y^0}~i_{y_*}\logwitt y r n.
\end{CD}
\end{equation}
The modified logarithmic Hodge-Witt sheaf $\mlogwitt Y r n$
       defined by Hyodo (\cite{hyodo}, (1.5))
     agrees with the image of this map.
See also Remark \ref{rem:hyodo} (4) below.
\end{itemize}
Now we state the theorems of Bloch-Kato \cite{bk}, (1.4)
       and Hyodo \cite{hyodo}, (1.6).
For local sections $x_i \in \iota^*j_*\O^{\times}_{X_K}$
      $(1 \leq i \leq n)$,
     we will denote the image of
       $\{x_1,x_2,\dotsc,x_n \} \in \K ^M_{n,X_K/Y}$
         under the symbol map \eqref{symbol}
         again by $\{x_1,x_2,\dotsc,x_n \}$,
              for simplicity.
\addtocounter{thm}{3}
\begin{thm}[{{\bf Bloch-Kato/Hyodo}}]\label{thm:hyodo}
\begin{enumerate}
\item[(1)]
The symbol map $\eqref{symbol}$ is surjective,
      that is,
      the subsheaf $U^0M_r^n$ is the full sheaf $M_r^n$
       for any $n \geq 0$ and $r>0$.
\item[(2)]
There are isomorphisms
$$
\begin{CD}
    \gr_{U/V}^0M_r^n @.~ \simeq  ~@. \mlogwitt Y r n ;@. \quad
      \{x_1,x_2,\dotsc,x_n \} \mod V^0M_r^n @. ~\mapsto ~ @. \dlog
         (\ol {x_1} \otimes \ol {x_2} \otimes  \dotsb \otimes \ol {x_n}),\\
    \gr_{V/U}^0M_r^n @. \simeq  @. \mlogwitt Y r {n-1} ;@. \quad
      \{x_1,\dotsc,x_{n-1},\pi\} \mod U^1M_r^n
        @. \mapsto @. \dlog
         (\ol {x_1} \otimes  \dotsb \otimes \ol {x_{n-1}}),
\end{CD}
$$
where for $x \in \iota^*j_*\O^{\times}_{X_K}$,
     $\ol x$ denotes its image into
       ${\cal L}_Y^{\gp}$
       via $\eqref{log:surj}$.
\item[(3)]
Let $e$ be the absolute ramification index of $K$,
    and let $r=1$.
Then for $q$ with $1 \leq q < e':=pe/(p-1)$,
       there are isomorphisms
$$
\begin{CD}
    \gr_{U/V}^qM_1^n @.~ \simeq  ~@.
      \begin{cases}
       \omega_Y^{n-1}/\BB_Y^{n-1} \qquad (p \hspace{-5pt}\not \vert q),\\
       \omega_Y^{n-1}/\ZZ_Y^{n-1} \qquad (p \vert q), \;
      \end{cases}
       \\
    \gr_{V/U}^qM_1^n @.~ \simeq  ~@.
          \hspace{-54pt}
    \omega_Y^{n-2}/\ZZ_Y^{n-2},
\end{CD}
$$
given by the following, respectively$:$
$$
{\small
\begin{CD}
\{1+\pi^q a,x_1,\dotsc,x_{n-1} \} \mod V^qM_1^n
         @. ~\mapsto ~ @.
      \begin{cases}
      \ol {a}
        \cdot \dlog(\ol {x_1}) \wedge \dotsb \wedge \dlog(\ol {x_{n-1}})
       \mod \BB_Y^{n-1} \quad (p \hspace{-5pt}\not \vert q)\\
      \ol {a}
        \cdot \dlog(\ol {x_1}) \wedge \dotsb \wedge \dlog(\ol {x_{n-1}})
       \mod \ZZ_Y^{n-1} \quad (p \vert q)\\
      \end{cases}
         \\
\{1+\pi^qa,x_1,\dotsc,x_{n-2},\pi \} \mod U^{q+1}M_1^n @. ~\mapsto ~ @.
        \hspace{-29pt}
         \ol {a}
          \cdot \dlog(\ol {x_1}) \wedge \dotsb \wedge \dlog(\ol {x_{n-2}})
            \mod \ZZ_Y^{n-2},\\
\end{CD}
}
$$
where $\BB_Y^{m}$ $($resp.\ $\ZZ_Y^m)$
     denotes the image of $d:\omega_Y^{m-1} \ra \omega_Y^{m}$
       $($resp.\ the kernel of $d:\omega_Y^m \ra \omega_Y^{m+1})$,
    $a$ denotes a local section of $\O_X$ and $\ol {a}$ denotes
     its residue class in $\O_Y$.
\item[(4)]
We have $U^qM_1^n=V^qM_1^n=0$
      for any $q \geq e'$.
\end{enumerate}
\end{thm}
\begin{rem}\label{rem:hyodo}
\begin{enumerate}
\item[(1)]
By Theorem $\ref{thm:hyodo}$ $(1)$ and $(2)$,
     the natural adjunction map
\addtocounter{equation}{2}
\begin{equation}\label{adjoint:inj}
\begin{CD}
M_r^n/U^1M_r^n   @>>>   {\bigoplus}_{y \in Y^0}~
     i_{y*}i_y^*(M_r^n/U^1M_r^n)
\end{CD}
\end{equation}
is injective.
We will use this injectivity
       to calculate the kernel of
         the map $\sigma_{X,r}^n$ defined in \eqref{sigma}.
See the proof of Theorem $\ref{thm:vcyc}$ below.
\item[(2)]
If $Y$ is smooth over $s=\Spec(k)$,
      then we have  $\mlogwitt Y r m=\logwitt Y r m$
       and $\omega_Y^m=\Omega_Y^m=\Omega_{Y/k}^m$,
     and the isomorphisms in Theorem $\ref{thm:hyodo}$ $(2)$
      yield the direct decomposition
\begin{equation}\label{smooth:direct}
\begin{CD}
    M_r^n/U^1M_r^n \simeq \logwitt Y r n \oplus
          \logwitt Y r {n-1}
\end{CD}
\end{equation}
         $($cf.\ \cite{bk}, {\rm (1.4.1.i))}.
By this decomposition, it is easy to see that
    the kernel of $\sigma_{X,r}^n$
         is generated by $U^1M_r^n$ and the image of
           $(\iota^*\O^{\times}_{X})^{\otimes n}$
              under \eqref{symbol}.
In the next subsection,
       we will extend the last fact
         to the regular semistable case,
          although the decomposition
            \eqref{smooth:direct} does not hold
              any longer in that case.
\item[(3)]
Theorem $\ref{thm:hyodo}$ $(3)$ and $(4)$
      will be used in later sections.
\item[(4)]
There are inclusions
     of \'etale sheaves $($cf.\ \cite{sato:ss}, $4.2.1)$
$$
\begin{CD}
      \lam_{Y,r}^n \subset \mlogwitt Y r n \subset \nu_{Y,r}^n.
\end{CD}
$$
These inclusions are not equalities, in general
     $($cf.\ loc.\ cit., $4.2.3)$.
If $n=\dim(Y)$, then we have
        $\mlogwitt {Y} r n =\nu_{Y,r}^n $
           by loc.\ cit., $1.3.2$.
\end{enumerate}
\end{rem}
\subsection{Structure of $\ker(\sigma_{X,r}^n)$}\label{sect3.4}
We define the \'etale subsheaf $FM_{r}^n$ of $M_r^n$
      as the part generated by $U^1M_{r}^n$
         and the image of $(\iota^*\O^{\times}_{X})^{\otimes n}$
           under \eqref{symbol}.
In the rest of this section, we are concerned with
     the map $\sigma_{X,r}^n$ in \eqref{sigma}
      and the filtration
$$
\begin{CD}
0 \subset U^1M_r^n \subset FM_r^n \subset M_r^n.
\end{CD}
$$
\begin{rem}
Clearly, $FM_r^n$ is contained in the kernel of
       $\sigma_{X,r}^n$.
\end{rem}
\noindent
The main result of this section is
    the following theorem,
      which plays an important role in later sections
       (see also Corollary \ref{cor:vcyc} below):
\begin{thm}\label{thm:vcyc}
Suppose that $X$ is a regular semistable family over $\Spec(O_K)$.
Then $\sigma_{X,r}^n$ induces an isomorphism
\addtocounter{equation}{2}
\begin{equation}\label{map_tame}
\begin{CD}
M_{r}^n/ F M_{r}^n @>{\simeq}>> \nu_{Y,r}^{n-1},
\end{CD}
\end{equation}
that is, $\sigma_{X,r}^n$ is surjective and
    $F M_{r}^n=\ker(\sigma_{X,r}^n)$.
Furthermore there is an isomorphism
\begin{equation}\label{map_diff}
\begin{CD}
F M_{r}^n/U^1 M_{r}^n @>{\simeq}>> \lam_{Y,r}^n
\end{CD}
\end{equation}
    sending the symbol $\{ x_1, x_2, \dotsc, x_n \}$
       $(x_i \in \iota^*\O^{\times}_{X})$
          to
          $\dlog(\ol{x_1} \otimes \ol{x_2}
            \otimes \dotsb \otimes \ol{x_n})$.
Here for a section $x \in \iota^*\O^{\times}_{X}$,
           $\ol x$ denotes its residue class in $\O^{\times}_{Y}$.
See $\S\ref{sect2.2}$ for the definition of $\lam_{Y,r}^n$.
\end{thm}
\addtocounter{thm}{2}
\begin{rem}\label{rem:vcyc}
Theorem $\ref{thm:vcyc}$ is not included in Theorem $\ref{thm:hyodo}$
    unless $X$ is smooth over $O_K$.
See also Remark $\ref{rem:hyodo}$ $(2)$.
In fact, $V^0M_r^n$ is not related to $FM_r^n$ directly.
However,
     Theorem $\ref{thm:hyodo}$ $(1)$ and $(2)$ play a key role
       in the proof of Theorem $\ref{thm:vcyc}$
       as the injectivity of \eqref{adjoint:inj}.
\end{rem}
\hspace{-16pt}
We first prove the following lemma, which is an essential step
       in the proof of Theorem \ref{thm:vcyc}:
\begin{lem}\label{lem:vcyc2}
Let
$$
\begin{CD}
\tau : \ker(\sigma_{X,r}^n)/U^1M_r^n @>>>
    {\bigoplus}_{y \in Y^0}~ i_{y*} \logwitt y r n
\end{CD}
$$
be the natural map induced by
     the first map in $\eqref{def:sigma}$
      and an exact sequence
$$
\begin{CD}
0 @>>> {\bigoplus}_{y \in Y^0}~ i_{y*} \logwitt {y} {r} {n} @>>>
    {\bigoplus}_{y \in Y^0}~ i_{y*}i_y^* (M_{r}^n/U^1M_{r}^n)
      @>{\partial^{\val}}>>
       {\bigoplus}_{y \in Y^0}~ i_{y*}\logwitt {y} {r} {n-1}
\end{CD}
$$
$($cf.\ Remark $\ref{rem:hyodo}$ $(2)$,
    see $\eqref{def:sigma}$ for $\partial^{\val})$.
Then $\tau$ is injective,
     and $\Image(\tau)$ is contained in $\lam_{Y,r}^n$.
\end{lem}
\begin{pf}
The injectivity of $\tau$ immediately follows from
        that of \eqref{adjoint:inj}.
We prove that $\Image(\tau)$ is contained in $\lam_{Y,r}^n$.
Since the problem is \'etale local on $Y$,
     we may assume that $Y$ has simple normal crossings on $X$.
For $y \in Y^0$, let
     $Y_y$ be the irreducible component of $Y$ whose generic point is $y$.
For $x \in Y$, let $i_x$ be the canonical map $x \hra Y$.
Let $Y^{(1)}$ (resp.\ $Y^{(2)}$)
     be the disjoint union of irreducible components
      of $Y$ (resp.\ the disjoint union of
        intersections of two distinct irreducible
          components of $Y$),
           and let $a_i:Y^{(i)} \ra Y$
             $(i=1,2)$ be the natural map.
Fix an ordering on the set $Y^0$.
There is a \v{C}ech restriction map
      $\check{r}: a_{1*}\logwitt {Y^{(1)}} r n
              \ra a_{2*}\logwitt {Y^{(2)}} r n$,
       and its kernel agrees with $\lam_{Y,r}^n$
     by \cite{sato:ss}, 3.2.1.
Our task is to prove the following two claims:
\begin{enumerate}
\item[(1)] {\it For arbitrary points $y \in Y^0$ and $x \in (Y_y)^1$,
    the composite map}
$$
\begin{CD}
\alpha_{y,x} : \ker(\sigma_{X,r}^n)
     @>{\tau_y}>> i_{y*} \logwitt {y} {r} {n}
              @>{\partial^{\val}_{y,x}}>> i_{x*} \logwitt {x} {r} {n-1}
\end{CD}
$$
{\it
is zero, where $\tau_y$ denotes the natural map induced by $\tau$.
   Consequently, $\tau$ induces a map
    $\tau':\ker(\sigma_{X,r}^n) \ra a_{1*}\logwitt {Y^{(1)}} r n$}.
\item[(2)] {\it The following composite map is zero}:
$$
\begin{CD}
\beta : \ker(\sigma_{X,r}^n)
     @>{\tau'}>> a_{1*}\logwitt {Y^{(1)}} r n
              @>{\check{r}}>> a_{2*}\logwitt {Y^{(2)}} r n.
\end{CD}
$$
\end{enumerate}
\medskip
\noindent
{\it Proof of Claim} (1).
It suffices to show that the stalk
    $(\alpha_{y,x})_{\ol x}$
        is the zero map.
Let $Y_{\sing}$
       be the singular locus of $Y$.
The case $x \not\in (Y_{\sing})^0$ immediately follows from
      the direct decomposition \eqref{smooth:direct}.
To show the case $x \in (Y_{\sing})^0$,
   we fix some notation.
Put $R:=\O_{X,\ol x}^{\sh}$,
     which is a strict henselian regular local ring of dimension $2$.
Let $T_1$ and $T_2$ be the irreducible components
       of $\Spec(\O_{Y,\ol x}^{\sh})$.
We suppose that $T_1$ lies above $Y_y$.
Fix a prime element $t_i \in R$ $(i=1,2)$
           defining $T_i$.
Put $w_1:= t_1 ~\mod ~(t_2) \in R/(t_2)$
     and $w_2:= t_2 ~\mod ~(t_1) \in R/(t_1)$.
Because the divisor $T_1 \cup T_2 \subset \Spec(R) $ has
     simple normal crossings,
       $R/(t_1)$ and $R/(t_2)$ are
         discrete valuation rings and
      $w_1$ (resp.\ $w_2$) is a prime element
           in $R/(t_2)$ (resp.\ $R/(t_1)$).
Let $\eta_i$ ($i=1,2$) be the generic point of $T_i$.
There is a commutative diagram with exact rows
\addtocounter{equation}{2}
\begin{equation}\label{vcyc:cd2}
\hspace{-80pt}
{\small
\begin{CD}
\ker(\partial)  @>>>  K_n^M(R[1/p])/p^r @>{\partial}>>
        K_{n-1}^M(\kappa(\eta_1))/p^r \oplus
           K_{n-1}^M(\kappa(\eta_2))/p^r  \\
@VVV @V{\eqref{symbol}}VV @V{\dlog}V{\simeq}V \\
\ker(\sigma_{X,r}^n)_{\ol x}  @>{\subset}>>
      (M_r^n)_{\ol x}  @>{(\sigma_{X,r}^n)_{\ol x}}>>
         (i_{y*}\logwitt {y} r {n-1})_{\ol x} \oplus
           (i_{y' *}\logwitt {y'} r {n-1})_{\ol x}.
\end{CD}
}
\hspace{-50pt}
\end{equation}
Here $\partial$ is the direct sum of
      the boundary maps of Milnor $K$-groups modulo $p^r$, and
         $y'$ denotes the generic point of $Y$
           corresponding to $T_2$.
In this diagram,
    the right vertical arrow is bijective by
                 a theorem of Bloch-Gabber-Kato \cite{bk}, (2.1), and
       the central vertical map is surjective
          by Theorem \ref{thm:hyodo} (1).
Hence the left vertical map,
        resulting from the right square, is surjective.
On the other hand, there is a composite map
\begin{equation}\label{vcyc:modmap}
\hspace{-50pt}
\begin{CD}
\sp_{t_1,R[1/p]}:K_n^M(R[1/p])/p^r @>>> K_n^M(\Frac(R))/p^r
       @>{\sp_{t_1}}>> K_n^M(\kappa(\eta_1))/p^r.
\end{CD}
\hspace{-50pt}
\end{equation}
See \S\ref{sect3.1} for $\sp_{t_1}$.
The restriction of this map to
           $\ker(\partial)$ fits into a commutative diagram
\begin{equation}\label{vcyc:cd3}
\begin{CD}
\ker(\partial)
     @>{\sp_{t_1,R[1/p]} \vert_{\ker(\partial)}}>>
      K_n^M(\kappa(\eta_1))/p^r @>{\partial^M_{(w_2)}}>>
      K_{n-1}^M(\kappa(\ol x))/p^r \\
@V{\mathrm{surj.}}VV @V{\dlog}VV  @V{\dlog}VV \\
   \ker(\sigma_{X,r}^n)_{\ol x} @>{(\tau_y)_{\ol x}}>>
     (i_{y*}\logwitt {y} r {n})_{\ol x}
        @>{(\partial^{\val}_{y,x})_{\ol x}}>> \logwitt {\ol x} r {n-1},
\end{CD}
\end{equation}
where the composite of the lower row gives
      $(\alpha_{y,x})_{\ol x}$.
The composite of
    the upper row is the zero map
      by a commutative diagram
      of Milnor $K$-groups modulo $p^r$
\begin{equation}\label{vcyc:cd1}
\begin{CD}
K_n^M(R[1/p])/p^r @>{\partial^M_{(t_2)}}>>
      K_{n-1}^M(\kappa(\eta_2))/p^r \\
     @V{\sp_{t_1,R[1/p]}}VV  @VV{\sp_{w_1}}V \\
K_n^M(\kappa(\eta_1))/p^r @>{\partial^M_{(w_2)}}>>
    K_{n-1}^M(\kappa(\ol x))/p^r,
\end{CD}
\end{equation}
whose commutativity is shown
      explicitly by the direct decomposition
$R[1/p]^{\times} \simeq
               R^{\times} \times \langle t_1 \rangle
                 \times \langle t_2 \rangle$.
Hence $(\alpha_{y,x})_{\ol x}$ is the zero map
    by the diagram \eqref{vcyc:cd3}, and we obtain the claim (1).
\par
\medskip
\noindent
{\it Proof of Claim} (2).
Let $Z$ be a connected component of $Y^{(2)}$.
Let $Y_1$ and $Y_2$
    be the irreducible components of $Y$
          such that
           $a_2(Z) \subset Y_1 \cap Y_2$.
Our task is to show that the composite map
$$
\begin{CD}
\beta_Z : \ker(\sigma_{X,r}^n)
     @>{\tau'}>> a_{1*}(\logwitt {Y_1} r n
                 \oplus \logwitt {Y_2} r n)
              @>{\check{r}}>> a_{2*}\logwitt {Z} r n
\end{CD}
$$
is zero, where the last map sends $(\omega_1,\omega_2)$
       ($\omega_i \in a_{1*}\logwitt {Y_i} r n$)
     to $\omega_1 \vert_Z - \omega_2 \vert_Z$.
Let $x$ be the generic point of $a_2(Z)$.
Since the canonical map
     $a_{2*}\logwitt {Z} r n \ra i_{x*}\logwitt x r n$ is injective,
we have only to show that the stalk $(\beta_Z)_{\ol x}$ is zero.
Put $R:=\O_{X,\ol x}^{\sh}$, and
     let the notation be as in
        the proof of Claim (1).
Suppose that $T_i$ ($i=1,2$) is the irreducible component of
       $\Spec(\O_{Y, \ol x}^{\sh})$
       lying above $Y_i$.
Let $N_i$ $(i=1,2)$ be the kernel of
       the boundary map $K_n^M(\kappa(\eta_i))/p^r \ra
         K_{n-1}^M(\kappa(\ol x))/p^r$.
By the commutative diagram \eqref{vcyc:cd1},
    $\sp_{t_i, R[1/p]}$ in \eqref{vcyc:modmap} induces a map
$$
\begin{CD}
f_i: \ker(\partial) @>>> N_i.
\end{CD}
$$
See \eqref{vcyc:cd2} for $\partial$.
This map fits into a commutative diagram
\begin{equation}\label{vcyc:cd5}
\hspace{-30pt}
\begin{CD}
\ker(\partial)
     @>{(f_1,f_2)}>>
        N_1 \oplus N_2 @>{f_3}>>
             K_{n}^M(\kappa(\ol x))/p^r \\
@V{\mathrm{surj.}}VV @V{\dlog}VV  @V{\dlog}VV \\
   \ker(\sigma_{X,r}^n)_{\ol x} @>{\tau'}>>
     (\logwitt {T_1} r {n} \oplus \logwitt {T_2} r {n})_{\ol x}
       @>{\check{r}}>> \logwitt {\ol x} r {n},
\end{CD}
\hspace{-30pt}
\end{equation}
where $f_3$ sends $(u_1,u_2)$
        ($u_i \in N_i$) to
            $\sp_{w_2}(u_1)-\sp_{w_1}(u_2)$.
See \eqref{vcyc:cd2} for the left vertical map.
The composite of the lower row gives $(\beta_Z)_{\ol x}$.
The composite of the upper row is zero by a commutative diagram
$$
\begin{CD}
K_n^M(R[1/p])/p^r @>{\sp_{t_2,R[1/p]} }>>
      K_{n}^M(\kappa(\eta_2))/p^r\\
     @V{ \sp_{ t_1, R[1/p] } }VV  @VV{\sp_{w_1}}V \\
K_n^M(\kappa(\eta_1))/p^r @>{\sp_{w_2}}>>
    K_{n}^M(\kappa(\ol x))/p^r,
\end{CD}
$$
whose commutativity is checked in the same way as
      for \eqref{vcyc:cd1}.
Thus $(\beta_Z)_{\ol x}$ is the zero map, and
     we obtain the claim (2) and Lemma \ref{lem:vcyc2}.
\end{pf}
\noindent
{\it Proof of Theorem \ref{thm:vcyc}.}
The surjectivity of \eqref{map_tame}
    follows from the same argument
       as for \cite{sato:ss}, 2.4.6
        (see the surjectivity of the map (2.4.9) in loc.\ cit.).
We prove the injectivity of \eqref{map_tame} and
      construct the bijection \eqref{map_diff}.
There are injective maps
$$
F M_r^n/U^1 M_r^n \hra
   \ker(\sigma_{X,r}^n)/U^1M_r^n
    \os{\tau}{\hra} \lam_{Y,r}^n
$$
(see Lemma \ref{lem:vcyc2} for $\tau$).
These two arrows are both bijective,
     because the sheaves $F M_r^n/U^1 M_r^n$ and $\lam_{Y,r}^n$ are
       generated by symbols from $(\iota^*\O^{\times}_{X})^{\otimes n}$
          and $(\O^{\times}_{Y})^{\otimes n}$, respectively.
Therefore we have $F M_r^n =  \ker(\sigma_{X,r}^n)$
      as subsheaves of $M_r^n$ and
     the composite of the above two maps gives the desired bijective map
      \eqref{map_diff}.
This completes the proof of Theorem \ref{thm:vcyc}.
\subsection{Tamely ramified case}\label{sect3.5}
Assume that $X$ satisfies the following condition over $O_K$:
\begin{cond}\label{cond1}
There exist a discrete valuation subring $O' \subset O_K$
         with $O_K/O'$ finite {\it tamely} ramified
       and a regular semistable family $X'$ over $O'$
         with $X \simeq X' \otimes_{O'} O_K$.
\end{cond}
\noindent
Let $Y$ and $M_r^n$ (resp.\ $U^1M_r^n$, $FM_r^n$)
     be as we defined in \S\ref{sect3.2} (resp.\ \S\ref{sect3.3},
\S\ref{sect3.4}).
By Theorems \ref{thm:hyodo} and \ref{thm:vcyc}, we obtain
\begin{cor}\label{cor:vcyc}
The map $\sigma_{X,r}^n$ induces an isomorphism
     $M_r^n/FM_r^n \simeq \nu_{Y,r}^{n-1}$,
and there is an isomorphism
       $FM_{r}^n/U^1M_{r}^n \simeq \lam_{Y,r}^n$
       described in the same way as $\eqref{map_diff}$.
\end{cor}
\begin{pf}
The second assertion is an immediate consequence of
     Remark \ref{rem:hyodo} (1) for $X$ and
           the definitions of $FM_r^n$ and $\lam_{Y,r}^n$
            (cf.\ Lemma \ref{lem:vcyc2}).
We prove the first assertion.
Since the problem is \'etale local on $Y$,
     we may assume that $O_K/O'$ is totally tamely ramified.
Then the divisor on $X'$
         defined by a prime element
            $\pi' \in O'$ agrees with $Y$.
Let $e_1$ be the ramification index of $O_K/O'$ and
      let $\iota'$ and $j'$ be as follows:
$$
\begin{CD}
        X'_{K'} @>{j'}>> X' @<{\iota'}<< Y,
\end{CD}
$$
where $K'$ denotes $\Frac(O')$.
Let $M_{r,X'}^n$ be the \'etale sheaf
       $\iota'{}^*R^nj'_*\mu_{p^r}^{\otimes n}$.
There is a commutative diagram
      with exact rows
$$
\begin{CD}
0 @>>> FM_{r,X'}^n/U^1M_{r,X'}^n
      @>>> M_{r,X'}^n/U^1M_{r,X'}^n
        @>{\sigma_{X',r}^n ~\mod~U^1M_{r,X'}^n}>> \nu_{Y,r} @>>> 0\\
@.  @VVV @V{\res}VV @V{\times e_{1}}VV \\
0 @>>> \ker(\sigma_{X,r}^n)/U^1M_{r}^n
       @>>> M_{r}^n/U^1M_{r}^n
         @>{\sigma_{X,r}^n~\mod~U^1M_{r}^n}>> \nu_{Y,r},
\end{CD}
$$
where the exactness of the upper row follows from
      Theorem \ref{thm:vcyc}.
Because $(e_{1},p)=1$ by assumption,
     the central vertical arrow is bijective
      by Theorem \ref{thm:hyodo} (1), (2)
         (for $X'$ and $X$), and
         the right vertical arrow is bijective as well.
Hence $\sigma_{X,r}^n$ is surjective
     and the left vertical map is bijective,
     which implies the equality $FM_r^n = \ker(\sigma_{X,r}^n)$.
\end{pf}

\medskip
\section{\bf $p$-adic \'etale Tate twists}\label{sect4}
\medskip
In this section, we define
       the objects $\T_r(n)_X$ ($n \geq 0$, $r \ge 1$)
        stated in Theorem \ref{thm0-2}
         and discuss their fundamental properties including
     {\bf T1}, {\bf T2}, {\bf T3}, {\bf T5} and {\bf T6}.
\subsection{Setting}\label{sect4.1}
Let $A$ be a Dedekind ring
      whose fraction field has characteristic zero and
        which has a maximal ideal of positive characteristic.
Let $p$ be a prime number which is not invertible in $A$,
     and we assume that the residue fields of $A$
       at maximal ideals
        of characteristic $p$ are {\it perfect}.
Put $B:=\Spec(A)$ and write $\Sigma$ for
          the set of the closed points on $B$ of characteristic $p$.
For a point $s$ on $B$, let $B_s$ be
                 the henselization of $B$ at $s$.
Let $X$ be a pure-dimensional scheme
      which is flat of finite type over $B$.
We assume that $X$
       satisfies the following condition,
          unless mentioned otherwise:
\begin{cond}\label{cond2}
$X[1/p]$ is regular. For any $s \in \Sigma$,
    each connected component $X'$ of
            $X \times_B B_s$ satisfies
              the condition $\ref{cond1}$
               over the integral closure of $B_s$
                in $\vG(X',\O_{X'})$.
\end{cond}
\noindent
We will often work under the following stronger assumption:
\begin{cond}\label{cond2'}
$X$ is regular.
For any $s \in \Sigma$,
        each connected component $X'$ of
            $X \times_B B_s$ is a regular semistable family
          over the integral closure of $B_s$
                in $\vG(X',\O_{X'})$.
\end{cond}
\noindent
Let $X$ be a pure-dimensional flat of finite type $B$-scheme
       satisfying \ref{cond2}.
Let $Y \subset X$ be the divisor
     defined by the radical of $(p) \subset \O_X$.
We always assume that $Y$ is non-empty.
Let $\nu_{Y,r}^n$ be as in \S\ref{sect2.2}.
Put $V : = X \setminus Y =X[1/p]$.
Let $\iota$ and $j$ be as follows:
$$
\begin{CD}
    V @>{j}>> X @<{\iota}<< Y.
\end{CD}
$$
Define the \'etale sheaf $M_r^n$ on $Y$ to be
       $\iota^*R^nj_*\mu_{p^r}^{\otimes n}$.
\subsection{Definition of $\T_r(n)_{X}$}\label{sect4.2}
Let $X$ and $p$ be as before.
We define $\T_r(0)_X := \Z/p^r\Z_X$.
For $n \geq 1$, let
\begin{equation}\label{sigma(n)}
\begin{CD}
\sigma_{X,r}(n)
     :\tau_{\leq n}Rj_{*}\mu_{p^r}^{\otimes n}
           @>>> \iota_*\nu_{Y,r}^{n-1}[-n]
            \quad \hbox { in } \; D^b(X_{\et},\Z/p^r\Z)
\end{CD}
\end{equation}
    be the morphism induced by the map
$\iota_*(\sigma_{X,r}^n):R^nj_{*}\mu_{p^r}^{\otimes n}
          = \iota_*M_r^n \ra \iota_*\nu_{Y,r}^{n-1}$
    of sheaves on $X_{\et}$ (cf.\ Lemma \ref{lem:CD}).
See \eqref{sigma} for $\sigma_{X,r}^n$.
\stepcounter{thm}
\begin{lem}\label{lem:defdef}
Suppose $n \geq 1$, and let
\stepcounter{equation}
\begin{equation}\label{defdef}
\begin{CD}
\iota_*\nu_{Y,r}^{n-1}[-n-1] @>{g}>>
     \K @>{t}>>
       \tau_{\leq n}Rj_{*}\mu_{p^r}^{\otimes n}
            @>{\sigma_{X,r}(n)}>>
            \iota_*\nu_{Y,r}^{n-1}[-n]
\end{CD}
\end{equation}
be a distinguished triangle in $D^b(X_{\et},\Z/p^r\Z)$.
Then $\K$ is concentrated in $[0,n]$,
    the triple $(\K,t,g)$ is unique up to a unique
     isomorphism and
       $g$ is determined by
       the pair $(\K,t)$.
\end{lem}
\begin{pf}
The map $\iota_*(\sigma_{X,r}^n)$ is surjective
     by Theorem \ref{thm:vcyc} and Corollary \ref{cor:vcyc}.
Hence $\K$ is acyclic outside of $[0,n]$ and
     there is no non-zero morphism from
       $\K$ to $\iota_*\nu_{Y,r}^{n-1}[-n-1]$
         by Lemma \ref{lem:CD}.
The uniqueness assertion
     follows from this fact and Lemma \ref{lem:CD2} (3).
\end{pf}
\stepcounter{thm}
\begin{defn}\label{def:Z/pZ}
For $n \geq 1$,
     we fix a pair $(\K,t)$ fitting into
     a distinguished triangle of the form \eqref{defdef},
     and define $\T_r(n)_X:=\K$.
The morphism $t$ determines an isomorphism $j^*\K \simeq \mu_{p^r}^{\otimes
n}$,
    and $\T_r(n)_X$ is concentrated in $[0,n]$,
     that is, $\T_r(n)_X$ satisfies {\bf T1} and {\bf T2} in $\ref{thm0-2}$.
Moreover, $t$ induces isomorphisms
\stepcounter{equation}
\begin{equation}\label{defdef2}
\begin{CD}
{\cal H}^q(\T_r(n)_X) \simeq
       \begin{cases}
          R^qj_{*}\mu_{p^r}^{\otimes n} \qquad & (0 \leq q < n),\\
          \iota_*FM_r^n  \qquad & (q=n),
       \end{cases}
\end{CD}
\end{equation}
where we have used Theorem $\ref{thm:vcyc}$ and Corollary $\ref{cor:vcyc}$
       for $q=n$.
\end{defn}
We prove here the existence of a natural product structure
    ({\bf T5} in $\ref{thm0-2}$).
\stepcounter{thm}
\begin{prop}[Product structure]\label{prop:prod}
For $m, n \geq 0$, there is a unique morphism
\stepcounter{equation}
\begin{equation}\label{def:prod2}
\begin{CD}
     \T_r(m)_X \otimes^{\L} \T_r(n)_X
         @>>> \T_r(m+n)_X
          \quad \hbox{ in } \; D^-(X_{\et},\Z/p^r\Z)
\end{CD}
\end{equation}
that extends
      the natural map $\mu_{p^r}^{\otimes m} \otimes \mu_{p^r}^{\otimes n}
           \ra \mu_{p^r}^{\otimes m+n}$ on $V_{\et}$.
\end{prop}
\begin{pf}
If $m=0$ or $n=0$, then the assertion is obvious.
Assume $m, n \geq 1$, and
   put ${\frak L}:=\T_r(m)_X {\otimes}^{\L}~ \T_r(n)_X$.
By the definition of $\T_r(m+n)_X$
      and Lemma \ref{lem:CD2} (1),
it suffices to show that the following composite morphism is zero
     in $D^-(X_{\et},\Z/p^r\Z)$:
$$
\begin{CD}
{\frak L} @.\; \lra \; @. \tau_{\leq m}Rj_*\mu_{p^r}^{\otimes m}
        \otimes^{\L}  \tau_{\leq n}Rj_*\mu_{p^r}^{\otimes n}
        @.\; \lra \; @. \tau_{\leq m+n}Rj_*\mu_{p^r}^{\otimes m+n}
        @>{\sigma_{X,r}(m+n)}>>
           \iota_*\nu_{Y,r}^{m+n-1}[-m-n],
\end{CD}
$$
where the second arrow is induced by the natural map
      $\mu_{p^r}^{\otimes m} \otimes \mu_{p^r}^{\otimes n}
           \ra \mu_{p^r}^{\otimes m+n}$ on $V$.
We prove this triviality.
Because ${\frak L}$ is concentrated in degrees $ \leq m+n$,
     this composite morphism is determined by
        the composite map of the $(m+n)$-th
        cohomology sheaves (cf.\ Lemma \ref{lem:CD})
$$
\begin{CD}
{\cal H}^{m+n}({\frak L}) @>>> \iota_*M_r^m \otimes
        \iota_*M_r^n
          @>>> \iota_*M_r^{m+n}
             @>{\iota_*(\sigma_{X,r}^{m+n})}>>  \iota_*\nu_{Y,r}^{m+n-1}.
\end{CD}
$$
The image of
${\cal H}^{m+n}({\frak L}) (\simeq
     \iota_*FM_r^m \otimes \iota_*FM_r^n)$
       into $\iota_*M_r^{m+n}$
     is contained in $\iota_*FM_r^{m+n}$.
Hence this composite map is zero
     and we obtain Proposition \ref{prop:prod}.
\end{pf}
The following proposition
    ({\bf T6} in $\ref{thm0-11}$) follows from a similar argument as for
      Proposition \ref{prop:prod}.
\stepcounter{thm}
\begin{prop}[Contravariant functoriality]
\label{prop:funct}
Let $X$ and $Z$ be flat $B$-schemes
     satisfying $\ref{cond2}$.
Let $f:Z \ra X$ be a morphism of schemes, and
   let $\psi : Z[1/p] \ra X[1/p]$ be the induced morphism.
Then there is a unique morphism
$$
\begin{CD}
f^* : f^*\T_r(n)_X @>>> \T_r(n)_Z
     \quad \hbox{ in } \; D^b(Z_{\et},\Z/p^r\Z)
\end{CD}
$$
that extends the natural isomorphism
      $\psi^*\mu_{p^r}^{\otimes n} \simeq \mu_{p^r}^{\otimes n}$
       on $(Z[1/p])_{\et}$.
Consequently, these pull-back morphisms satisfy the transitivity property.
\end{prop}
\subsection{Bockstein triangle}\label{sect4.3}
We prove the following proposition:
\begin{prop}\label{prop:bock}
For $r, s \geq 1$, the following holds$:$
\begin{enumerate}
\item[(1)]
There is a unique morphism
    $\ul p : \T_r(n)_X \ra \T_{r+1}(n)_X$ in $D^b(X_{\et},\Z/p^{r+1}\Z)$
   that extends the natural inclusion
    $\mu_{p^r}^{\otimes n} \hra \mu_{p^{r+1}}^{\otimes n}$ on $V_{\et}$.
\item[(2)]
There is a unique morphism ${\cal R} : \T_{r+1}(n)_X \ra \T_r(n)_X$
   in $D^b(X_{\et},\Z/p^{r+1}\Z)$ that extends the natural projection
   $\mu_{p^{r+1}}^{\otimes n} \ra \mu_{p^{r}}^{\otimes n}$ on $V_{\et}$.
\item[(3)]
There is a canonical Bockstein morphism
   $\delta_{s,r}: \T_s(n)_X \ra \T_r(n)_X[1]$
   in $D^b(X_{\et})$
satisfying
\begin{itemize}
\item[(3-1)]
   $\delta_{s,r}$ extends the Bockstein morphism
     $\mu_{p^s}^{\otimes n} \ra \mu_{p^{r}}^{\otimes n}[1]$
   in $D^b(V_{\et})$ associated with the short exact sequence
        $0 \ra \mu_{p^{r}}^{\otimes n}
           \ra \mu_{p^{r+s}}^{\otimes n}
            \ra \mu_{p^{s}}^{\otimes n} \ra 0$.
\item[(3-2)]
    $\delta_{s,r}$ fits into a distinguished triangle
$$
\begin{CD}
       \T_{r+s}(n)_X
          @>{{\cal R} ^r}>> \T_s(n)_X
       @>{\delta_{s,r}}>> \T_r(n)_X[1]
          @>{\ul p ^s[1]}>>  \T_{r+s}(n)_X[1].
\end{CD}
$$
\end{itemize}
\end{enumerate}
\end{prop}
\begin{pf}
The claims (1) and (2) follow from
      the fact that
           $\T_r(n)_X$ concentrated in $[0,n]$
      and Lemma \ref{lem:CD2} (1).
The details are straight-forward and left to the reader.
We prove (3).
For two complexes
    $M^{\bullet}=(\{M^u \}_{u \in \Z}, \{ d_M^u:M^u \to M^{u+1}\}_{u \in \Z})$,
    $N^{\bullet}=(\{N^v \}_{v \in \Z}, \{ d_N^v \}_{v \in \Z})$
    and a map $h^{\bullet}:M^{\bullet} \ra N^{\bullet}$ of complexes,
     let $\cone(h)^{\bullet}$ be the mapping cone
      (cf.\ \cite{sga4}, XVII)
$$
     \cone(h)^q:= M^{q+1} \oplus N^{q}, \quad
     d_{\cone(h)}^q:= (-d_M^{q+1},h^{q+1}+d_N^{q}).
$$
We construct a morphism $\delta_{s,r}$
        satisfying (3-1) and (3-2)
          in a canonical way.
Take injective resolutions
       $\mu_{p^v}^{\otimes n} \ra I_v^{\bullet}$
          $(v=r,r+s)$ and an injective resolution
            $\mu_{p^s}^{\otimes n} \ra J_s^{\bullet}$
        in the category of sheaves on $V_{\et}$
          for which there is a short exact sequence of complexes
            of the form
$$
\begin{CD}
0 @>>> I_r^{\bullet} @>>> I_{r+s}^{\bullet} @>>> J_s^{\bullet} @>>> 0.
\end{CD}
$$
Let $a_v: \tau_{\leq n}j_*I_v^{\bullet} \ra
         \iota_*\nu_{Y,v}^{n-1} [-n]$ $(v=r, r+s)$
and
$b_s:\tau_{\leq n}j_*J_s^{\bullet} \ra
         \iota_*\nu_{Y,s}^{n-1} [-n]$
    be the natural maps of complexes
          that represent
       $\sigma_{X,v}(n):\tau_{\leq n}Rj_*\mu_{p^v}^{\otimes n}
             \ra \iota_*\nu_{Y,v}^{n-1} [-n]$
          with $v=r,r+s$ and $s$, respectively (cf.\ \S\ref{sect4.2}).
The complexes $\cone^{\bullet}(a_v)$ $(v=r,r+s)$ and
      $\cone^{\bullet}(b_s)$
       represent $\T_v(n)_X$ with
          $v=r,r+s$ and $s$, respectively.
We show that the sequence of complexes
\stepcounter{equation}
\begin{equation}\label{short}
\hspace{-30pt}
\begin{CD}
0  @>>> \cone^{\bullet}(a_r)
      @>>> \cone^{\bullet}(a_{r+s})
      @>{f}>> \cone^{\bullet}(b_s) @>>> 0
\end{CD}
\hspace{-30pt}
\end{equation}
   is exact.
Indeed, this exactness follows from that of the sequence
       $0 \ra \nu_{Y,r}^{n-1} \ra \nu_{Y,r+s}^{n-1}
           \ra \nu_{Y,s}^{n-1} \ra 0$
            (\cite{sato:ss}, 2.2.5 (2)) and
     that of the sequence of complexes
\begin{equation}\label{short2}
\hspace{-30pt}
\begin{CD}
0  @>>> \tau_{\leq n}\, j_*I_r^{\bullet}
          @>>> \tau_{\leq n} \, j_*I_{r+s}^{\bullet}
          @>>> \tau_{\leq n} \, j_*J_s^{\bullet}
          @>>> 0
\end{CD}
\hspace{-30pt}
\end{equation}
(cf.\ Theorem \ref{thm:hyodo} (1)).
Finally, we define $\delta_{s,r}$ as the composite
     $\cone^{\bullet}(b_s)
      \ra \cone^{\bullet}(f) \simeq \cone^{\bullet}(a_r)[1]$ in $D^b(X_{\et})$,
         i.e., connecting morphism associated with \eqref{short}.
By definition, $\delta_{s,r}$ is canonical and
      satisfies the properties (3-1) and (3-2).
This completes the proof.
\end{pf}
\addtocounter{thm}{2}
\begin{rem}\label{rem:flabby}
One can construct a map
       $\delta_{s,r}':\T_s(n)_X \ra \T_r(n)_X[1]$
         in $D^b(X_{\et},\Z/p^{s+r}\Z)$
        satisfying $(3$-$1)$ and $(3$-$2)$
         in the same way as above.
Clearly,
         $\delta_{s,r}'=\delta_{s,r}$
              in $D^b(X_{\et})$.
\end{rem}
\subsection{Gysin morphism and purity}\label{sect4.4}
We define Gysin morphisms
    for closed subschemes of $X$ contained in $Y$
     and prove {\bf T3} in \ref{thm0-2}.
See \S\ref{sect6} below for a purity result for horizontal subschemes.
\begin{lem}\label{lem:gysin2}
\begin{enumerate}
\item[(1)]
There is a unique morphism
$$
\begin{CD}
    g' : \nu_{Y,r}^{n-1}[-n-1]  @>>>
                       R\iota^!\T_r(n)_X
                \quad \hbox{ in } \; D^b(Y_{\et},\Z/p^r\Z)
\end{CD}
$$
fitting into a commutative diagram
      with distinguished rows
\stepcounter{equation}
\begin{equation}\label{DT:loc}
\hspace{-60pt}
{\small
\begin{CD}
\T_r(n)_X @>{t}>> \tau_{\leq n}Rj_*\mu_{p^r}^{\otimes n}
      @>{-\sigma_{X,r}(n)}>> \iota_* \nu_{Y,r}^{n-1}[-n] @>{g[1]}>>
         \T_r(n)_X[1]\\
     @| @VVV @V{R\iota_*(g')[1]}VV @| \\
\T_r(n)_X @>{j^*}>> Rj_*\mu_{p^r}^{\otimes n}
      @>{\delta^{\loc}_{V,Y}(\Z/p^r\Z(n)_X)}>>
         R\iota_*R\iota^!\T_r(n)_X[1]
      @>{\iota_*}>> \T_r(n)_X[1].
\end{CD}
}
\hspace{-60pt}
\end{equation}
Here $t$ and $g$ denote the same morphisms as in Lemma $\ref{lem:defdef}$,
   and the lower row is the localization distinguished triangle
          \eqref{DT:local}.
\item[(2)]
$g'$ induces an isomorphism
$$
\begin{CD}
\tau_{\leq n+1}(g'): \nu_{Y,r}^{n-1}[-n-1] @>{\simeq}>>
         \tau_{\leq n+1} R\iota^!\T_r(n)_X
                \quad \hbox{ in } \; D^b(Y_{\et},\Z/p^r\Z).
\end{CD}
$$
\end{enumerate}
\end{lem}
\begin{pf}
We first calculate the cohomology sheaves of
      $R\iota^!\T_r(n)_X$.
In the lower row of \eqref{DT:loc},
     the map of the $q$-th cohomology sheaves of $\alpha \circ t$
        is bijective (resp.\ injective) if
          $q<n$ (resp.\ $q=n$), by \eqref{defdef2}.
Hence by {\bf T2}, we obtain
\begin{equation}\label{p-sup}
\hspace{-30pt}
\begin{CD}
R^q\iota^!\T_r(n)_X \simeq
       \begin{cases}
          0         & (q < n+1),\\
          \iota^*R^{q-1}j_*\mu_{p^r}^{\otimes n}
           \qquad & (q>n+1),
       \end{cases}
\end{CD}
\hspace{-30pt}
\end{equation}
and a short exact sequence
\begin{equation}\label{p-sup'}
\hspace{-50pt}
\begin{CD}
0 @>>> FM_r^n
     @>>> M_r^n
     @>{\iota^*{\cal H}^n(\delta^{\loc}_{V,Y}(\T_r(n)_X))}>>
       R^{n+1}\iota^!\T_r(n)_X @>>> 0.
\end{CD}
\hspace{-50pt}
\end{equation}
By Lemma \ref{lem:CD}, \eqref{p-sup} and {\bf T2}, we have
$$
\begin{CD}
\Hom_{D^b(Y_{\et},\Z/p^r\Z)}
    (\T_r(n)_X[1], R\iota_*R\iota^!\T_r(n)_X[1])=0.
\end{CD}
$$
Hence the first assertion of the lemma follows from
      Lemma \ref{lem:CD2} (2).
The second assertion follows from \eqref{p-sup'}.
\end{pf}
\addtocounter{thm}{3}
\begin{defn}\label{def:gysin}
Let $\phi: Z \hra Y$ be a closed immersion of pure codimension.
Put $c:=\codim_{X}(Z)$, and
     let $i$ be the composite map $Z \hra Y \hra X$.
We define the morphism
\stepcounter{equation}
\begin{equation}\label{p-gysin}
\hspace{-50pt}
\begin{CD}
\gys_{i}^n:\nu_{Z,r}^{n-c}[-n-c]
      @>>>
       Ri^! \T_r(n)_X
       \quad \hbox{ in } \; D^b(Z_{\et},\Z/p^r\Z)
\end{CD}
\hspace{-50pt}
\end{equation}
    as follows, where $\nu_{Z,r}^{n-c}$ means the zero sheaf if $n < c$.
If $Z=Y$ $($hence $c=1$ and $i=\iota)$,
     then we define $\gys_{\iota}^n$ as the
      morphism $g'$ in Lemma $\ref{lem:gysin2}$.
This morphism agrees with
   the adjoint of $g$ in Lemma $\ref{lem:defdef}$
           by the commutativity of the right square of \eqref{DT:loc}.
For a general $Z$,
    we define $\gys_{i}^n$ as the composite
$$
\begin{CD}
\nu_{Z,r}^{n-c}[-n-c]
       @>{\gys_{\phi}^{n-1}[-n-1]}>>
     R\phi^!\nu_{Y,r}^{n-1} [-n-1]
       @>{R\phi^!(\gys_{\iota}^n)}>>
     R\phi^!R\iota^! \T_r(n)_X=Ri^! \T_r(n)_X.
\end{CD}
$$
See Definition $\ref{rem:nugysin}$ for $\gys_{\phi}^{n-1}$.
\end{defn}
\stepcounter{thm}
\begin{thm}[Purity]\label{thm:purity}
The morphism
$$
\begin{CD}
\tau_{\leq n+c}(\gys_{i}^n):
    \nu_{Z,r}^{n-c}[-n-c] @>>> \tau_{\leq n+c}Ri^! \T_r(n)_X
\end{CD}
$$
      is an isomorphism.
\end{thm}
\begin{pf}
By the definition of $\gys_{i}^n$,
   the morphism $\tau_{\leq n+c}(\gys_{i}^n)$ is decomposed as follows:
$$
\begin{CD}
\nu_{Z,r}^{n-c}[-n-c]
       @>{\tau_{\le n+c}(\gys_{\phi}^{n-1}[-n-1])}>>
        \tau_{\leq n+c}(R\phi^!\nu_{Y,r}^{n-1} [-n-1]) \\
       @>{\tau_{\leq n+c}R\phi^!\{\tau_{\leq n+1}(\gys_{\iota}^n)\}}>>
         \tau_{\leq n+c}\{R\phi^!(\tau_{\leq n+1}R\iota^! \T_r(n)_X)\}\\
       @>{{\mathrm{canonical}}}>>
         \tau_{\leq n+c}Ri^! \T_r(n)_X.
\end{CD}
$$
The first two arrows are isomorphisms by
      Theorem \ref{thm:wpurity} and Lemma \ref{lem:gysin2}.
We show that the last arrow is an isomorphism as well.
There is a distinguished triangle of the form
$$
     \tau_{\leq n+1}R\iota^!\T_r(n)_X
     \lra R\iota^!\T_r(n)_X
     \lra \tau_{\geq n+2}R\iota^!\T_r(n)_X
     \lra (\tau_{\leq n+1}R\iota^!\T_r(n)_X)[1]
$$
   and we have
     $\tau_{\geq n+2}R\iota^! \T_r(n)_X
        \simeq (\tau_{\geq n+1}\iota^*Rj_*\mu_{p^r}^{\otimes n})[-1]$
         (cf.\ \eqref{p-sup}).
Hence it suffices to show
\stepcounter{equation}
\begin{equation}\label{van:hagihara}
\tau_{\leq n+c-1}R\phi^!(\tau_{\geq n+1}\iota^*Rj_*\mu_{p^r}^{\otimes n})
     =0.
\end{equation}
By the exactness of \eqref{short2},
     there is a distinguished triangle of the form
$$
\begin{CD}
      \tau_{\geq n+1}Rj_*\mu_{p^{r-1}}^{\otimes n}
      \lra \tau_{\geq n+1}Rj_*\mu_{p^r}^{\otimes n}
      \lra \tau_{\geq n+1}Rj_*\mu_{p}^{\otimes n}
      \lra (\tau_{\geq n+1}Rj_*\mu_{p}^{\otimes n})[1]
\end{CD}
$$
Hence \eqref{van:hagihara} is reduced to the case $r=1$ and
  then to the following semi-purity due to
  Hagihara (cf.\ Theorem \ref{lem:fin} below):
$$
\begin{CD}
R^q\phi^!(\iota^*R^mj_*\mu_{p}^{\otimes n}) = 0 \;\; \hbox{ for any }
     \; m,\, q \;\hbox{ with } \;q \leq c-2,
\end{CD}
$$
where one must note $c= \codim_Y(Z) +1$.
This completes the proof.
\end{pf}
\stepcounter{thm}
\begin{cor}\label{cor:purity2}
Let $i: Z \ra X$ be a closed immersion of codimesion $\ge n+1$.
Then we have $R^qi^!\T_r(n)_X=0$ for any $q \le 2n+1$.
\end{cor}
\begin{pf}
If $Z[1/p]$ is empty,
     then we have $R^qi^!\T_r(n)_X=0$ for $q \le 2n+1$
        by Theorem \ref{thm:purity}.
We next prove the case that $Z[1/p]$ is non-empty.
Put $U :=Z[1/p]$ and $T := Z \sm U$.
Let $\alpha:T \hra Z$, $\beta : U \hra Z$ and $\gamma : T \hra X$
         be the natural immersions.
There is a long exact sequence of sheaves on $Z_{\et}$
$$
    \dotsb \lra
   \alpha_*R^q\gamma^!\T_r(n)_X \lra R^qi^!\T_r(n)_X  \lra
         R^q\beta_*\beta^*Ri^!\T_r(n)_X \lra \dotsb,
$$
where $\alpha_*R^q\gamma^!\T_r(n)_X$ is zero for $q \le 2n+1$
     by the previous case.
We show that $R^q\beta_*\beta^*Ri^!\T_r(n)_X$ is zero
    for $q \le 2n+1$.
Indeed, we have $\beta^*Ri^!\T_r(n)_X= R\psi^!\mu_{p^r}^{\otimes n}$
            with $\psi$ the closed immersion $U \hra V$,
    and it is concentrated in degrees $\ge 2n+2$
    by the absolute purity of Thomason and Gabber \cite{Th}, \cite{fujiwara}
    and by the assumption that $\codim_{X}(Z) \ge n+1$.
\end{pf}
We next prove a projection formula,
     which will be used later in \S\ref{sect5} and \S\ref{sect6}.
\begin{prop}[Projection formula]\label{prop:proj}
Let $i:Z \hra X$ be as in $\ref{def:gysin}$.
We define the morphism
$i^* : \T_r(n)_X  \to i_*\lam_{Z,r}^{n}[-n]$ in $D^b(X_{\et},\Z/p^r\Z)$
by the natural pull-back of symbols on the $n$-th
      cohomology sheaves $($cf.\ \eqref{defdef2}$)$.
Then the square
$$
\begin{CD}
       i_* \nu_{Z,r}^{m-c} [-m-c] \otimes^{\L} \T_r(n)_X
        @>{\gys^{m}_i \otimes^{\L} \id}>>
          \T_r(m)_X \otimes^{\L} \T_r(n)_X\\
     @V{(\sharp)}VV @VV{\eqref{def:prod2}}V\\
       i_*\nu_{Z,r}^{m+n-c}[-m-n-c]
           @>{\gys^{m+n}_i}>> \T_r(m+n)_X
\end{CD}
$$
commutes in $D^-(X_{\et},\Z/p^r\Z)$.
Here the left vertical arrow $(\sharp)$
is the composite map
$$
{\small
\begin{CD}
   i_*\nu_{Z,r}^{m-c} [-m-c] \otimes^{\L} \T_r(n)_X
@>{\id \otimes^{\L} i^*}>>
   i_*\nu_{Z,r}^{m-c} [-m-c] \otimes^{\L} i_*\lam_{Z,r}^{n}[-n]
    @.~ \lra ~@.
    i_*\nu_{Z,r}^{m+n-c}[-m-n-c],
\end{CD}
}
$$
and the last arrow is induced by
   the pairing $\eqref{def:prod}$ on the $(m+n+c)${\rm th} cohomology sheaves.
\end{prop}
\begin{pf}
One can easily check the case $Z=Y$ by the commutativity of
     the central square in \eqref{DT:loc}.
The general case is, by the previous case, reduced to
   the commutativity of a diagram
$$
\begin{CD}
   \phi_*\nu_{Z,r}^{m-c} [-m-c] \otimes^{\L} \lam_{Y,r}^n[-n]
   @>{\gys^{m-1}_{\phi}[-m-1] \otimes^{\L} \id[-n]}>>
         \nu_{Y,r}^{m-1}[-m-1]
            \otimes^{\L} \lam_{Y,r}^n[-n]\\
   @VVV @VVV\\
    \phi_*\nu_{Z,r}^{m+n-c}[-m-n-c]
   @>{\gys^{m+n-1}_{\phi}[-m-n-1]}>>
         \nu_{Y,r}^{m+n-1}[-m-n-1]
\end{CD}
$$
in $D^-(Y_{\et},\Z/p^r\Z)$ with $\phi:Z \hra Y$.
Here the vertical arrows are defined in a similar way
    as for $(\sharp)$.
We prove the commutativity of this square.
For two complexes
   $M^{\bullet}=(\{M^u \}_{u \in \Z}, \{ d_M^u:M^u \to M^{u+1}\}_{u \in \Z})$
   and $N^{\bullet}=(\{N^v \}_{v \in \Z}, \{ d_N^v \}_{v \in \Z})$,
    we define the double complex
     $M^{\bullet} \otimes N^{\bullet}$ as
$$
(M^{\bullet} \otimes N^{\bullet})^{u,v}:= M^{u} \otimes N^{v}, \quad
   \partial_1^{u,v}:=d_M^u \otimes \id_{N^v}, \quad
   \partial_2^{u,v}:=(-1)^u \, \id_{M^u} \otimes d_N^{v}.
$$
We write $(M^{\bullet} \otimes N^{\bullet})^{\tot}$
     for the associated total complex,
   whose image into the derived category gives
     $M^{\bullet} \otimes^{\L} N^{\bullet}$
   if either $M^{\bullet}$ or $N^{\bullet}$ is bounded above and
   consists of flat objects.
Now for $T \in \{Y,Z \}$ and
    $a \geq 0$, let $C^{\bullet}_r(T,a)$
   be the complex of sheaves defined in \S\ref{sect2.2}.
Because $\lam_{Y,r}^n$ is flat over $\Z/p^r\Z$
       by \cite{sato:ss}, 3.2.3,
   the commutativity in question follows from that of a diagram
   of complexes on $Y_{\et}$
$$
\begin{CD}
   (\phi_* C^{\bullet}_r(Z,m-c) [-m-c] \otimes
\lam_{Y,r}^n[-n])^{\tot}
       @. \;\; \hra \;\; @.
   (C^{\bullet}_r(Y,m-1)[-m-1] \otimes \lam_{Y,r}^n[-n])^{\tot}\\
   @V{\mathrm{product}}VV @. @VV{\mathrm{product}}V\\
   \phi_* C^{\bullet}_r(Z,m+n-c) [-m-n-c]
          @. \; \hra \; @.
   C^{\bullet}_r(Y,m+n-1)[-m-n-1],
\end{CD}
$$
where the vertical arrows are induced by the pairings \eqref{def:prod}
   and the horizontal arrows are natural inclusions of complexes.
This completes the proof.
\end{pf}
\subsection{Kummer sequence for $\Gm$ and purity of Brauer groups}\label{sect4.5}
We study the case $n=1$.
\begin{prop}\label{prop:modp}
Put $\Gm:=\O^{\times}_X$. Then
there is a unique morphism
$$
\begin{CD}
   \Gm \otimes^{\L}\Z/p^r\Z [-1]
     @>>> \T_r(1)_X \quad
     \hbox{ in } \; D^b(X_{\et},\Z/p^r\Z)
\end{CD}
$$
that extends the canonical isomorphism
    $j^*(\Gm \otimes^{\L}\Z/p^r\Z[-1]) \simeq \mu_{p^r}$.
Moreover it is an isomorphism, if $X$ satisfies $\ref{cond2'}$.
\end{prop}
\begin{pf}
Put ${\cal M}:= \Gm \otimes^{\L}\Z/p^r\Z[-1]$.
By definition,
    (i) ${\cal M}$ is concentrated in $[0,1]$, and
        (ii) there are natural isomorphisms
$$
\begin{CD}
{\cal H}^0({\cal M}) \simeq \ker(\Gm \os{\times p^r}{\lra} \Gm)
     \quad \hbox{ and } \quad
{\cal H}^1({\cal M}) \simeq \Gm/{p^r}.
\end{CD}
$$
Because
     $j^*{\cal M} \simeq \mu_{p^r}$ canonically in $D^b(V_{\et},\Z/p^r\Z)$,
       there is a natural morphism
${\cal M}  \lra
         \tau_{\leq 1}Rj_* \mu_{p^r}$
in $D^b(X_{\et},\Z/p^r\Z)$ by (i).
The composite morphism
$$
\begin{CD}
{\cal M}  @>>>
         \tau_{\leq 1}Rj_* \mu_{p^r}
           @>{\sigma_{X,r}(1)}>> \iota_*\nu_{Y,r}^0[-1]
\end{CD}
$$
is zero by (ii) and Lemma \ref{lem:CD}.
Hence by Lemma \ref{lem:CD2} (1),
    we obtain a unique morphism
     ${\cal M} \lra \T_r(1)_X$ that extends
      the isomorphism $j^*{\cal M} \simeq \mu_{p^r}$.
Next we prove that this morphism is
      bijective on cohomology sheaves,
        assuming that $X$ satisfies $\ref{cond2'}$.
By the standard purity for $\Gm$
   (\cite{g:brauer}, (6.3)--(6.5)),
       there is an exact sequence
$$
\begin{CD}
   0 @>>>
\Gm @>>> j_*j^* \Gm @>>>
     \bigoplus{}_{y \in Y^0}~ \iota_{y*}\Z
     @>>> 0,
\end{CD}
$$
where for $x \in X$,
    $\iota_x$ denotes the canonical map $x \hra X$.
Since $\bigoplus_{y \in Y^0}~ \iota_{y*}\Z$ is torsion-free,
   we have ${\cal H}^0({\cal M}) \simeq j_*\mu_{p^r}$
     and there is an exact sequence
$$
\begin{CD}
0 @>>> {\cal H}^1({\cal M}) @>>> R^1j_*\mu_{p^r}
      @>>> \bigoplus{}_{y \in Y^0}~ \iota_{y*}\Z/p^r\Z
         (=\iota_*\nu_{Y,r}^0)
\end{CD}
$$
by (ii) and the snake lemma.
Here we have used the isomorphism $(j_*j^* \Gm)/p^r \simeq R^1j_*\mu_{p^r}$
          obtained from Hilbert's theorem 90: $R^1j_*j^* \Gm=0$.
Now the assertion follows from \eqref{defdef2}.
\end{pf}
As an application of
     Corollary \ref{cor:purity2} and Proposition \ref{prop:modp},
         we prove the $p$-primary part of the purity of Brauer groups
          (cf.\ \cite{g:brauer}, \S6).
\begin{cor}[Purity of Brauer groups]\label{cor:purity}
Assume that $X$ satisfies $\ref{cond2'}$.
Let $i:Z \hra X$ be a closed immersion
         with $\codim_X(Z) \geq 2$.
Then the $p$-primary torsion part of $R^3i^!\Gm$
           is zero.
\end{cor}
\noindent
If $\dim(X) \le 3$, then the full sheaf $R^3i^!\Gm$ is zero
    by a theorem of Gabber \cite{gabber}.
\begin{pf}
By Proposition \ref{prop:modp}, there is a distinguished
      triangle
\addtocounter{equation}{2}
\begin{equation}\label{DT:weight1}
\hspace{-50pt}
\begin{CD}
\Gm[-1]  @>>> \T_r(1)_X @>>>
\Gm @>{\times p^r}>> \Gm
     \quad \hbox{ in } \; D^b(X_{\et}),
\end{CD}
\hspace{-50pt}
\end{equation}
which yields an exact sequence
    $R^3i^!\T_r(1)_X \to R^3i^!\Gm \os{\times p^r}{\to} R^3i^!\Gm$.
Hence the corollary follows from
    the vanishing result in Corollary \ref{cor:purity2}.
\end{pf}

\medskip
\section{\bf Cycle class and intersection property}\label{sect5}
\medskip
Throughout this section, we work with the setting in \S\ref{sect4.1} and
   assume that $X$ satisfies the condition \ref{cond2'}.
In this section we define the cycle class
  $\cl_X(Z) \in \H^{2n}_Z(X,\T_r(n))$ for an integral closed subscheme
    $Z \subset X$ of codmension $n \ge 0$,
     and prove an `intersection formula'
$$
   \cl_X(Z) \cup \cl_X(Z') = \cl_X(Z \cap Z')
   \quad \hbox{ in } \; \H^{2(m+n)}_Z(X,\T_r(m+n)),
$$
assuming that $Z$ of codimension $m$ and $Z'$ of codimension $n$
  are regular and meet transversally.
In \S\ref{sect6},
we will prove {\bf T4} in Theorem \ref{thm0-2} using this result.
\subsection{Cycle class}\label{sect5.1}
We first note a standard consequence of
   Corollary \ref{cor:purity2}.
\begin{lem}\label{vpurity}
Let $Z$ be a closed subscheme of $X$ of pure codimension $n \ge 0$.
Let $Z'$ be a dense open subset of $Z$, and let
  $T$ be the complement $Z \setminus Z'$.
Then the natural map
$$
   \H^{2n}_Z(X,\T_r(n)_X) \lra \H^{2n}_{Z'}(X \setminus T,\T_r(n)_{X 
\setminus T})
$$
is bijective.
\end{lem}
\begin{pf}
There is a long exact sequence of cohomology groups with supports
$$
\dotsb \ra \H^{2n}_{T}(X,\T_r(n)_X) \ra
  \H^{2n}_Z(X,\T_r(n)_X) \ra \H^{2n}_{Z'}(X \setminus T,\T_r(n)_X)
   \ra \H^{2n+1}_{T}(X,\T_r(n)_X) \ra \dotsb.
$$
Since $\codim_X(T) \ge n+1$, we have
   $\H^{2n}_{T}(X,\T_r(n)_X)=\H^{2n+1}_{T}(X,\T_r(n)_X)=0$
     by Corollary \ref{cor:purity2},
      which shows the lemma.
\end{pf}
\begin{defn}\label{def:cycleclass}
For an integral closed subscheme $Z \subset X$ of codimension $n \ge 0$,
  we define the cycle class
  $\cl_X(Z) \in \H^{2n}_Z(X,\T_r(n)_X)$ as follows.
\begin{enumerate}
\item[(1)]
If $Z$ is regular and contained in $Y$,
    then we define $\cl_X(Z)$ to be the image of $1 \in \Z/p^r\Z$
      under the Gysin map
$$
   \gys_{i}^n: \Z/p^r\Z \lra \H^{2n}_Z(X,\T_r(n)_X)
$$
induced by the Gysin morphism defined in Definition $\ref{def:gysin}$.
\item[(2)]
If $Z$ is regular and {\it not} contained in $Y$,
  then we have Gabber's refined cycle class
    $\cl_{V}(U) \in \H^{2n}_U(V,\mu_{p^r}^{\otimes n})$
     $($cf.\ \cite{fujiwara}$)$,
     where we put $U:=Z[1/p]$ and $V:=X[1/p]$.
We define $\cl_X(Z)$ as the inverse image of $\cl_{V}(U)$ under the natural map
$$
   \H^{2n}_Z(X,\T_r(n)_X) \lra \H^{2n}_{U}(V,\mu_{p^r}^{\otimes n}).
$$
This map is bijective by Lemma $\ref{vpurity}$ and excision,
   and hence $\cl_X(Z)$ is well-defined.
Note that Gabber's refined cycle class agrees with
   Deligne's cycle class $($\cite{sga4}, Cycle$)$ in
    any situation where the latter is defined $($cf.\ \cite{fujiwara}, $1.1.5)$.
\item[(3)]
For a general $Z$,
  we take a dense open regular subset $Z' \subset Z$
      and define $\cl_X(Z)$ to be the inverse image of
        $\cl_{X'}(Z') \in \H^{2n}_{Z'}(X',\T_r(n)_X)$
     $(X':=X \setminus (Z \setminus Z'))$
      under the natural map
$$
   \H^{2n}_Z(X,\T_r(n)_X) \lra \H^{2n}_{Z'}(X',\T_r(n)_{X'}),
$$
which is bijective by Lemma $\ref{vpurity}$ and $\cl_X(Z)$ is well-defined.
\end{enumerate}
\end{defn}
We prove the following result:
\begin{prop}[Intersection property]\label{prop:int}
Let $Z$ and $Z'$ be integral regular closed subschemes of $X$
  of codimension $a$ and $b$, respectively.
Assume that $Z$ and $Z'$ meet transversally on $X$.
Then we have
$$
   \cl_X(Z) \cup \cl_X(Z') = \cl_X(Z \cap Z')
    \quad \hbox{ in } \; \H^{2(a+b)}_{Z \cap Z'}(X,\T_r(a+b)_X).
$$
Here, if $Z \cap Z'$ is not connected,
   then $\cl_X(Z \cap Z')$ means the sum of the cycle classes of
     the connected components.
\end{prop}
\smallskip
\subsection{Proof of Proposition \ref{prop:int}}\label{sect5.2}
Without loss of generality, we may assume that
   $Z \cap Z'$ is connected (hence integral and regular).
If $Z \cap Z'$ is not contained in $Y$,
   the assertion follows from Lemma \ref{vpurity} and
    the corresponding property of
     Gabber's refined cycle classes \cite{fujiwara}, 1.1.4.
We treat the case that $Z \cap Z' \subset Y$.
Let $x$ be the generic point of $Z \cap Z'$.
By Lemma \ref{vpurity}, we may replace $X$ by $\Spec(\O_{X,x})$.
Because $Z$ and $Z'$ are regular and meet transversally at $x$,
   there is a normal crossing divisor $D=\cup_{i=1}^{a+b} \, D_i$
   with each $D_i$ integral regular
    such that $\cap_{i=1}^a D_i = Z$ and
       $\cap_{i=a+1}^{a+b} D_i = Z'$.
Therefore we are reduced to the following local assertion:
\begin{lem}\label{class}
Suppose that $X$ is local with closed point $x$ of characteristic $p$.
Put $n:=\codim_X(x) \ge 1$.
Let $D=\cup_{i=1}^n \, D_i$
   be a normal crossing divisor on $X$ with each $D_i$ integral regular
    such that $\cap_{i=1}^n \, D_i = x$.
Then the cohomology class
$$
   \cl_X(x;D):=\cl_X(D_1) \cup \cl_X(D_2) \cup \dotsb
     \cup \cl_X(D_n) \in \H^{2n}_x(X,\T_r(n)_X)
$$
depends only on the flag\,$:$
$D_1 \supset D_1 \cap D_2 \supset \dotsb \supset
     D_1 \cap \dotsb \cap D_{n-1} \supset x$,
and agrees with $\cl_X(x)$.
\end{lem}
\noindent
We prove this lemma by induction on $n \ge 1$.
The case $n=1$ is clear.
Suppose that $n \ge 2$ and put $S:=\cap_{i=1}^{n-1} D_i$.
Let $\psi$ (resp.\ $i_x$) be the closed immersion $S \hra X$
   (resp.\ $x \to S$).
Note that $S$ is regular, local and of dimension $1$.

We first show the case that $S \subset Y$.
By the induction hypothesis and Lemma \ref{vpurity},
  we have $\cl_X(D_1) \cup \dotsb \cup \cl_X(D_{n-1})= \cl_X(S)$,
   and hence
$$
\cl_X(x;D) =\cl_X(S)\cup\cl_X(D_n)
     =\gys_{\psi}^n(\gys_{i_x}^{n-c}(1))=\cl_X(x).
$$
Here the second equality follows from
  Proposition \ref{prop:proj} for $\psi$ and the last equality
    follows from Remark \ref{rem:newtrace} (1).
In particular, $\cl(x;D)$ depends only on the flag of $D$.

We next show the case that $S \not\subset Y$.
Let $y$ be the generic point of $S$.
Since $\ch(y)=0$,
  we have $\cl_X(D_1) \cup \dotsb \cup \cl_X(D_{n-1})= \cl_X(S)$
    by Lemma \ref{vpurity} and \cite{fujiwara}, 1.1.4.
We have to show
\begin{sublem}\label{lem:int}
Let $E$ and $E'$ be regular connected divisors on $X$
   each of which meets $S$ transversally at $x$.
Then we have $\cl_X(S) \cap \cl_X(E)=\cl_X(S) \cap \cl_X(E')$.
\end{sublem}
\noindent
We first finish the proof of the lemma, admitting this sublemma.
It implies that
   $\cl_X(x;D)$ depends only on the flag of $D$,
  and moreover that $\cl_X(x;D)$ is independent of $D$
    by \cite{sga4.5}, Cycle, 2.2.3.
Hence we obtain $\cl_X(x;D)=\cl_X(x)$ by the computation
   in the previous case.
\begin{pf*}{\it Proof of Sublemma \ref{lem:int}.}
Let $E$ be a regular divisor on $X$ as in the sublemma.
The map
$$
\H^{2}_E(X,\T_r(1)_X) \lra
  \H^{2n}_x(X,\T_r(n)_X), \qquad \alpha \mapsto \cl_X(S) \cup \alpha
$$
factors through a natural pull-back map
$$
\psi^* : \H^{2}_E(X,\T_r(1)_X) \lra \H^2_x(S,\psi^*\T_r(1)_X).
$$
We compute $\psi^*(\cl_X(E))$ as follows.
Since $\ch(y)=0$, we have
   $(\psi^*\T_r(1)_X) \vert _y \simeq \mu_{p^r}$ on $y_{\et}$ and
  there is a commutative diagram with exact rows
$$
\begin{CD}
\H^1(X,\T_r(1)_X)
  @. ~ \lra ~ @.
\H^1(X \setminus E,\T_r(1)_{X \setminus E})
  @. ~ \os{\delta_1}{\lra} ~ @.
\H^{2}_E(X,\T_r(1)_X) \\
@V{\psi^*}VV @. @V{\psi^*}VV @. @V{\psi^*}VV \\
\H^1(S,\psi^*\T_r(1)_X)
  @. ~ \lra ~ @.
\H^1(y,\mu_{p^r})
  @. ~ \os{\delta_2}{\lra} ~ @.
\H^2_x(S,\psi^*\T_r(1)_X),
\end{CD}
$$
where $\delta_1$ denotes $\delta_{X\setminus E, E}^{\loc}(\psi^*\T_r(1)_X)$
  and $\delta_2$ denotes $\delta_{y,x}^{\loc}(\psi^*\T_r(1)_X)$
   (cf.\ \eqref{DT:local}).
Take a prime element $\pi_E \in \vG(X,\O_X)$
   which defines $E$, and let $\{ \pi_E \}
    \in \H^1(X \setminus E,\T_r(1)_{X \setminus E})$ be
   the image of $\pi_E$ under the boundary map of Kummer theory
     (cf.\ \eqref{DT:weight1})
$$
   \vG(X \setminus E,\O_{X\setminus E}^{\times})
      \lra \H^1(X \setminus E,\T_r(1)_{X \setminus E}).
$$
We have $\cl_X(E)=-\delta_1(\{ \pi_E^{-1} \})$
   by \cite{sga4.5}, Cycle, 2.1.3 (cf.\ \eqref{sign}).
By the diagram,
$$
\psi^*(\cl_X(E))=-\delta_2(\psi^*\{ \pi_E\}))
   = -\delta_2(\{ \ol {\pi_E}\}).
$$
Here $\ol {\pi_E}$ denotes the residue class of $\pi_E$
  in $\O_{S,x}$ and it is a prime element by the assumption that
    $E$ meets $S$ transversally at $x$.
Moreover we have $\delta_2(\{u\})=0$ for any unit $u \in \O_{S,x}^{\times}$,
   because every $u \in \O_{S,x}^{\times}$ lifts to
      $\O_{X,x}^{\times}$.
Hence for a fixed prime $\pi_x \in \O_{S,x}$
  we have $\psi^*(\cl_X(E))=-\delta_2(\{\pi_x\})$,
     which shows the sublemma.
\end{pf*}
%
This completes the proof of Lemma \ref{class} and
   Proposition \ref{prop:int}.

\medskip
\section{\bf Compatibility and purity for horizontal subschemes}\label{sect6}
\medskip
In this section,
   we prove {\bf T4} in Theorem \ref{thm0-2}.
This result is rather technical, but
   we will need its consequence, Theorem \ref{thm:gysin},
    to prove the covariant functoriality {\bf T7} in \S\ref{sect7}.
\subsection{Gysin maps}\label{sect6.1}
We work with the setting in \S\ref{sect4.1},
     and assume that $X$ satisfies $\ref{cond2}$.
Let $b$ and $n$ be integers with $n \geq b \geq 0$.
For $x \in X^b$,
      we define the complex $\Z/p^r\Z(n)_x$ on $x_{\et}$ as
$$
\begin{CD}
\Z/p^r\Z(n)_x :=
        \begin{cases}
            \mu_{p^r}^{\otimes n}
               \quad & (\hbox{if }\, \ch(x) \not=p), \\
            \logwitt x r {n}[-n]
               \quad & (\hbox{if }\, \ch(x) =p). \\
        \end{cases}
\end{CD}
$$
We define the Gysin map
$$
\gys_{i_x}^{n}:
     \H^{n-b}(x,\Z/p^r\Z(n-b)_x) \lra \H^{n+b}_x(X,\T_r(n)_X)
      := \H^{n+b}_x(\Spec(\O_{X,x}),\T_r(n)_X)
$$
as the map induced by
         the Gysin morphism for $i_x:x \hra X$,
           if $\ch(x)=p$ (cf.\ Definition \ref{def:gysin}).
If $\ch(x) \not=p$, we define $\gys_{i_x}^{n}$
by sending $\alpha \in \H^{n-b}(x,\mu_{p^r}^{\otimes n-b})$
    to $\cl_V(x) \cup \alpha$,
            where $\cl_V(x) \in \H^{2b}_x(V,\mu_{p^r}^{\otimes b})$
          denotes Gabber's refined cycle class
        we mentioned in Definition \ref{def:cycleclass} (2).
The aim of this section is to prove the following two theorems:
\begin{thm}[{Compatibility}]\label{thm:comp}
Let $x$ and $y$ be points
    with $x \in \ol{ \{ y \} } \cap Y \cap X^b$
     and $y \in X^{b-1}$
      $($hence $\ch(y)=0$ or $p)$.
Then the diagram
\stepcounter{equation}
\begin{equation}\label{p-diagram}
\hspace{-50pt}
\begin{CD}
   \H^{n-b+1}(y,\Z/p^r\Z(n-b+1)_y)
   @>{-\partial^{\val}_{y,x}}>> \H^{n-b}(x,\Z/p^r\Z(n-b)_x)\\
@V{\gys_{i_y}^n}VV @VV{\gys_{i_x}^n}V \\
\H^{n+b-1}_y(X,\T_r(n)_X)
    @>{\delta^{\loc}_{y,x}(\T_r(n)_X)}>>
   \H^{n+b}_x(X,\T_r(n)_X)
\end{CD}
\hspace{-50pt}
\end{equation}
is commutative
$($see \eqref{note.3.3} for the definition of the bottom arrow$)$.
\end{thm}
\noindent
Sheafifying this commutative diagram,
    we obtain {\bf T4} in Theorem \ref{thm0-2}.
As for the case $x \not \in Y$,
   the corresponding commutativity is proved in \cite{jss}, \S1.
\stepcounter{thm}
\begin{thm}[Purity]\label{thm:gysin}
Let $Z$ be an integral locally closed subscheme of $X$
      which is flat over $B$ and satisfies $\ref{cond2}$.
Put $c:=\codim_X(Z)$ and
     $U:=Z[1/p]$.
Let $i$ and $\psi$ be the locally closed immersions
    $Z \hra X$ and $U \hra V$, respectively.
Then for $n \ge c$, there is a unique morphism
\begin{equation}\notag
\begin{CD}
\hspace{-50pt}
\gys_{i}^n:
    \T_r(n-c)_Z[-2c] @>>> Ri^!\T_r(n)_X
       \quad \hbox{ in } \; D^b(Z_{\et},\Z/p^r\Z)
\hspace{-50pt}
\end{CD}
\end{equation}
that extends the purity isomorphism
     $($cf.\ \cite{Th}, \cite{fujiwara}$)$
\begin{equation}\notag
\begin{CD}
\gys_{\psi}^n:
      \mu_{p^r}^{\otimes n-c}[-2c]
    @>{\simeq}>>   R\psi^!\mu_{p^r}^{\otimes n}
             \quad \hbox{ in } \; D^b(U_{\et},\Z/p^r\Z).
\end{CD}
\end{equation}
Moreover, $\gys_{i}^n$ induces an isomorphism
$$
\begin{CD}
    \tau_{\le n+c}(\gys_i^n)
     :  \T_r(n-c)_Z[-2c] @>{\simeq}>> \tau_{\le n+c} Ri^!\T_r(n)_X.
\end{CD}
$$
\end{thm}
\noindent
This result extends Theorem \ref{thm:purity}
     to horizontal situations.
Before starting the proof of these theorems,
     we state a consequence of Theorem \ref{thm:comp}.
For a point $x \in X^n$ and a closed subscheme $S \subset X$
           containing $x$
     there is a natural map
$$
\begin{CD}
     \H^{2n}_x(X,\T_r(n)_X) @>>> \H^{2n}_S(X,\T_r(n)_X)
\end{CD}
$$
by Lemma \ref{vpurity}.
By Theorem \ref{thm:comp} and \cite{jss}, Theorem 1.1,
       we obtain
\begin{cor}[Reciprocity law]
\label{cor:recip}
Let $y$ be a point with $y \in X^{n-1}$, and
    put $S:= \ol{\{ y \}} \subset X$.
Then for any $\alpha \in \H^1(y,\Z/p^r\Z(1)_y)$,
      we have
$$
\begin{CD}
     \sum{}_{x \in X^n \cap S}~
        \gys_{i_x}^n(\partial_{y,x}^{\val}(\alpha)) =0
        \quad \hbox{ in }\; \H^{2n}_S(X,\T_r(n)_X).
\end{CD}
$$
Consequently, the sum of Gysin maps
$$
\begin{CD}
    \sum_{x \in X^n}~\gys_{i_x}^n :
     \bigoplus_{x \in X^n}~ \Z/p^r\Z  @>>> \H^{2n}(X,\T_r(n)_X)
\end{CD}
$$
factors through the Chow group of algebraic cycles modulo rational
equivalence$:$
$$
\begin{CD}
     \cl_{X,r}^n : \CH^n(X)/p^r  @>>>
\H^{2n}(X,\T_r(n)_X).
\end{CD}
$$
\end{cor}
\begin{rem}\label{rem:gysin0}
\begin{enumerate}
\item[(1)]
The case $\ch(y)=p$
     of Theorem $\ref{thm:comp}$ follows from the definition of Gysin maps
       $($cf.\ $\ref{rem:nugysin}$, $\ref{def:gysin})$
      and a similar arguments as for \cite{sato:ss}, $2.3.1$
        $($see also \eqref{sign}$)$.
On the other hand,
   the case $\ch(y)=0$ of Theorem $\ref{thm:comp}$
     is closely related to Theorem $\ref{thm:gysin}$.
\item[(2)]
Corollary $\ref{cor:recip}$ is not a new result if $X$ is smooth over $B$.
In fact, by an argument of Geisser $\cite{geisser}$, $\S6$, Proof of $1.3$,
     one can construct a canonical map
       from higher Chow groups of $X$
          to $\H^*(X,\T_r(n)_X)$.
A key ingredient in his argument is
      the localization exact sequences for higher Chow groups
          due to Levine \cite{levine:mot}.
In this paper, we give a more elementary proof of Theorem $\ref{thm:comp}$
          without using Levine's localization sequences.
\end{enumerate}
\end{rem}
\noindent
In what follows, we refer the case $\ch(y)=0$
     of Theorem \ref{thm:comp} as Case (M).
We will proceed the proof of Theorems \ref{thm:comp}
      and \ref{thm:gysin} in three steps.
In \S\ref{sect6.2},
      we will prove Case (M) of Theorem \ref{thm:comp}
         assuming that $X$ satisfies \ref{cond2'}
           and that $S:=\ol {\{ y \}}$ is normal at $x$.
In \S\ref{sect6.3},
      we will prove Theorem \ref{thm:gysin}
       assuming that $X$ satisfies \ref{cond2'}
        and then reduce Case (M) of Theorem \ref{thm:comp}
            to the case where $X$ is {\it smooth} over $B$
            (and $S$ is arbitrary).
The last case will be proved in \S\ref{sect6.4},
    which will complete the proof of Theorems \ref{thm:comp} and 
\ref{thm:gysin}.
\subsection{Proof of the theorems, Step 1}\label{sect6.2}
In this step, we prove Case (M) of Theorem \ref{thm:comp}
     assuming that $X$
and $S(=\ol{ \{ y \} })$
         are regular at $x$.
Replacing $X$ by $\Spec(\O_{X,x}^{\h})$
     and replacing $y$ by the point on $\Spec(\O_{X,x}^{\h})$
         lying above $y$,
      we suppose that $X$ is regular henselian local with closed point $x$.
Note that it suffices to show the desired compatibility
     in this situation, and that
     $\O_{S,x}$ is a henselian discrete valuation ring.
By the Bloch-Kato theorem \cite{bk}, (5.12),
      $\H^{n-b+1}(y,\mu_{p^r}^{\otimes n-b+1})$
         is generated by symbols of the forms
$$
    {\mathrm{(i)}} \quad \{ \beta_1,\dotsc, \beta_{n-b+1}\}
     \qquad \hbox{ and } \qquad
   {\mathrm{(ii)}} \quad \{\pi_x, \beta_{1}, \dotsc, \beta_{n-b}\},
$$
where each $\beta_{\lam}$ belongs to $\O_{S,x}^{\times}$, and
   $\pi_x$ denotes a prime element of $\O_{S,x}$.
We show that
      the diagram \eqref{p-diagram}
        commutes for these two kinds of symbols.
Recall that $\gys_{i_y}^n$
             is given by the cup product with the cycle class
$\cl_V(y) \in \H^{2b-2}_y(V,\mu_{p^r}^{\otimes b-1})$,
    and that this cycle class extends to the cycle class
$\cl_{X}(S) \in \H^{2b-2}_S(X,\T_r(b-1)_X)$
(cf.\ Definition \ref{def:cycleclass}).
We first show that the diagram \eqref{p-diagram}
      commutes for symbols of the form (i).
Because a symbol $\omega$ of this form lifts to
$\wt \omega \in \H^{n-b+1}(X, \T_r(n-b+1)_X)$,
    its image $\gys_{i_y}^n(\omega)$ lifts to
      $\cl_{X}(S) \cup \wt{\omega} \in \H^{n+b-1}_S(X, \T_r(n)_X)$.
Hence we have $\delta_{y,x}^{\loc} \circ \gys_{i_y}^n(\omega)=0$,
         which implies the assertion.
We next consider symbols of the form (ii).
The map $\delta_{y,x}^{\loc}(\T_r(n)_X) \circ \gys_{i_y}^n$
      sends a symbol $\{\pi_x,\beta_{1}, \dotsc, \beta_{n-b}\}$ to
\begin{equation}\label{gysin:hor7}
\hspace{-50pt}
\begin{CD}
      \delta (\cl_V(y) \cup \{ \pi_x \})
            \cup \omega
              \in \H^{n+b}_x(X,\T_r(n)_X).
\end{CD}
\hspace{-50pt}
\end{equation}
Here $\omega$ denotes a lift of $\{\beta_{1}, \dotsc, \beta_{n-b}\}$
      to $\H^{n-b}(X,\T_r(n-b)_X)$, and
       $\delta$ denotes the connecting map
    $\delta_{y,x}^{\loc}(\T_r(b)_X)$.
Since $\cl_V(y)$ extends to $\cl_X(S)$,
    we have $\delta (\cl_V(y) \cup \{ \pi_x \})=-\cl_X(x)$
      by Proposition \ref{prop:int}
     (see also the proof of Sublemma \ref{lem:int}).
Hence we have
$$
      \eqref{gysin:hor7} =
         - \cl_{X}(x) \cup \omega
         = - \gys_{i_x}^n(\ol{\omega})
         = - \gys_{i_x}^n \circ \partial_{y,x}^{\val}
          (\{\pi_x, \beta_{1}, \dotsc, \beta_{n-b}\}),
$$
    where $\ol{\omega}$ denotes
        the residue class of $\omega$
           in $\H^{0}(x,\logwitt x r {n-b})$
        and the second equality follows from
           Proposition \ref{prop:proj} for $i_x$.
Thus we obtain the desired commutativity.
This completes Step 1.
\subsection{Proof of the theorems, Step 2}\label{sect6.3}
In this step we prove Theorem \ref{thm:gysin}
        assuming that $X$ satisfies \ref{cond2'}
           (see also Remark \ref{rem:gysin} below).
Let $i : Z \hra X$ and $\psi:U \hra V$
     be as in Theorem \ref{thm:gysin}.
Let $T$ be the divisor on $Z$
       defined by the radical of $(p) \subset \O_Z$.
We obtain a commutative diagram of schemes
$$
\begin{CD}
     U  @>{\beta}>> Z @<{\alpha}<< T \\
    @V{\psi}V{\quad\;\;\;\square}V
       @V{i}VV  @VVV\\
     V   @>{j}>>  X @<{\iota}<<  Y.
\end{CD}
$$
Put $\phi:= i \circ \alpha: T \hra X$,
     ${\frak L}:= \T_r(n-c)_Z[-2c]$
         and ${\frak M}:= Ri^!\T_r(n)_X$.
By Theorem \ref{thm:purity} for $\phi$,
      there is an isomorphism
\begin{equation}\label{gysin:hor3}
\nu_{T,r}^{n-c-1}[-n-c-1]
   \simeq \tau_{\leq n+c+1} R\phi^!\T_r(n)_X
       = \tau_{\leq n+c+1}R\alpha^!{\frak M}.
\end{equation}
Consider a diagram
   with distinguished rows in $D^+(Z_{\et},\Z/p^r\Z)$
\begin{equation}\label{gysin:hor4}
\hspace{-90pt}
\begin{CD}
    {\frak L} @>{t[-2c]}>>
   (\tau_{\leq n-c}R\beta_*\mu_{p^r}^{\otimes n-c})[-2c]
   @>{-\sigma[-2c]}>> \alpha_*\nu_{T,r}^{n-c-1}[-n-c]
   @>{g[-2c]}>> {\frak L}[1]\\
     @.
    @V{R\beta_*(\gys_{\psi}^n)}VV
     @VV{R\alpha_*(\gys_{\phi}^{n})[1]}V @.\\
   {\frak M}
   @>{\beta^*}>> R\beta_*\beta^*{\frak M}
   @>{\delta^{\loc}_{U,T}({\frak M})}>>
   R\alpha_*R\alpha^!{\frak M}[1] @>{\alpha_*}>> {\frak M}[1].
\end{CD}
\hspace{-70pt}
\end{equation}
Here the upper low is the distinguished
    triangle defining $\T_r(n-c)_Z$
     shifted by degree $-2c$, and
    we wrote $\sigma$ for $\sigma_{Z,r}(n-c)$.
The lower row is the localization distinguished triangle for ${\frak M}$
   (cf.\ \eqref{DT:local}).
We show that the square of \eqref{gysin:hor4} commutes.
Indeed, by \eqref{gysin:hor3} and Lemma \ref{lem:CD},
      it is enough to show that the induced diagram of
        the $(n+c)$-th cohomology sheaves is commutative
           at the generic points of $T$, which was shown in Step 1.
Hence the square commutes, and
   there is a unique morphism
       $\gys_i^n : {\frak L} \to {\frak M}$
    that extends $\gys_{\psi}^n$ by Lemma \ref{lem:CD2} (1),
     because
$$
\begin{CD}
     \Hom_{D^+(Z_{\et},\Z/p^r\Z)}
        ({\frak L}, R\alpha_*R\alpha^!{\frak M})=0
\end{CD}
$$
by \eqref{gysin:hor3} and Lemma \ref{lem:CD}.

We next show that $\tau_{\le n+c}(\gys_i^n)$ is an isomorphism.
By the commutativity of the square of \eqref{gysin:hor4},
      the morphism $\delta^{\loc}_{U,T}({\frak M})$
   is surjective on the $(n+c)$-th cohomology sheaves, and
   there is a distinguished triangle
\begin{equation}\notag
\begin{CD}
    \tau_{\leq n+c}{\frak M}
       @>{\beta^*}>> \tau_{\leq n+c}R\beta_*\beta^*{\frak M}
        @>{\delta^{\loc}_{U,T}({\frak M})}>>
             (\tau_{\leq n+c+1}R\alpha_*R\alpha^!{\frak M})[1]
              @>{(\alpha_*)'}>>
                (\tau_{\leq n+c}{\frak M})[1],
\end{CD}
\end{equation}
where the arrow $(\alpha_*)'$ is obtained by decomposing
    $\alpha_*: R\alpha_*R\alpha^!{\frak M})[1] \to {\frak M}[1]$.
Replacing the lower row of \eqref{gysin:hor4}
        with this distinguished triangle,
    we see that $\tau_{\le n+c}(\gys_i^n)$ is an isomorphism.
This completes Step 2.
\par
\medskip
\addtocounter{thm}{2}
\begin{cor}\label{cor:trans}
Let $i: Z \hra X$ be as in Theorem $\ref{thm:gysin}$,
    and assume further that
        $X$ satisfies $\ref{cond2'}$.
Let $h : Z' \hra Z$ be a closed immersion
         of pure codimension with $\ch(Z')=p$.
Put $g:=i \circ h$ and $c':=\codim_X(Z')$.
Then we have
     $\gys_g^n=Rh^!(\gys_i^n) \circ (\gys_h^{n-c}[-2c])$
      as morphisms
$$
\begin{CD}
\nu_{Z'}^{n-c'}[-n-c'] @>>>
       Rg^!\T_r(n)_X \quad \hbox{ in } \; D^b(Z'_{\et},\Z/p^r\Z).
\end{CD}
$$
\end{cor}
\begin{pf}
Because $\tau_{\le n+c'-1}Rg^!\T_r(n)_X = 0$
     by Theorem \ref{thm:purity},
a morphism $\nu_{Z'}^{n-c'}[-n-c'] \to Rg^!\T_r(n)_X$
     is determined by a map on the $(n+c')$-th cohomology sheaves
      (cf.\ Lemma \ref{lem:CD}).
Because $R^{n+c'}g^!\T_r(n)_X$ is isomorphic to
     the sheaf $\nu_{Z'}^{n-c'}$ by Theorem \ref{thm:purity},
    we are reduced to
      the case that $X$ and $Z$ are local with closed point $Z'$,
    and moreover,
    to the case that $Z'$ is a generic point of $Z \otimes_{\Z} {\Bbb F}_p$
          (that is, $c'=c+1$).
This last case follows from
    the commutativity of \eqref{p-diagram} proved in Step 1.
\end{pf}
\begin{rem}\label{rem:gysin}
\begin{enumerate}
\item[(1)]
By the results in this step and the bijectivity of $\gys_{i_x}^n$
       in \eqref{p-diagram}
         $($cf.\ Theorem $\ref{thm:purity})$,
       Case $\mathrm{(M)}$ of Theorem $\ref{thm:comp}$
         $($with $X$ and $\ol{ \{ y \} }$ arbitrary$)$
           is reduced to the case where $X$ is smooth over $B$.
We will prove this case in the next step.
\item[(2)]
Once we finish the proof of Theorem $\ref{thm:comp}$,
       we will obtain Theorem $\ref{thm:gysin}$
         by repeating the same arguments as for Step $2$.
\end{enumerate}
\end{rem}
\subsection{Proof of the theorems, Step 3}\label{sect6.4}
Assume that $X$ is smooth over $B$.
In this step, we prove Case (M) of Theorem \ref{thm:comp} for this $X$,
    which will complete the proof of Theorems \ref{thm:comp} and \ref{thm:gysin}
              (cf.\ Remark \ref{rem:gysin}).
We first show Lemma \ref{lem:trace2} below.
Let $m$ be a positive integer, and
    put $\P:=\P_X^m$.
Consider cartesian squares of schemes
$$
\begin{CD}
     \P_V  @>{\beta}>> \P  @<{\alpha}<< \P_Y \\
     @V{\psi}V{\quad\;\;\; \square}V
      @V{f}V{\quad\;\;\; \square}V  @V{g}VV\\
      V @>{j}>> X @<{\iota}<< Y.
\end{CD}
$$
For this diagram, we prove
\begin{lem}\label{lem:trace2}
\begin{enumerate}
\item[(1)]
There is a unique morphism
$$
\begin{CD}
      \tr_f^{n} : Rf_*\T_r(n+m)_{\P}[2m]
                     @>>> \T_r(n)_X \quad \hbox { in } \;
                       D^b(X_{\et},\Z/p^r\Z)
\end{CD}
$$
that extends the trace morphism for $\psi$
       $($\cite{sga4}, $\mathrm{XVIII}.2.9$, $\mathrm{XII}.5.3)$
$$
\begin{CD}
      \tr_{\psi}^{n} : R\psi_{*}\mu_{p^r}^{\otimes n+m}[2m]
                     @>>> \mu_{p^r}^{\otimes n} \quad \hbox { in } \;
                       D^b(V_{\et},\Z/p^r\Z).
\end{CD}
$$
\item[(2)] $\tr_f^n$ fits into a commutative diagram
$$
\begin{CD}
    R\iota_*Rg_*\logwitt {\P_Y} r {n+m-1} [m-n-1]
        @>{(\natural)}>> Rf_*\T_r(n+m)_{\P}[2m] \\
     @V{R\iota_*(\tr_g^{n-1})[-n-1]}VV  @V{\tr_f^n}VV\\
    \iota_*\logwitt Y r {n-1} [-n-1] @>{\gys_{\iota}^{n}}>> \T_r(n)_X.
\end{CD}
$$
Here $\tr_g^{n-1}$ denotes
  $(-1)^m$-times of the Gysin morphism of Gros
       $($\cite{gros:purity}, $\mathrm{II}.1.2.7)$.
The arrow $(\natural)$ is induced by
    the isomorphism $R\iota_*Rg_*=Rf_*R\alpha_*$ and
     the Gysin morphism $\gys_{\alpha}^{m+n}$.
\end{enumerate}
\end{lem}
\begin{pf}
Because $Rf_*\T_r(n+m)_{\P}[2m]$
     is concentrated in derees $\le n$
       (\cite{sga4}, XII.5.2, X.5.2),
      it is enough to show that the square
$$
\begin{CD}
Rf_*(\tau_{\leq n+m }R\beta_*\mu_{p^r}^{\otimes n+m})[2m]
       @>{Rf_*(\sigma_{\P,r}(n+m))[2m]}>>
           Rf_*(\alpha_*\logwitt {\P_Y} r {n+m-1})[m-n]\\
     @VVV @VV{R\iota_*(\tr_g)[-n]}V \\
\tau_{\leq n} Rj_*\mu_{p^r}^{\otimes n}
       @>{\sigma_{X,r}(n)}>> \iota_*\logwitt {Y} r {n-1}[-n]
\end{CD}
$$
is commutative in $D^b(X_{\et},\Z/p^r\Z)$
        (cf.\ Lemma \ref{lem:CD2} (1)).
Here the left vertical arrow is defined as the composite
      of the natural morphism
$$
Rf_*(\tau_{\leq n+m }R\beta_*\mu_{p^r}^{\otimes n+m})[2m]
        \lra \tau_{\leq n}
           (Rf_*R\beta_*\mu_{p^r}^{\otimes n+m}[2m])
         = \tau_{\leq n}
           (Rj_*R\psi_*\mu_{p^r}^{\otimes n+m}[2m])
$$
and $\tau_{\leq n}Rj_*(\tr_{\psi})$.
The vertices of this diagram are concentrated in degrees $\leq n$.
Hence we are reduced to the commutativity of
       the diagram of the $n$-th cohomology sheaves,
         which one can check by taking
            a section $s : X \hra \P_X$ of $f$ and
              using the compatibility proved in Step 1
                (see also Remark \ref{rem:newtrace} (2)).
More precisely,
    using \eqref{van:hagihara} one can construct a Gysin map
$$
\begin{CD}
    \gys_s :  R^nj_*\mu_{p^r}^{\otimes n} @>>>
       R^mf_*(R^{m+n}\beta_*\mu_{p^r}^{\otimes n+m})
    = {\cal H}^n(Rf_*(\tau_{\leq n+m }R\beta_*\mu_{p^r}^{\otimes n+m})[2m]),
\end{CD}
$$
   induced by that for $s_V : V \hra \P_V$.
One can further check that it is surjective by Theorem \ref{thm:hyodo} with
$r=1$
          and \cite{gros:purity}, I.2.1.5, I.2.2.3.
The details are straight-forward and left to the reader.
\end{pf}
Now we turn to the proof of Theorem \ref{thm:comp}.
Replacing $X$ by $\Spec(\O_{X,x})$,
      we assume that $X$ is local with closed point $x$.
Suppose that $S(=\ol {\{ y \}})$ is not normal, and
    let ${\frak n} : T \ra S$ be the normalization of $S$.
Since ${\frak n}$ is finite,
    the composite $T \ra S \hra X$ is projective, i.e.,
    factors as $T \os{i}{\hookrightarrow} \P_X^m=:\P \os{f}{\ra} X$
          with $i$ closed, for some $m \geq 1$.
Let $\psi : \P_V^m \to V$ be the morphism induced by $f$.
Let $T_x$ be the fiber ${\frak n}^{-1}(x) \subset T$ with reduced structure,
    and let $h : T_x \to x$ be the natural map.
Consider the diagram in Figure 1 below,
       where we wrote $\gys$ for Gysin maps
         for simplicity.
\begin{figure}[htpn]
\setlength{\unitlength}{.6mm}
{\scriptsize
\begin{picture}
(220,90)(-125,-45)
\thinlines
          \put(-88,-32){\vector(2,1){13}}
          \put(57,-32){\vector(-2,1){13}}
          \put(57,32){\vector(-2,-1){13}}
          \put(-88,32){\vector(2,-1){13}}
          \put(-110,34){\line(0,-1){68}}
          \put(-109,34){\line(0,-1){68}}
          \put(80,34){\vector(0,-1){68}}
          \put(-65,40){\vector(1,0){100}}
          \put(-65,-39){\vector(1,0){100}}
          \put(-65,13){\vector(0,-1){26}}
          \put(35,13){\vector(0,-1){26}}
          \put(-35,18){\vector(1,0){36}}
          \put(-35,-18){\vector(1,0){36}}
          \put(-77,-32){{\tiny $\gys$}}
          \put(38,-32){{\tiny $\gys$}}
          \put(38,29){{\tiny $\gys$}}
          \put(-77,29){{\tiny $\gys$}}
          \put(84,0){{\tiny $\tr_{h}^{n-b}$}}
          \put(-22,44){{\tiny $-\partial^{\val}_{y,T_x}$}}
          \put(-22,-45){{\tiny $-\partial^{\val}_{y,x}$}}
          \put(-63,0){{\tiny $\tr_{\psi}^n$}}
          \put(26,0){{\tiny $\tr_{f}^n$}}
          \put(-38,12){{\tiny $\delta^{\loc}_{y,T_x}(\T_r(n+m)_{\P})$}}
          \put(-32,-14){{\tiny $\delta^{\loc}_{y,x}(\T_r(n)_X)$}}
          \put(-18,28){{\tiny $(1)$}}
          \put(-90,0){{\tiny $(2)$}}
          \put(-18,0){{\tiny $(3)$}}
          \put(53,0){{\tiny $(4)$}}
          \put(-18,-30){{\tiny $(5)$}}
\put(59,39){$\H^{0}(T_x,\logwitt {T_x} r {n-b})$}
\put(-130,39){$\H^{n-b+1}(y,\mu_{p^r}^{\otimes n-b+1})$}
\put(-130,-41){$\H^{n-b+1}(y,\mu_{p^r}^{\otimes n-b+1})$}
\put(61,-41){$\H^{0}(x,\logwitt x r {n-b})$}
\put(6,17){$\H^{n+2m+b}_{T_x}(\P,\T_r(n+m)_{\P})$}
\put(-98,17){$\H^{n+2m+b-1}_y(\P_V,\mu_{p^r}^{\otimes m+n})$}
\put(-84,-20){$\H^{n+b-1}_y(V,\mu_{p^r}^{\otimes n})$}
\put(12,-20){$\H^{n+b}_x(X,\T_r(n)_X)$}
\end{picture}
\caption{diagram for Step 3}
}
\end{figure}
In this diagram, the outer large square
     commutes by the definition of
       $\delta_{y,x}^{\val}$ and \cite{jss}, Lemma A.1.
The square (1) commutes by Step 1.
The square (2) commutes by \cite{sga4.5}, Cycle, 2.3.8 (i).
The square (3) commutes by the property $j^*(\tr_f^n)=\tr_{\psi}^n$ of
$\tr_f^n$.
The square (4) commutes by Lemma \ref{lem:trace2} (2)
            and Remark \ref{rem:newtrace} (1).
Hence the square (5) commutes,
      which is the commutativity of \eqref{p-diagram}.
This completes the proof of Theorems \ref{thm:comp} and \ref{thm:gysin}.
    \qed
\medskip
\section{\bf Covariant functoriality and relative duality}\label{sect7}
\medskip
In this section, we prove the covariant functoriality {\bf T7}
    in Theorem \ref{thm0-11} and prove a relative duality result
     (see Theorem \ref{thm:relative} below).
Throughout this section, we work with the setting in \S\ref{sect4.1}.
Let $X$ and $Z$ be integral schemes which are flat of finite type over $B$
      and satisfy \ref{cond2}, and
   let $f : Z \ra X$ be a separated morphism of finite type.
Put $c:=\dim(X)-\dim(Z)$, and
    let $\psi:Z[1/p] \ra X[1/p]=V$ be the morphism induced by $f$.
By the absolute purity \cite{Th}, \cite{fujiwara},
     there is a trace morphism
$$
\begin{CD}
    \tr_{\psi}^n: R\psi_!\mu_{p^r}^{\otimes n-c}[-2c] @>>> \mu_{p^r}^{\otimes n}
      \quad \hbox{ in } \; D^b(V_{\et},\Z/p^r\Z),
\end{CD}
$$
which extends the trace morphisms for flat morphisms due to Deligne
     \cite{sga4}, XVIII.2.9 and satisfies the transitivity property.
\subsection{Covariant functoriality}\label{sect7.1}
The first result of this section is the following:
\begin{thm}[Covariant funtoriality]
\label{prop:trace}
For $f : Z \ra X$ as before,
    there is a unique morphism
$$
\begin{CD}
\tr_f^n: Rf_!\T_r(n-c)_Z[-2c] @>>> \T_r(n)_X
      \quad \hbox{ in } \; D^b(X_{\et},\Z/p^r\Z)
\end{CD}
$$
that extends $\tr_{\psi}^n$.
Consequently, these trace morphisms satisfy the transitivity property.
\end{thm}
\noindent
This theorem will be proved in the next subsection.
In this subsection,
    we prove the following:
\begin{lem}
\label{lem:semipurity}
Let $k$ be a perfect field of characteristic $p>0$, and
let $Y$ be a normal crossing variety over $\Spec(k)$.
Let $g : T \to Y$ be a separated morphism of finite type of schemes,
    and assume that $T$ has dimension $\le a$.
Put $c:= \dim(Y) - a$.
Assume that ${\frak L} \in D^b(W_{\et},\Z/p^r\Z)$
    and ${\frak M} \in D^-(Y_{\et},\Z/p^r\Z)$ are concentrated in degrees
$\le \ell$
    and $\le m$, respectively.
Let $\FF$ be a locally free $(\O_Y)^p$-module of finite rank.
Then for an integer $q < c - \ell - m$,
     we have
\begin{align*}
\Hom_{D^-(Y_{\et},\Z/p^r\Z)}((Rg_!{\frak L}) \otimes^{\L}{\frak M},\,
   \nu_{Y,r}^n[q]) &=0 \quad \hbox{ and }\\
\Hom_{D^-(Y_{\et},\Z/p^r\Z)}((Rg_!{\frak L}) \otimes^{\L}{\frak M},\,
   \FF[q]) &=0.
\end{align*}
\end{lem}
\begin{pf}
We prove the assertion only for $\nu_{Y,r}^n$.
One can check the assertion for $\FF$
    by repeating the same arguments as for $\nu_{Y,r}^n$,
      using Lemma \ref{hclaim1} below.

We first prove the case that ${\frak L}$ has constructible cohomology sheaves.
If there are closed subschemes $\phi_i: T_i \hra T$ ($i=1,2$)
    such that $T=T_1 \cup T_2$ and $\dim(T_1 \cap T_2) \le a-1$,
      then there is a distinguished triangle of the form
$$
   R\phi_{12*}\phi_{12}^*{\frak L}[-1] \lra
   {\frak L} \lra R\phi_{1*}\phi_1^*{\frak L} \oplus R\phi_{2*}\phi_2^*{\frak L}
    \lra
     R\phi_{12*}\phi_{12}^*{\frak L}
$$
in $D^b(T_{\et},\Z/p^r\Z)$, where
    $\phi_{12}$ denotes the closed immersion $T_1 \cap T_2 \hra T$.
Hence by induction on $a \ge 0$, we may assume that $T$ is irreducible.
Let $b$ be the dimension of $g(T) \subset Y$.
Noting that $Rg_!{\frak L}$ has constructible cohomology sheaves
     by \cite{sga4}, XIV.1.1, we prove
\begin{sublem}\label{sublem:support}
For $i \in \Z$,
    the support of $R^ig_!{\frak L}$ has dimension
     $\le \min\{b,a+\ell-i\}$,
        i.e.,
         there is a closed subscheme $\phi : W \hra Y$ of dimension
           $\le \min\{b,a+\ell-i\}$ for which we have
       $R^ig_!{\frak L}=\phi_*\phi^*R^ig_!{\frak L}$.
\end{sublem}
\begin{pf*}
{\it Proof of Sublemma \ref{sublem:support}}
Without loss of generality, we may assume that $g$ is proper
    (hence $Rg_!=Rg_*$).
Since $R^ig_*{\frak L}$ is zero outside of $g(T)$,
    the support of $R^ig_*{\frak L}$ has dimension $\le b$.
We show that the support of $R^ig_*{\frak L}$ is at most
   $(a+\ell-i)$-dimensional.
Let $y$ be a point on $g(T)$.
Put $T_y:=T \times _Y y=g^{-1}(y)$.
Since $\dim \ol {\{ y \}} + \dim T_y$
    equals the dimension of the closure of $T_y$ in $T$
      (\cite{sga4}, XIV.2.3 (iii)), we have
$\dim \ol {\{ y \}} + \dim T_y \le a$.
Now suppose that $(R^ig_*{\frak L})_{\ol y} \ne 0$.
Because $(R^qg_*{\frak L})_{\ol y}$
   is zero for $q > \dim T_y + \ell$ (loc.\ cit., XII.5.2, X.5.2),
    we have
$$
i \le \dim T_y + \ell \le a - \dim \ol {\{ y \}} + \ell,
$$
    that is, $\dim \ol {\{ y \}} \le a + \ell -i$.
Thus we obtain the sublemma.
\end{pf*}
We now turn to the proof of the lemma and
   compute a spectral sequence
$$
E_2^{u,v}=\Ext_{Y,\Z/p^r\Z}^u
   ((R^{-v}g_!{\frak L})\otimes^{\L} {\frak M},\nu_{Y,r}^n)
     \Lra \Ext_{Y,\Z/p^r\Z}^{u+v}
   ((Rg_!{\frak L})\otimes^{\L} {\frak M},\nu_{Y,r}^n).
$$
For $(u,v)$ with $u+v < c-\ell-m$ and
      $b \le a + \ell + v$, we have $u+m < \dim(Y)-b$ and
\begin{align*}
    E_2^{u,v}&=
     \Ext_{Y,\Z/p^r\Z}^u
    ((\phi_*\phi^*R^{-v}g_!{\frak L})\otimes^{\L} {\frak M},
     \nu_{Y,r}^n) \\
    &= \Ext_{Y,\Z/p^r\Z}^u
    (R\phi_*\phi^*(R^{-v}g_!{\frak L}\otimes^{\L} {\frak M}),
     \nu_{Y,r}^n) \\
    &= \Ext_{W,\Z/p^r\Z}^u
    (\phi^*(R^{-v}g_!{\frak L}\otimes^{\L} {\frak M}),R\phi^!\nu_{Y,r}^n)
     = 0
\end{align*}
by Theorem \ref{thm:wpurity} and Lemma \ref{lem:CD}.
Here $W$ denotes the closure of $g(T)$ and
    $\phi$ denotes the closed immersion $W \hra Y$.
For $(u,v)$ with $u+v < c-\ell-m$ and
      $b > a + \ell + v$, we have $u+m < \dim(Y)-(a+\ell+v)$.
There is a closed subscheme $\phi : W \hra Y$ of codimension
           $\ge \dim(Y)-(a+\ell+v)$ such that
       $R^{-v}g_!{\frak L}=\phi_*\phi^*R^{-v}g_!{\frak L}$,
        by the sublemma.
Hence
$$
    E_2^{u,v}=
     \Ext_{Y,\Z/p^r\Z}^u
   (\phi^*(R^{-v}g_!{\frak L}\otimes^{\L} {\frak M}),R\phi^!\nu_{Y,r}^n)
     = 0
$$
again by Theorem \ref{thm:wpurity} and Lemma \ref{lem:CD}.
Thus we obtain the assertion.

We next prove the case that ${\frak L}$ is general.
Take a bounded complex of $\Z/p^r\Z$-sheaves ${\frak L}^{\bullet}$
    which is concentrated in degrees $\le \ell$
     and represents ${\frak L}$.
Take a filtered inductive system
    $\{{\frak L}^{\bullet}_{\lam} \}_{\lam \in \Lam}$
      consisting of bounded complexes
        of constructible $\Z/p^r\Z$-sheaves
    which are concentrated in degrees $\le \ell$
     and whose limit is ${\frak L}^{\bullet}$.
Then for $q<c-\ell-m$, we have
\begin{align*}
\Ext_{Y,\Z/p^r\Z}^q
   ((Rg_!{\frak L})\otimes^{\L} {\frak M},\nu_{Y,r}^n)
    &= \Ext_{T,\Z/p^r\Z}^q
     ({\frak L}, Rg^!R\sHom({\frak M},\nu_{Y,r}^n)) \\
    &= \varprojlim_{\lam \in \Lam} \; \Ext_{T,\Z/p^r\Z}^q
   ({\frak L}^{\bullet}_{\lam}, Rg^!R\sHom({\frak M},\nu_{Y,r}^n)) \\
    &= \varprojlim_{\lam \in \Lam} \; \Ext_{Y,\Z/p^r\Z}^q
    ((Rg_!{\frak L}^{\bullet}_{\lam})\otimes^{\L} {\frak M},\nu_{Y,r}^n) =0
\end{align*}
by the previous case.
Here $Rg^!$ denotes the twisted inverse image functor
    of Deligne \cite{sga4}, {\rm XVIII},
     and we have used the adjointness between $Rg^!$ and $Rg_!$.
The second equality follows from the vanishing of
    the groups $\Ext_{Y,\Z/p^r\Z}^{q-1}
    ((Rg_!{\frak L}^{\bullet}_{\lam})\otimes^{\L} {\frak M},\nu_{Y,r}^n)$
    for all $\lam \in \Lam$
    and a standard argument which is similar as for \eqref{eq:projlim} below.
This completes the proof of the lemma.
\end{pf}
As a sepecial case of Lemma \ref{lem:semipurity},
    we obtain
\begin{cor}[Semi-purity]
\label{cor:semipurity}
Under the same setting as in Lemma $\ref{lem:semipurity}$,
   we have
$$
\tau_{ \le c-1}Rg^!\nu_{Y,r}^n = \tau_{ \le c-1}Rg^!\FF = 0.
$$
\end{cor}
\begin{pf}
For $T'$ \'etale separated of finite type over $T$, ${\cal G} \in
\{\nu_{Y,r}^n,\FF\}$
    and $q \le c-1$, we have
$$
   \Hom_{D^+(T'_{\et},\Z/p^r\Z)}(\Z/p^r\Z,
    Rh^!{\cal G}[q])=
    \Hom_{D^b(Y_{\et},\Z/p^r\Z)}(Rh_!\Z/p^r\Z,
    {\cal G}[q]) = 0
$$
by Lemma \ref{lem:semipurity},
where $h$ denotes the composite map $T' \to T \to Y$.
Hence $\tau_{ \le c-1}Rg^!{\cal G} = 0$.
\end{pf}
\subsection{Proof of Theorem \ref{prop:trace}}\label{sect7.2}
Let $j : V \hra X$ and $\iota : Y \hra X$ be as in \S\ref{sect4.1}.
Put ${\frak L}:=\T_r(n-c)_Z[-2c]$, for simplicity.
We first show
\begin{equation}\label{van:semipurity}
\Hom_{D^b(X_{\et},\Z/p^r\Z)}(Rf_!{\frak L},
   R\iota_*R\iota^!\T_r(n)_X)=0.
\end{equation}
Indeed, for $g: T := Z \times_X Y \to Y$ induced by $f$,
    we have
$$
\Hom_{D^b(X_{\et},\Z/p^r\Z)}(Rf_!{\frak L},
   R\iota_*R\iota^!\T_r(n)_X)
    =
     \Hom_{D^b(Y_{\et},\Z/p^r\Z)}(Rg_!\alpha^*{\frak L},R\iota^!\T_r(n)_X)
$$
by the adjointness between $\iota^*$ and $R\iota_*$
    and the proper base-change theorem: $\iota^*Rf_!=Rg_!\alpha^*$,
      where $\alpha$ denotes the closed immersion $T \hra Z$.
The latter group is zero by
    Lemma \ref{lem:semipurity} and a similar argument
    as for the vanishing \eqref{van:hagihara}.
By \eqref{van:semipurity} and Lemma \ref{lem:CD2} (1),
     it remains to show that the composite morphism
\begin{equation}\label{trace:comp}
\hspace{-50pt}
\begin{CD}
Rf_!{\frak L} @>{Rj_*(\tr_{\psi}^n)}>>
     Rj_*\mu_{p^r}^{\otimes n}
       @>{\delta^{\loc}_{V,Y}(\T_r(n)_X)}>> R\iota_*R\iota^!\T_r(n)_X[1]
\end{CD}
\hspace{-50pt}
\end{equation}
is zero in $D^b(X_{\et},\Z/p^r\Z)$.
We show the following:
\addtocounter{thm}{2}
\begin{lem}\label{sublem:semipurity}
\begin{enumerate}
\item[(1)]
Let $\{Z_{\lam}\}_{\lam \in \Lam}$ be an open covering of $Z$ with $\Lam$ 
finite,
  and let $f_{\lam}:Z_{\lam} \ra X$ be the composite map
      $Z_{\lam} \hra Z \ra X$ for each $\lam \in \Lam$.
Then the adjunction map
\begin{align*}
\Hom_{D^b(X_{\et},\Z/p^r\Z)}&(Rf_!{\frak L},
   R\iota_*R\iota^!\T_r(n)_X[1]) \\
    &\lra
    \bigoplus{}_{\lam \in \Lam}~
    \Hom_{D^b(X_{\et},\Z/p^r\Z)}(Rf_{\lambda!}({\frak L}\vert_{Z_{\lambda}}),
     R\iota_*R\iota^!\T_r(n)_X[1])
\end{align*}
   is injective.
\item[(2)]
Assume that $f$ is flat.
Let $Y' \subset Y$ be a closed subscheme of codimension $\ge 1$.
Put $U := X \setminus Y'$.
Then the following natural restriction map is injective$:$
\begin{align*}
\Hom_{D^b(X_{\et},\Z/p^r\Z)}&(Rf_!{\frak L},
    R\iota_*R\iota^!\T_r(n)_X[1]) \\
   & \lra
    \Hom_{D^b(U_{\et},\Z/p^r\Z)}\big( (Rf_!{\frak L})\vert_U,
      (R\iota_*R\iota^!\T_r(n)_X[1])\vert_U \big).
\end{align*}
\end{enumerate}
\end{lem}
\begin{pf*}
{\it Proof of Lemma \ref{sublem:semipurity}.}
(1) It suffices to consider the case that $\Lam=\{1,2\}$.
    Put $Z_{12}:=Z_1 \cap Z_2$.
   Let $f_{12}$ be the composite map $Z_{12} \hra Z \ra X$.
There is a distinguished triangle of the form
$$
\begin{CD}
   Rf_{12!}({\frak L}\vert_{Z_{12}}) @>>>
    Rf_{1!}({\frak L}\vert_{Z_1}) \oplus
     Rf_{2!}({\frak L}\vert_{Z_2})
     @>>> Rf_!{\frak L}
      \lra Rf_{12!}({\frak L}\vert_{Z_{12}})[1].
\end{CD}
$$
Hence the assertion follows from the vanishing \eqref{van:semipurity}
    for $f_{12}: Z_{12} \to X$.
\par
(2)
Let $\phi$ be the closed immersion $Y' \hra X$.
The kernel of the map in question is a quotient of
$\Hom_{D^b(X_{\et},\Z/p^r\Z)}(Rf_!{\frak L},
      R\phi_*R\phi^!\T_r(n)_X[1])$.
One can check that it is zero by a similar argument as for 
\eqref{van:semipurity},
   noting that $Y' \times_X Z$ has codimension $\ge 1$ in $T$ by
   the flatness of $f$.
\end{pf*}
We show that the composite morphism \eqref{trace:comp} is zero.
We first assume that $Z=\P_X^m$ ($m \ge 1$) and
   that $f$ is the natural projection.
By Lemma \ref{sublem:semipurity} (2),
  we are reduced to the case that $Y$ is smooth.
In this case, \eqref{trace:comp} is zero by Lemma \ref{lem:trace2}.
We next prove the general case.
By Lemma \ref{sublem:semipurity} (1), we may assume that $f$ is affine.
Take a decomposition
     $Z  \os{i}{\hra} \P_X^m=:\P \os{h}{\ra} X$
       of $f$ for some integer $m \geq 0$, where $i$ is a locally closed 
immersion.
We have morphisms
$$
\begin{CD}
   Rf_!{\frak L} @>{Rh_*(\gys_i^{m+n})[2m]}>>
   Rh_*\T_r(m+n)_{\P}[2m] @>{\tr_h^n}>> \T_r(n)_X,
\end{CD}
$$
where $\tr_h^n$ is obtained from the vanishing of
   \eqref{trace:comp} for $h$.
See Theorem \ref{thm:gysin} for $\gys_i^{m+n}$.
Since this composite morphism extends $\tr_{\psi}^n$,
  we see that \eqref{trace:comp} is zero.
This completes the proof of Theorem \ref{prop:trace}.
\qed
\par
\medskip
\noindent
The following corollary is a horizontal variant of
      Proposition \ref{prop:proj}:
\begin{cor}[Projection formula]\label{cor:proj}
For $f:Z \to X$ as before, the diagram
$$
{\small
\begin{CD}
Rf_!\T_r(m-c)_Z[-2c]\otimes^{\L} \T_r(n)_X
       @>{\tr^{m}_f \otimes^{\L} \id}>>
           \T_r(m)_X \otimes^{\L} \T_r(n)_X
            @>{\eqref{def:prod2}}>>\T_r(m+n)_X \\
     @V{\id \otimes^{\L} f^*}VV @. @|\\
Rf_!\T_r(m-c)_Z[-2c] \otimes^{\L} Rf_*\T_r(n)_Z
   @>{\eqref{def:prod2}}>>  Rf_!\T_r(m+n-c)_Z[-2c]
            @>{\tr^{m+n}_f}>> \T_r(m+n)_X
\end{CD}
}
$$
commutes in $D^-(X_{\et},\Z/p^r\Z)$.
See Proposition $\ref{prop:funct}$ for
   $f^*: \T_r(n)_X \to Rf_*\T_r(n)_Z$.
\end{cor}
\begin{pf}
Because the diagram in question commutes on $X[1/p]$,
    the assertion follows from Lemma \ref{lem:CD2} (1) and a vanishing
$$
\Hom_{D^-(X_{\et},\Z/p^r\Z)}
   (Rf_!\T_r(m-c)_Z[-2c]\otimes^{\L} \T_r(n)_X,
    R\iota_*R\iota^!\T_r(m+n)_X)=0,
$$
which one can check by Lemma \ref{lem:semipurity}
    and a similar argument as for \eqref{van:semipurity}.
\end{pf}
\subsection{Relative duality}\label{sect7.3}
Let $f: Z \to X$ be as before.
Let $j : V \hra X$ and $\iota: Y \hra X$ be as in \S\ref{sect4.1}.
Let $T$ be the divisor on $Z$ defined by the radical of $(p) \subset \O_Z$.
There is a commutative diagram of schemes
$$
\begin{CD}
     Z[1/p]  @>{\beta}>> Z @<{\alpha}<< T \\
    @VV{\quad\quad\;\;\,\square}V
       @V{f}VV  @V{g}VV\\
     V   @>{j}>>  X @<{\iota}<<  Y.
\end{CD}
$$
Put $d:=\dim(X)$, $b:=\dim(Z)$ and $c:=d-b$.
We prove the following result, which was included in the earlier version
  of \cite{jss}:
\begin{thm}[Relative duality]
\label{thm:relative}
\begin{enumerate}
\item[(1)]
$\tr_f^{d}$ induces an isomorphism
$$
\begin{CD}
     \tr^f : \T_r(b)_Z[-2c] @>{\simeq}>> Rf^!\T_r(d)_X
      \quad \hbox{ in } \;
       D^b(Z_{\et},\Z/p^r\Z).
\end{CD}
$$
\item[(2)]
There is a commutative diagram in $D^b(X_{\et},\Z/p^r\Z)$
$$
\begin{CD}
    R\iota_*Rg_!\nu_{T,r}^{b-1}[b-1]
     @>{(\natural)}>> Rf_!\T_r(b)_Z[2b] \\
     @V{R\iota_*(\tr_g)}VV @VV{\tr_f^{d} [2d]}V\\
    \iota_*\nu_{Y,r}^{d-1}[d-1] @>{\gys_{\iota}^{d}[2d]}>>
       \T_r(d)_X[2d],\\
\end{CD}
$$
where $\tr_g$ denotes the trace morphism
   in Remark $\ref{rem:newtrace}$ $(3)$, and
   the arrow $(\natural)$ is induced by
    the isomorphism $R\iota_*Rg_!=Rf_!R\alpha_*$ and
     the Gysin morphism $\gys_{\alpha}^{b}$.
\end{enumerate}
\end{thm}
\noindent
To prove this theorem,
  we first note a standard fact
       (cf.\ \cite{jss}, Lemma 3.8).
\begin{lem}\label{lem4''}
For a torsion sheaf $\FF$ on $V_{\et}$
    and an integer $q > d$,
   we have $R^qj_*\FF=0$.
\end{lem}
\noindent
As an immediate consequence of
     Lemma \ref{lem4''} and \eqref{p-sup},
    we obtain
\begin{lem}\label{lem4'}
The natural morphism
$
\tau_{\leq d+1}R\iota^!\T_r(d)_X \to R\iota^!\T_r(d)_X
$
      is an isomorphism in $D^b(Y_{\et},\Z/p^r\Z)$.
Consequently, $\gys_{\iota}^{d} :
    \nu_{Y,r}^{d-1}[-d-1] \to R\iota^!\T_r(d)_X$ is an isomorphism.
\end{lem}
\begin{pf*}
{\it Proof of Theorem \ref{thm:relative}}
We first show (2).
By Lemma \ref{lem:semipurity}
  and a similar argument as for the vanishing of \eqref{trace:comp},
   one can reduce the assertion to
    Lemma \ref{lem:trace2} (2) and Corollary \ref{cor:trans}
      (see also Remark \ref{rem:newtrace} (3)).
The details are straight-forward and left to the reader.
We next show (1).
Let $\tr^f : \T_r(b)_Z[-2c] \ra Rf^!\T_r(d)_X$
   be the adjoint morphism of $\tr_f^{d}$.
Because $\beta^*(\tr^f)$ is an isomorphism by
   the absolute purity (\cite{Th}, \cite{fujiwara}),
   we have only to show that $R\alpha^!(\tr^f)$ is an isomorphism.
By (2), there is a commutative diagram in $D^+(T_{\et},\Z/p^r\Z)$
$$
\begin{CD}
    \nu_{T,r}^{b-1}[b-1]
     @>{\gys_{\alpha}^b[2b]}>> R\alpha^!\T_r(b)_Z[2b] \\
     @V{\tr^g}VV @VV{R\alpha^!(\tr^f)[2d]}V\\
    Rg^!\nu_{Y,r}^{d-1}[d-1] @>{Rg^!(\gys_{\iota}^{d})[2d]}>>
       R\alpha^!Rf^!\T_r(d)_X[2d].\\
\end{CD}
$$
The horizontal arrows are isomorphisms by Lemma \ref{lem4'}
    for $Z$ and $X$, respectively.
The left vertical arrow, defined as the adjoint morphism
    of $\tr_g$, is an isomorphism by \cite{jss}, Theorem 2.8.
Consequently, $R\alpha^!(\tr^f)$ is an isomorphism.
This completes the proof of Theorem \ref{thm:relative}.
\end{pf*}
\begin{rem}\label{rem:jss}
By Theorem $\ref{thm:relative}$,
    $\T_r(d)_X[2d]$ is canonically isomorphic to
    the object ${\mathcal D}_{X,p^r} \in D^b(X_{\et},\Z/p^r\Z)$
      considered in \cite{jss}, Theorem $4.4$.
\end{rem}
\medskip
\section{\bf Explicit formula for $p$-adic vanishing cycles}\label{sect8}
\medskip
In this section we construct a canonical pairing on the sheaves
    of $p$-adic vanishing cycles in the derived category,
   and prove an explicit formula for that pairing,
   which will be used in \S\ref{sect9}.
\subsection{Setting}\label{sect8.0}
The setting is the same as in \S\ref{sect3.2}.
Note the condition \ref{cond0} assumed there.
We further assume that
            {\it $K$ contains a primitive $p$-th root of unity}
     and that {\it $k$ is finite}.
We put
\begin{equation}\label{def:simple}
\begin{CD}
      \nu_Y^n := \nu_{Y,1}^n,
                    \quad
        \mu':=\iota^*j_*\mu_p \quad \hbox { and }\quad
      \mu:=\mu_p(K)
\end{CD}
\end{equation}
      for simplicity.
Note that $\mu'$ is the constant \'etale sheaf on $Y$
           associated with the abstract group $\mu (\simeq \Z/p\Z)$,
      because the sheaf $\mu_p$ on $X_K$ is constant
        and the sheaf $j_*\mu_p$ on $X$ is also constant
             by the normality of $X$ (cf.\ \cite{ts}, 1.5.1).
Now let $N$ be the relative dimension
      $\dim(X/O_K)$.
Let $n$ be a positive integer with $1 \leq n \leq N+1$.
Put $n' := N+2-n$,
       $M^q:=M^q_1= \iota^*R^qj_{*}\mu_{p}^{\otimes q}$,
     and let $U^{\bullet}$ be the filtration on $M^q$
       defined in Definition \ref{def:vcyc}.
The aim of this section
     is to construct a morphism
$$
\begin{CD}
\Theta^n : U^1M^n \otimes U^1M^{n'}[-N-2] @>>>
           \mu' \otimes \nu_{Y}^N [-N-1]
             \quad \hbox{ in } \; D^b(Y_{\et},\Z/p\Z)
\end{CD}
$$
   and to prove an explicit formula for this morphism
    (cf.\ Theorem \ref{lem:trace} below).
\subsection{Construction of $\Theta^n$}\label{sect8.1}
Because $\mu'$
is (non-canonically) isomorphic to
    the constant sheaf $\Z/p\Z$,
      we will write
        $\mu' \otimes \K$ ($\K \in D^-(Y_{\et},\Z/p\Z)$)
   for the derived tensor product $\mu' \otimes^{\L} \K$
            in $D^-(Y_{\et},\Z/p\Z)$.
For $q$ with $1 \leq q \leq N+1$,
    fix a distinguished triangle
$$
\begin{CD}
   (M^q/U^1M^q) [-q-1] @>{g'}>>
     \AA(q)  @>{t'}>> \tau_{\leq q}\iota^*Rj_*\mu_{p}^{\otimes q}
        @>>> (M^q/U^1M^q) [-q],
\end{CD}
$$
where the last arrow is defined as
    the composite
   $\tau_{\leq q}\iota^*Rj_*\mu_{p}^{\otimes q}
      \ra M^q[-q] \ra (M^q/U^1M^q) [-q]$.
Clearly, $\AA(q)$ is concentrated in $[0,q]$, and
   the triple $(\AA(q),t',g')$ is unique up to a unique isomorphism
      (and $g'$ is determined by $(\AA(q),t')$)
        by Lemma \ref{lem:CD2} (3).
We construct $\Theta^n$ by decomposing the morphism
\begin{equation}\label{tau1}
\begin{CD}
\AA(n) \otimes^{\L}
    \AA(n')
     @>>> (\tau_{\leq n}\iota^*Rj_*\mu_{p}^{\otimes n})
   \otimes^{\L} (\tau_{\leq n'}\iota^*Rj_*\mu_{p}^{\otimes n'}) \\
   @>>>  \iota^*Rj_*\mu_{p}^{\otimes N+2}
\end{CD}
\end{equation}
   induced by the natural isomorphism $\mu_{p}^{\otimes n}
            \otimes \mu_{p}^{\otimes n'} \simeq \mu_p^{\otimes N+2}$
    in characteristic zero.
By Lemma \ref{lem4''} and the assumption
    that $\zeta_p \in K$, there is a morphism
\begin{equation}\label{tau2}
\hspace{-50pt}
\begin{CD}
   \iota^*Rj_*\mu_{p}^{\otimes N+2}
      \;{\simeq}\;
        \mu' \otimes (\tau_{\leq N+1}\iota^*Rj_*\mu_{p}^{\otimes N+1})
       @>{\id \otimes \eqref{sigma(n)}}>> \mu' \otimes \nu_Y^N [-N-1],
\end{CD}
\hspace{-50pt}
\end{equation}
which, together with \eqref{tau1},
    induces a morphism
\begin{equation}\label{tau3}
\hspace{-50pt}
\begin{CD}
\AA(n) \otimes^{\L} \AA(n')
         @>>> \mu' \otimes \nu_Y^N [-N-1].
\end{CD}
\hspace{-50pt}
\end{equation}
Noting that $\AA(q)$ is concentrated in $[0,q]$
     with ${\cal H}^q(\AA(q)) \simeq U^1M^q$,
    we show the following:
\addtocounter{thm}{3}
\begin{lem}
There is a unique morphism
\stepcounter{equation}
\begin{equation}\label{tau4}
\hspace{-50pt}
\begin{CD}
\AA(n) \otimes^{\L} (U^1M^{n'}[-n'])
     @>>> \mu' \otimes \nu_Y^N
[-N-1]  \quad \hbox { in } \; D^-(Y_{\et},\Z/p\Z)
\end{CD}
\hspace{-50pt}
\end{equation}
that the morphism \eqref{tau3} factors through.
\end{lem}
\begin{pf}
There is a distinguished triangle of the form
$$
{\small
\begin{CD}
   \AA(n) \otimes^{\L} (\tau_{\leq n'-1}\AA(n'))
     \ra \AA(n) \otimes^{\L}\AA(n')
     \ra \AA(n) \otimes^{\L}(U^1M^{n'}[-n']) \ra
       \AA(n) \otimes^{\L}(\tau_{\leq n'-1}\AA(n'))[1].
\end{CD}
}
$$
By Lemma \ref{lem:CD2} (2),
      it suffices to show that (i) the morphism
$$
\begin{CD}
\AA(n) \otimes^{\L} (\tau_{\leq n'-1}\AA(n'))
   @>>> \mu' \otimes \nu_Y^N[-N-1]
\end{CD}
$$
induced by \eqref{tau3} is zero and that (ii) we have
$$
\begin{CD}
\Hom_{D^-(Y_{\et},\Z/p\Z)}
   (\AA(n) \otimes^{\L}(\tau_{\leq n'-1}\AA(n'))[1],
       \mu' \otimes \nu_Y^N [-N-1])=0.
\end{CD}
$$
The claim (ii) follows from Lemma \ref{lem:CD}.
As for the claim (i),
because $\AA(n) \otimes^{\L} (\tau_{\leq n'-1}\AA(n'))$
   is concentrated in degrees $\leq N+1$,
        the problem is reduced to
   the triviality of the induced map
            on the $(N+1)${\it st} cohomology sheaves
             (cf.\ Lemma \ref{lem:CD}).
One can check this by a similar argument as for Proposition \ref{prop:prod}.
\end{pf}
Applying a similar argument as for this lemma
     to the morphism \eqref{tau4},
    we obtain a morphism
\begin{equation}\label{tau5}
   \begin{CD}
(U^1M^{n}[-n]) \otimes^{\L} (U^1M^{n'}[-n'])
     @>>> \mu' \otimes \nu_Y^N [-N-1].
\end{CD}
\end{equation}
Finally because $\Z/p\Z$-sheaves are flat over $\Z/p\Z$,
   there is a natural isomorphism
$$
\hspace{-50pt}
\begin{CD}
   (U^1M^{n}[-n]) \otimes^{\L} (U^1M^{n'}[-n'])
     \simeq U^1M^{n} \otimes U^1M^{n'}[-N-2]
    \quad \hbox{ in } \; D^-(Y_{\et},\Z/p\Z)
\end{CD}
\hspace{-50pt}
$$
induced by the identity map on the $(n+n')${\it th} cohomology sheaves.
We thus define $\Theta^n$ by
    composing the inverse of this isomorphism and
    the morphism \eqref{tau5}.
\subsection{Explicit formula for $\Theta^n$}\label{sect8.2}
We formulate an explicit formula
        (see Theorem \ref{lem:trace} below)
        to calculate the morphism $\Theta^n$.
Let
$$
\begin{CD}
\chi: \mu' \otimes (\omega_Y^N/\BB_Y^N)
           @>>> \mu' \otimes \nu_Y^N[1]
\quad \hbox { in } \; D^b(Y_{\et},\Z/p\Z)
\end{CD}
$$
be the connecting morphism
    associated with a short exact sequence
   (\cite{hyodo}, (1.5.1))
$$
\begin{CD}
0 @>>> \mu' \otimes \nu_Y^N  @>{\id \otimes \mathrm{incl} }>>
     \mu' \otimes\omega_Y^{N}
       @>{\id \otimes (1-C^{-1})}>>
       \mu' \otimes (\omega_Y^N/\BB_Y^N) @>>> 0.
\end{CD}
$$
Here $\BB_Y^N$ denotes the image of
   $d:\omega_Y^{N-1} \ra \omega_Y^N$,
   $C^{-1}$ denotes the inverse Cartier operator
         defined in loc.\ cit., (2.5)
   (cf.\ \eqref{def:cartier} below)
   and we have used the isomorphism
    $\omega_{Y,\log}^N \simeq \nu_Y^N$
          in Remark \ref{rem:hyodo} (4).
We next construct a key map $f^{q,n}$
   (cf.\ Definition \ref{def:fqn} (2) below).
Let $e$ be the absolute ramification
         index of $K$ and put $e':=pe/(p-1)$.
Because $K$ contains primitive $p$-th roots of unity
    by assumption,
    $e'$ is an integer divided by $p$.
Fix a prime element $\pi \in O_K$.
Put $s:=\Spec(k)$.
Let ${\cal L}_Y$
   (resp.\ ${\cal L}_s$)
    be the log structure on $Y$
     (resp.\ on $s$)
       defined in \S\ref{sect3.3}.
We use the trivial log structure $s^{\times}$ on $s$
    and a map on $Y_{\et}$ analogous to \eqref{dlog:log2}
\begin{equation}\label{dlog:log3}
\begin{CD}
\dlog:  {\cal L}_Y^{\gp} @>>> \omega_{(Y,{\cal L}_Y)/(s,s^{\times})}^1.
\end{CD}
\end{equation}
\addtocounter{thm}{1}
\begin{rem}\label{rem:log}
\begin{enumerate}
\item[(1)]
The composite of \eqref{dlog:log3} with the
canonical projection
    $\omega_{(Y,{\cal L}_Y)/(s,s^{\times})}^1 \ra
   \omega_Y^1$
       agrees with the map $\dlog$ in \eqref{dlog:log2}.
\item[(2)]
Let $\ol {\pi}$ be the residue class of $\pi$ in
   ${\cal L}_Y$ under \eqref{log:surj}.
Then we have $\dlog(\ol {\pi})=0$
   in $\omega_Y^1$, but not in
             $\omega_{(Y,{\cal L}_Y)/(s,s^{\times})}^{1}$.
Indeed, by the definition of relative differential modules
      $($\cite{kk:log}, \cite{kf:log}$)$,
    there is a short exact sequence on $Y_{\et}$
$$
\begin{CD}
0 @>>> \O_Y @>{a
\mapsto a \cdot \dlog(\ol {\pi})}>>
        \omega_{(Y,{\cal L}_Y)/(s,s^{\times})}^{1}
           @>>> \omega_Y^1 @>>> 0.
\end{CD}
$$
The isomorphism \eqref{def:rho} below follows from this fact.
\end{enumerate}
\end{rem}
\noindent
\smallskip
Now let $n$ and $q$ satisfy
     $1 \leq n \leq N+1$ and $1 \leq q \leq e'-1$.
Put $n':=N+2-n$.
Let $U_{X_K}^q$ be the \'etale subsheaf of
   $\iota^*j_*\O^{\times}_{X_K}$ defined
    in Definition \ref{def:vcyc}, and put
\stepcounter{equation}
$$
\begin{CD}
   \Symb^{q,n}:= U_{X_K}^q
        \otimes \big( \iota^*j_*\O^{\times}_{X_K} \big)^{\otimes n-1}
    \otimes U_{X_K}^{e'-q}
        \otimes \big(\iota^*j_*\O^{\times}_{X_K} \big)^{\otimes n'-1}.
\end{CD}
$$
The sheaf $U^qM^n \otimes U^{e'-q}M^{n'}$
   is a quotient of $\Symb^{q,n}$
        (cf.\ Definition \ref{def:vcyc}):
$$
\begin{CD}
   U^qM^n \otimes U^{e'-q}M^{n'} = \Image \left(
    \Symb^{q,n} \lra U^1M^n \otimes U^1M^{n'}
       \right).
\end{CD}
$$
We define the homomorphism of \'etale sheaves
$$
\begin{CD}
F^{q,n} : \Symb^{q,n}
   \lra \omega_Y^N/\BB_Y^N
\end{CD}
$$
   by sending a local section
$(1+\pi^q
\alpha_1) \otimes
          (\otimes_{i=1}^{n-1}~\beta_i)
    \otimes
(1+\pi^{e'-q}\alpha_2)
        \otimes (\otimes_{i=n}^{N}~\beta_i)$
    with $\alpha_1, \alpha_2 \in \iota^*\O_X$
      and $\beta_1,\dotsc,\beta_N \in \iota^*j_*\O^{\times}_{X_K}$,
     to the following:
$$
\begin{CD}
   q \cdot \ol{\alpha_1\alpha_2} \cdot
          (\wedge_{i=1}^{N} ~\dlog \ol{\beta_i})
   + g^{-1}\left(\ol{\alpha_2} \cdot d \ol{\alpha_1} \wedge
   (\wedge_{i=1}^{N}~ \dlog \ol{\beta_i}) \right)
   \mod \BB_Y^N,
\end{CD}
$$
where for $x \in \iota^* \O_X$
            (resp.\ $x \in \iota^*j_* \O^{\times}_{X_K}$),
              $\ol x$ denotes its residue class in $\O_Y$
                (resp.\ in ${\cal L}_Y^{\gp}$ under \eqref{log:surj})
       and $g$ denotes the following $\O_Y$-linear
   isomorphism
          (cf.\ Remark \ref{rem:log} (2)):
\begin{equation}\label{def:rho}
\begin{CD}
g : \omega_Y^N @>{\simeq}>>
        \omega_{(Y,{\cal L}_Y)/(s,s^{\times})}^{N+1},
         \quad \omega \mapsto \dlog(\ol {\pi}) \wedge \omega.
\end{CD}
\end{equation}
\addtocounter{thm}{1}
\begin{lem}\label{lem:f}
Let $n$ and $q$ be as above.
Then $F^{q,n}$ factors through $U^qM^n \otimes U^{e'-q}M^{n'}$.
\end{lem}
\begin{pf}
Let $Y_{\sing}$ be the singular locus of $Y$,
    and let $j_Y$ be the open immersion $Y \sm Y_{\sing} \hra Y$.
Replacing $X$ by $X \sm Y_{\sing}$,
        we may assume that
      $Y$ is smooth over $s=\Spec(k)$,
       because $\omega_Y^N/\BB_Y^N$ is a locally free $(\O_Y)^p$-module
       and the canonical map
        $\omega_Y^N/\BB_Y^N \ra j_{Y*}j_Y^*(\omega_Y^N/\BB_Y^N)$
         is injective.
We show that $F^{q,n}$
       factors through
           $\gr_U^qM^n \otimes \gr_U^{e'-q}M^{n'}$,
         assuming that $Y$ is smooth.
For $m \geq 1$ and $\ell$ with
    $1 \leq \ell \leq e'-1$, let
\def\BK{\mathrm{BK}}
$$
\begin{CD}
\rho^{\ell,m} : \Omega_Y^{m-2} \oplus \Omega_Y^{m-1}
      @>>> \gr_U^{\ell}M^m
\end{CD}
$$
be the Bloch-Kato map (cf.\ \cite{bk}, (4.3)) defined as
$$
\begin{cases}
(\alpha
\cdot \dlog \beta_1 \wedge \dotsb \wedge \dlog \beta_{m-2},~0)
      \mapsto \{1+ \pi^{\ell} \wt{\alpha}, \wt{\beta_1}, \dotsc,
   \wt{\beta_{m-2}},\pi \} \hspace{-15pt} &\mod U^{{\ell}+1}M^m, \\
(0,~ \alpha \cdot \dlog \beta_1 \wedge \dotsb \wedge \dlog \beta_{m-1})
   \mapsto \{1+ \pi^{\ell} \wt{\alpha}, \wt{\beta_1}, \dotsc,
   \wt{\beta_{m-1}} \} &\mod U^{{\ell}+1}M^m \\
\end{cases}
$$
for $\alpha \in \O_Y$ and each $\beta_i \in \O^{\times}_Y$,
    where $\wt {\alpha} \in \O_X$ (resp.\ each $\wt{\beta_i} \in \O^{\times}_X$)
      denotes a lift of $\alpha$ (resp.\ $\beta_i$).
There are short exact sequences
\stepcounter{equation}
\begin{equation}\label{exact:gr}
\hspace{-60pt}
\begin{CD}
    0  @.~\lra~@. \Omega_Y^{m-2} @>{\theta^{{\ell},m}}>>
   \Omega_Y^{m-2} \oplus \Omega_Y^{m-1}
       @>{\rho^{{\ell},m}}>> \gr_U^{\ell}M^m @.~\lra~@. 0 @. \qquad @.
        (\hbox{if } p \hspace{-5pt}\not\vert {\ell}),\\
    0  @.~\lra~@. \ZZ_Y^{m-2} \oplus \ZZ_Y^{m-1} @>{{\mathrm{incl.}}}>>
       \Omega_Y^{m-2} \oplus \Omega_Y^{m-1}
       @>{\rho^{{\ell},m}}>> \gr_U^{\ell}M^m @.~\lra~@. 0 @. @.
        (\hbox{if } p \vert {\ell}),\phantom{i}
\end{CD}
\hspace{-50pt}
\end{equation}
where $\theta^{{\ell},m}$ is given by
      $\omega \mapsto ((-1)^m \cdot {\ell} \cdot \omega, d\omega)$
     (cf.\ \cite{bk}, Lemma (4.5)).
Let
$$
\begin{CD}
h^{{\ell},m}:
   U_{X_K}^{\ell} \otimes
     \big( \iota^*j_*\O^{\times}_{X_K} \big)^{\otimes m-1}
          @>>>  \Omega_Y^{m-2} \oplus \Omega_Y^{m-1}
\end{CD}
$$
be the map that sends $(1+\pi^{\ell} \alpha) \otimes
   (\otimes_{i=1}^{m-1}~\beta_i)$
    with $\alpha \in \iota^*\O_X$
      and $\beta_i \in \iota^*\O^{\times}_{X} \cup \{ \pi \}$ to
$$
{
\begin{CD}
\begin{cases}
     (0,\; \ol{\alpha} \cdot \wedge_{1 \leq i \leq m-1}~
       \dlog \ol{\beta_i})
        \quad &
         (\hbox{if }\beta_i \in \iota^*\O^{\times}_{X} \hbox{ for all } i), \\
   ((-1)^{m-1-i'} \cdot \ol{\alpha}
       \cdot \wedge_{1 \leq i \leq m-1, i \not= i'}~
       \dlog \ol{\beta_i}, \; 0)
         \qquad & (\hbox{if } \beta_{i} = \pi
           \hbox{ for exactly one } i=:i'),\\
     (0, \; 0) \quad & (\hbox{otherwise}).
\end{cases}
\end{CD}
}
$$
Here for $x \in \iota^* \O_X$
            (resp.\ $x \in \iota^*\O^{\times}_{X}$),
   $\ol x$ denotes its residue class in $\O_Y$
                (resp.\ in $\O^{\times}_Y$).
Now there is a diagram
$$
{\small
\begin{CD}
\Symb^{q,n}  @>{h^{q,n} \otimes h^{e'-q,n'}}>>
     (\Omega_Y^{n-2} \oplus \Omega_Y^{n-1})
        \otimes (\Omega_Y^{n'-2} \oplus \Omega_Y^{n'-1})
           @>{\rho^{q,n} \otimes \rho^{e'-q,n'}}>>
             \gr_U^qM^n \otimes \gr_U^{e'-q}M^{n'} \\
    @. @V{\varphi^{q,n}}VV @.\\
           @. \Omega_{Y}^N/\BB_Y^N, @.
\end{CD}
}
$$
where $\varphi^{q,n}$ is defined as
$$
\begin{CD}
   (\omega_1,\omega_2)\otimes(\omega_3,\omega_4) \mapsto
        q \cdot \omega_2 \wedge \omega_4 +
   (-1)^{n-1}\cdot(d\omega_1) \wedge \omega_4
       + (-1)^{n'-1} \cdot \omega_2 \wedge d\omega_3 \qquad \mod \BB_Y^N.
\end{CD}
$$
In this diagram, the composite of the top row agrees with the symbol map,
   and the composite of $h^{q,n} \otimes h^{e'-q,n'}$ and $\varphi^{q,n}$
   agrees with $F^{q,n}$ (see also Remark \ref{rem:log} (2)).
Hence to prove the assertion,
   it suffices to show that the subsheaf
$$
\begin{CD}
    \ker(\rho^{q,n} \otimes \rho^{e'-q,n'})
   \subset (\Omega_Y^{n-2} \oplus \Omega_Y^{n-1}) \otimes
   (\Omega_Y^{n'-2} \oplus \Omega_Y^{n'-1})
\end{CD}
$$
   has trivial image under $\varphi^{q,n}$,
       which follows from \eqref{exact:gr}
         with $(\ell,m)=(q,n), (e'-q,n')$
        (note that $p \vert q
   \Leftrightarrow p \vert (e'-q)$, because $p \vert e'$).
Thus we obtain Lemma \ref{lem:f}.
\end{pf}
\stepcounter{thm}
\begin{defn}\label{def:fqn}
\begin{enumerate}
\item[(1)]
For $\zeta \in \mu=\mu_p(K)$ with $\zeta \not = 1$,
    let $v(\zeta) \in k^{\times}$ be the residue class of
      $(1-\zeta)/\pi^{e/(p-1)} \in O_K^{\times}$.
We define $u:=\zeta \otimes v(\zeta)^{-p} \in \mu \otimes k$,
which is independent of the choice of $\zeta \not =1$.
\item[(2)]
Let ${\ul k}$ be the constant sheaf on $Y_{\et}$ associated with $k$.
We define the homomorphism
$$
\begin{CD}
f^{q,n}: U^qM^n \otimes U^{e'-q}M^{n'} @>>>
   (\mu' \otimes {\ul k}) \otimes_{\ul k} (\omega_Y^N/\BB_Y^N)
          \simeq \mu' \otimes (\omega_Y^N/\BB_Y^N)
\end{CD}
$$
   as $u \otimes_{k} (-1)^{N+n}\, \ol {F}{}^{q,n}$.
Here we regarded $u \in \mu \otimes k$
     as a global section of $\mu' \otimes {\ul k}$,
    and $\ol {F}{}^{q,n}$ denotes the map
    induced by $F^{q,n}$ $($cf.\ Lemma $\ref{lem:f})$.
\end{enumerate}
\end{defn}
\begin{rem}
    $f^{q,n}$ is independent of the choice of $\pi$
    by the definitions of $F^{q,n}$ and $u$.
\end{rem}
\noindent
Now we state the main result of this section.
\par
\begin{thm}[Explicit formula]\label{lem:trace}
Assume that $X$ is proper over $B$.
Then for $(q,n)$ with $1 \leq q \leq e'-1$ and $1 \leq n \leq N+1$,
  the following square commutes in $D^b(Y_{\et},\Z/p\Z)$$:$
$$
\begin{CD}
U^qM^n \otimes U^{e'-q}M^{n'} @>{f^{q,n}}>>
     \mu' \otimes (\omega_Y^N/\BB_Y^N) \\
@V{{\mathrm{canonical}}}VV @VV{\chi}V \\
U^1M^n \otimes U^1M^{n'} @>{\Theta^n[N+2]}>>
         \mu' \otimes \nu_Y^N[1].
\end{CD}
$$
\end{thm}
\medskip
\noindent
We will prove Theorem \ref{lem:trace} in
    \S\S\ref{sect8.3}--\ref{sect8.6} below.
We will first reduce the problem
     to an induced diagram of cohomology groups of $Y$ in \S\ref{sect8.3},
    and then to an induced diagram of cohomology groups
      of higher local fields in \S\ref{sect8.4}.
In \S\ref{sect8.5}, we will prove a Galois descent of
     invariant subgroups of Galois modules.
We will finish the proof of Theorem \ref{lem:trace}
   in \S\ref{sect8.6} by computing symbols,
    whose details are standard in higher local class field theory
      (cf.\ \cite{kk:cft2})
    but will be included for the convenience of the reader.
\subsection{Reduction to cohomology groups}\label{sect8.3}
In this step,
    we reduce Theorem \ref{lem:trace}
     to the equality \eqref{trace:eq} below.
We first show the following:
\begin{prop}\label{prop1-1}
Assume that $Y$ is proper over $\Spec(k)$.
Let $\FF$ be a $\Z/p\Z$-sheaf on $Y_{\et}$.
Then for $i \in \Z$, the Yoneda pairing
$$
\begin{CD}
\H^i(Y,\FF) \times
   \Ext^{N-i+1}_{Y,\Z/p\Z}(\FF,\nu_Y^N)
       @>>> \H^{N+1}(Y,\nu_Y^N) @>{\tr_Y}>> \Z/p\Z
\end{CD}
$$
$($see Theorem $\ref{thm:milne}$ for $\tr_Y)$
induces an isomorphism
$$
\Ext^{N-i+1}_{Y,\Z/p\Z}(\FF,\nu_Y^N)
      \simeq \Hom(\H^i(Y,\FF),\Z/p\Z).
$$
\end{prop}
\begin{pf}
If $\FF$ is constructible,
  then the isomorphism in question is an isomorphism of finite groups
      by the duality theorem of Moser \cite{moser}
    (note that the complex $\wt{\nu}_{r,Y}^{n}$ defined in loc.\ cit.\
     is quasi-isomorphic to the sheaf $\nu_{Y,r}^n$
      by \cite{sato:ss}, 2.2.5 (1)).
We prove the general case.
Write $\FF$ as a filtered inductive limit
   $\varinjlim_{\lam \in \Lam}~\FF_{\lam}$,
    where $\Lam$ is a filtered small category and
    each $\FF_{\lam}$ is a constructible $\Z/p\Z$-sheaf.
Replacing $\{\FF_{\lam}\}_{\lam \in \Lam}$ with
     their images into $\FF$ if necessary, we suppose that
        the transition maps are injective.
Since $\H^i(Y,\FF) \simeq
    \varinjlim_{\lam \in \Lam}~\H^i(Y,\FF_{\lam})$ and
$$
   \Hom(\H^i(Y,\FF),\Z/p\Z) \simeq
   \varprojlim_{\lam \in \Lam}~\Hom(\H^i(Y,\FF_{\lam}),\Z/p\Z),
$$
it is enough to show that
\stepcounter{equation}
\begin{equation}\label{eq:projlim}
\Ext^{N-i+1}_{Y,\Z/p\Z}(\FF,\nu_Y^N) \simeq
   \varprojlim_{\lam \in \Lam}~
   \Ext^{N-i+1}_{Y,\Z/p\Z}(\FF_{\lam},\nu_Y^N).
\end{equation}
Take an injective resolution $\nu_Y^N \ra I^{\bullet}$
    in the category of $\Z/p\Z$-sheaves on $Y_{\et}$.
The group $\Ext^m_{Y,\Z/p\Z}(\FF,\nu_Y^N)$ ($m \in \Z$)
    is the $m$-th cohomology group of the complex
      $\Hom_{Y}(\FF,I^{\bullet}) \simeq
        \varprojlim_{\lam \in \Lam}~\Hom_{Y}(\FF_{\lam},I^{\bullet})$.
Noting that $\Ext^m_{Y,\Z/p\Z}(\FF_{\lam},\nu_Y^N)$
    is finite for any $\lam \in \Lam$ and that the transition maps
      $\FF_{\lam} \ra \FF_{\lam'}$ ($\lam < \lam'$) are injective,
      we are reduced to the following standard fact on projective limits:
\par
\vspace{8pt}
\noindent
{\it Fact.
Let $\Lam$ be a cofiltered small category, and let
    $\{C_{\lam}^{\bullet}\}_{\lam \in \Lam}$
      be a projective system of complexes of abelian groups.
For $m \in \Z$ and $\lam \in \Lam$, put
    $H_{\lam}^m:=\H^m(C_{\lam}^{\bullet})$,
     the $m$-th cohomology group of $C_{\lam}^{\bullet}$.
Now fix $m \in \Z$, and
   assume that
    $\{C_{\lam}^{m-1} \}_{\lam \in \Lam}$,
    $\{ C_{\lam}^m \}_{\lam \in \Lam}$ and
    $\{ H_{\lam}^m \}_{\lam \in \Lam}$
   satisfy the Mittag-Leffler condition.
Then we have
   $\varprojlim_{\lam \in \Lam}\; H_{\lam}^{m+1} \simeq
      \H^{m+1}\big(\varprojlim_{\lam \in \Lam}\; C_{\lam}^{\bullet}\big)$.
}
\par
\vspace{8pt}
\noindent
This completes the proof of Proposition \ref{prop1-1}.
\end{pf}
We now turn to the proof of Theorem \ref{lem:trace}.
Without loss of generality, we may assume that
     $X$ is connected.
Then by Proposition \ref{prop1-1} for $i=N$,
    we have
\begin{align*}
    & \, \Hom_{D^b(Y_{\et},\Z/p\Z)}
    (U^{q}M^n \otimes U^{e'-q}M^{n'},\mu' \otimes\nu_Y^N[1] )\\
    \simeq & \, \Hom\big(\H^N(Y,U^{q}M^n \otimes U^{e'-q}M^{n'}),
     \mu \otimes \H^{N+1}(Y,\nu_Y^N)\big).
\end{align*}
Hence we are reduced to the equality of induced maps
     on cohomology groups
\begin{equation}\label{trace:eq}
\begin{CD}
\H^N(Y,\Theta^{q,n}) = \H^N(Y,\chi \circ f^{q,n}),
\end{CD}
\end{equation}
   where we wrote $\Theta^{q,n}$ for the composite morphism
$$
\begin{CD}
\hspace{-50pt}
\Theta^{q,n} : U^{q}M^n \otimes U^{e'-q}M^{n'} @>{{\mathrm{canonical}}}>>
      U^1M^n \otimes
      U^1M^{n'} @>{\Theta^n[N+2]}>> \mu' \otimes \nu_Y^N [1].
\hspace{-50pt}
\end{CD}
$$
\subsection{Reduction to higher local fields}\label{sect8.4}
In this step, \eqref{trace:eq} will be reduced to \eqref{trace:eq2} below.
We define {\it a chain on} $Y$ to be a sequence $(y_0,y_1,y_2,\dotsc,y_N)$
   of points (=spectra of fields) {\it over} $Y$
   such that $y_0$ is a closed point on $Y$
     and such that for each $m$ with $1 \leq m \leq N$,
       $y_m$ is a closed point on the scheme
$$
\begin{CD}
    \left[ \Spec  \left( \cdots
    \left( \left( \O_{Y,y_0}^{\h} \right) {}_{y_1}^{\h} \right)
    \cdots  \right) {}_{y_{m-1}}^{\h} \right]  \sm \{ y_{m-1} \},
\end{CD}
$$
where the superscript $\h$ means the henselization at
    the point given on subscript.
For a chain $(y_0,y_1,\cdots,y_N)$ on $Y$,
     each $\kappa(y_m)$ ($0 \leq m \leq N$)
      is an $m$-dimensional local field.
We write $\Chain(Y)$ for the set of chains on $Y$.
Now for ${\cal K} \in D^b(Y_{\et},\Z/p\Z)$
       and ${\frak d}=(y_0,y_1,\cdots,y_N) \in \Chain(Y)$,
   there is a composite map
$$
\begin{CD}
\H^0(y_N,{\cal K}) \ra
      \H^1_{y_{N-1}}(Y_{{\frak d}, N-1},{\cal K}) \ra
      \cdots
       \ra \H^{N-1}_{y_{1}}(Y_{{\frak d},1},{\cal K}) \ra
        \H^{N}_{y_{0}}(Y_{{\frak d}, 0},{\cal K}) \ra
   \H^{N}(Y,{\cal K}).
\end{CD}
$$
Here $Y_{{\frak d},m}$ ($0 \leq m \leq N$) denotes
     the henselian local scheme
$$
\begin{CD}
    \Spec \left( \cdots
    \left( \left( \O_{Y,y_0}^{\h} \right) {}_{y_1}^{\h} \right)
    \cdots  \right) {}_{y_{m}}^{\h}
\end{CD}
$$
and the map $\H^{N-m}_{y_m}(Y_{{\frak d},m},{\cal K}) \ra
    \H^{N-m+1}_{y_{m-1}}(Y_{{\frak d}, m-1},{\cal K})$
    ($1 \le m \le N$) is defined as the composite
$$
\begin{CD}
\H^{N-m}_{y_m}(Y_{{\frak d},m},{\cal K}) =
\H^{N-m}_{y_m}(Y_{{\frak d},m-1} \setminus \{ y_{m-1} \},{\cal K})
    @>>>
    \H^{N-m}(Y_{{\frak d},m-1} \setminus \{ y_{m-1} \},{\cal K}) \\
   @>{\delta^{\loc}({\cal K})}>>
\H^{N-m+1}_{y_{m-1}}(Y_{{\frak d}, m-1},{\cal K}).
\end{CD}
$$
Taking the direct sum with respect to all chains on $Y$,
    we obtain a map
$$
\begin{CD}
\delta_Y({\cal K}):
     \us{(y_0,y_1,\cdots,y_N) \in \Chain(Y)}\bigoplus
        \H^0(y_N,{\cal K})
           @>>> \H^N(Y,{\cal K}).
\end{CD}
$$
\begin{lem}\label{lem:trace:surj}
The map
$\delta_Y(U^{q}M^n \otimes U^{e'-q}M^{n'})$ is surjective.
\end{lem}
\begin{pf}
By Theorem \ref{thm:hyodo},
   the sheaf $U^qM^n \otimes U^{e'-q}M^{n'}$
       is a finitely successive extension
          of \'etale sheaves of the form $\FF \otimes \GG$,
     where $\FF$ and $\GG$ are locally free $(\O_Y)^p$-modules
        of finite rank and
          the tensor product is taken over $\Z/p\Z$.
We are reduced to the following sublemma.
\begin{sublem}\label{sublem:vanish}
Let $Z$ be a noetherian scheme which is of pure-dimension
   and essentially of finite type over $\Spec(k)$.
Put $d:=\dim(Z)$.
Let $\FF$ and $\GG$ be locally free $(\O_Z)^p$-modules
        of finite rank.
Then$:$
\begin{enumerate}
\item[(1)]
For any $x \in Z$ and $i > \codim_Z(x)$,
      $\H^i_x(Z,\FF \otimes \GG)$ is zero.
\item[(2)]
We have $\H^i(Z,\FF \otimes \GG)=0$ for $i > d$, and
the natural map
   $\bigoplus_{x \in Z^d}\ \H^d_x(Z,\FF \otimes \GG)
    \to \H^d(Z,\FF \otimes \GG)$ is surjective.
\item[(3)]
If $Z$ is henselian local, then
    $\H^i(Z,\FF \otimes \GG)$ is zero for $i > 0$.
\end{enumerate}
\end{sublem}
\medskip
\noindent
{\it Proof of Sublemma \ref{sublem:vanish}}.
Since the absolute Frobenius morphism $\F_Z:Z \to Z$ is finite
    by assumption,
    we have $\H^*(Z,\FF \otimes \GG) \simeq \H^*(Z,\F_{Z*}(\FF \otimes \GG))$
      and $\H^*_x(Z,\FF \otimes \GG) \simeq \H^*_x(Z,\F_{Z*}(\FF \otimes \GG))$
        for any $x \in Z$.
Hence we are reduced to the case where $\FF$ and $\GG$ are
    locally free $\O_Z$-modules of finite rank.

We first show (3).
Let $R$ be the affine ring of $Z$,
   which is a henselian local ring by assumption.
Let $R^{\sh}$ be the strict henselization of $R$.
Without loss of generality, we may assume that
    $\FF=\GG=\O_Z$.
By the isomorphism $\H^q(Z,\O_Z \otimes \O_Z)
    \simeq \H^q_{\Gal}(G_R,R^{\sh} \otimes R^{\sh})$
      with $G_R:=\Gal(R^{\sh}/R)$,
      our task is to show that the right hand side is zero
        for $q > 0$.
We show that for a finite \'etale galois extension $R'/R$
    with Galois group $G:=\Gal(R'/R)$,
   we have $\H^q(G,R' \otimes R')=0$ for $q > 0$.
Indeed, by taking a normal basis,
   we have $R' \simeq R[G]$ as left $R[G]$-modules, and
    there is an isomorphism of left $G$-modules
$$
\begin{CD}
R[G] \otimes R[G] @>{\simeq}>> R[G] \otimes (R[G]^{\circ}), \quad
    a[g] \otimes b[h] \mapsto
    a[g] \otimes b[ g^{-1}h],
\end{CD}
$$
where $a$ and $b$ (resp.\ $g$ and $h$)
   are elements of $R$ (resp.\ of $G$), and
     $R[G]^{\circ}$ denotes the abelian group $R[G]$ with trivial
     $G$-action.
Hence $R' \otimes R'$ is an induced $G$-module in the sense of
    \cite{se}, I.2.5 and we obtain the assertion.

We next prove (1) and (2) by induction on $d$ and a standard
    local-global argument (cf.\ \cite{ras}, 1.22).
The case $d=0$ follows from (3).
Assume $d \ge 1$ and
   that (1) and (2) hold true for schemes of dimension $\le d-1$.
We first show (1).
Indeed, the case $\codim_Z(x)=0$ follows from the case $d=0$.
If $\codim_Z(x) \ge 1$ and $i \ge 1$,
   then the connecting homomorphism
$$
\begin{CD}
\delta^{\loc}({\cal K}) :
\H^{i-1}(\Spec(\O_{Z,x}^{\h})\setminus \{ x \},\FF \otimes \GG) @>>>
    \H^i_x(\Spec(\O_{Z,x}^{\h}),\FF \otimes \GG)
     = \H^i_x(Z,\FF \otimes \GG)
\end{CD}
$$
is surjective by (3), and
    the left hand side is zero for $i > \codim_Z(x)$
     by the induction hypothesis.
Thus we obtain (1).
Finally one can easily check (2) by (1) and
   a local-global spectral sequence
$$
E_1^{u,v}=\bigoplus{}_{x \in Z^u} \ \H^{u+v}_x(Z,\FF \otimes \GG)
   \Lra \H^{u+v}(Z,\FF \otimes \GG).
$$
This completes the proof of Sublemma \ref{sublem:vanish}
   and Lemma \ref{lem:trace:surj}.
\end{pf}
By Lemma \ref{lem:trace:surj},
   \eqref{trace:eq} is reduced to
     the formula
\addtocounter{equation}{2}
\begin{equation}\label{trace:eq2}
\begin{CD}
\H^0(y_N ,\Theta^{q,n}) = \H^0(y_N,\chi \circ f^{q,n})
\end{CD}
\end{equation}
for all chains $(y_0,y_1,\cdots,y_N) \in \Chain(Y)$,
         which will be proved in \S\ref{sect8.6} below.
\subsection{Galois descent by corestriction maps}\label{sect8.5}
We prove here the following lemma:
\begin{lem}\label{lem:H^0}
Let $F$ be a field of characteristic $p>0$.
Let $V_1$ and $V_2$ be discrete $G_F$-$\Z/p\Z$-modules
      which are finitely successive extensions of
            finite direct sums of copies of
               $\ol F$ as $G_F$-modules.
Then $(V_1 \otimes V_2)^{G_F}$ agrees with
$$
\begin{CD}
\us{F'/F:~ \text{finite galois}}{\cup}
    \Image \big( (V_1)^{G_{F'}} \otimes (V_2)^{G_{F'}}
                \os{\subset}{\lra} (V_1 \otimes V_2)^{G_{F'}}
                 @>{\cores_{F'/F}}>> (V_1 \otimes V_2)^{G_F} \big),
\end{CD}
$$
where all tensor products are taken over $\Z/p\Z$,
     and $F'$ runs through all finite galois
       field extensions of $F$ contained in $\ol F$.
\end{lem}
\begin{pf}
It suffices to show the case
     $V_1 = V_2 = \ol F$.
We prove that
       the corestriction map
$$
\begin{CD}
   \cores_{F'/F}: F' \otimes F'
                 @>>> (F' \otimes F')^{G}, \quad
                   x \otimes y \mapsto
   \Sigma_{g \in G}~ gx \otimes gy
\end{CD}
$$
is surjective for a finite galois extension
          $F'/F$ with $G:=\Gal(F'/F)$,
         which implies the assertion
             by a limit argument.
Since $F' \simeq F[G]$ as $F[G]$-modules,
    we have
$$
\begin{CD}
(F' \otimes F')^G \simeq (F \otimes F) \otimes
    (\Z/p\Z[G] \otimes \Z/p\Z[G])^G
\end{CD}
$$
by the finiteness of $G$ and
   the flatness of $\Z/p\Z$-modules over $\Z/p\Z$.
Hence the surjectivity of $\cores_{F'/F}$
   follows from that of the map
$$
\begin{CD}
\Z/p\Z[G] \otimes \Z/p\Z[G] @>>>
     (\Z/p\Z[G] \otimes \Z/p\Z[G])^{G}, \quad
                   x \otimes y \mapsto
   \Sigma_{g \in G}~ gx \otimes gy.
\end{CD}
$$
Thus we obtain the lemma.
\end{pf}
\subsection{Proof of $\eqref{trace:eq2}$}\label{sect8.6}
In this step,
    we finish the proof of Theorem \ref{lem:trace}.
Fix an arbitrary chain $(y_0,y_1,\cdots,y_N) \in \Chain(Y)$.
Put $F_N:=\kappa(y_N)$ and
$$
\begin{CD}
L_{N+1} := \Frac \left[ \left( \cdots
    \left( \left(
\O_{X,y_0}^{\h} \right) {}_{y_1}^{\h} \right)
    \cdots  \right)
{}_{y_{N}}^{\h} \right],
\end{CD}
$$
where
    $L_{N+1}$ is a henselian discrete valuation field
       (of characteristic $0$)
        with residue field $F_N$, that is,
          $L_{N+1}$ is an $(N+1)$-dimensional local field.
Now let $F/F_N$ be a finite separable field extension.
Put $y:=\Spec(F)$ and
$$
\begin{CD}
A^{q,n}(F) : = \H^0(y,U^qM^n) \otimes
   \H^0(y,U^{e'-q}M^{n'})
         \subset
      \H^0(y,U^qM^n \otimes U^{e'-q}M^{n'}).
\end{CD}
$$
By Lemma \ref{lem:H^0} (for the subfield $(F_N)^p \subset F_N$)
     and the naturality of corestriction maps,
    the formula \eqref{trace:eq2} for $y_N$ is reduced to
        the formula
\begin{equation}\label{galois:cor4}
\hspace{-70pt}
\begin{CD}
   \H^0(y,\Theta^{q,n}) \vert_{A^{q,n}(F)}
      =
       \H^0(y,\chi \circ f^{q,n})
\vert_{A^{q,n}(F)}.
\end{CD}
\hspace{-50pt}
\end{equation}
To prove this equality,
   we compute the left hand side, i.e., the composite map
\begin{equation}\label{galois:cor5}
\hspace{-50pt}
\begin{CD}
A^{q,n}(F)
    \; \hra \;
      \H^0(y,U^qM^n \otimes U^{e'-q}M^{n'})
    @>{\H^0(y,\Theta^{q,n})}>>
       \mu' \otimes
   \H^1(y,\Omega_{y,\log}^N).
\end{CD}
\hspace{-50pt}
\end{equation}
Let $L/L_{N+1}$ be the finite unramified extension
     corresponding to $F/F_N$.
For $i>0$, put $k^M_i(L):=K^M_i(L)/pK^M_i(L)$.
By a similar argument as for \cite{bk}, (5.15), we have
$$
\begin{CD}
A^{q,n}(F) = \{U^qk^M_n(L)/U^{e'}k^M_n(L) \}
    \otimes \{U^{e'-q}k^M_{n'}(L)/U^{e'}k^M_{n'}(L) \}.
\end{CD}
$$
Let us recall that $1 \leq q \leq e'-1$ by assumption.
In view of the construction of $\Theta^{n}$ (cf.\ \S\ref{sect8.1})
    and the fact that $U^{e'+1}k^M_{N+2}(L)=0$
     (\cite{bk}, (5.1.i)),
     the map \eqref{galois:cor5} is written by
     the product of Milnor $K$-groups and
      boundary maps of Galois cohomology groups:
{\allowdisplaybreaks
\begin{equation}\label{galois:cor7}
\hspace{-50pt}
\begin{CD}
A^{q,n}(F)
         @>{{\mathrm{product}}}>> U^{e'}k^M_{N+2}(L)
         @>{{\mathrm{Galois~symbol}}}>>
   \H^{N+2}(L,\mu_p^{\otimes N+2}) \\
        @>{\simeq}>>
           \mu \otimes
   \H^{1}(F,
              \H^{N+1}(L^{\ur},\mu_p^{\otimes N+1}))
    @>{\id \otimes \eqref{isom:bk}}>>
     \mu \otimes \H^1(y,\Omega_{y,\log}^N),
\end{CD}
\hspace{-50pt}
\end{equation}
}\par
\noindent
where $L^{\ur}$ denotes the maximal unramified extension of
     $L$, and the third arrow
      is obtained from a Hochschild-Serre
         spectral sequence together with
          the facts that $\cd_p(F)=1$ and
$\cd_p(L^{\ur})=N+1$
           (cf.\ Lemma \ref{lem4''}).
Here we compute the product of symbols:
\addtocounter{thm}{3}
\begin{lem}\label{lem:prod}
For $\alpha_1$, $\alpha_2 \in O_L \sm \{ 0 \}$
     and $\beta_1,\dotsc,\beta_N \in L^{\times}$, we have
$$
\begin{CD}
\{1+\pi^q\alpha_1,\beta_1,\dotsc,\beta_{n-1} ,
    1+\pi^{e'-q}\alpha_2,
      \beta_{n},\dotsb,\beta_{N} \} @.~ =  ~@.
   (-1)^{N+n} \cdot q \cdot \{ 1+\pi^{e'}\alpha_1\alpha_2,
   \beta_1,\dotsc,\beta_{N},\pi\} \\
      @.  @.+ (-1)^n \cdot \{1+\pi^{e'}\alpha_1\alpha_2,
   \alpha_1,\beta_1,\dotsc,\beta_{N}\}
\end{CD}
$$
   in $k^M_{N+2}(L)$.
The second term on the right hand side is zero if
     $\beta_i$ belongs to
       $O_L^{\times}$ for all $i$.
\end{lem}
\begin{pf}
We compute the symbol
    $\{1+\pi^q \alpha_1,1+\pi^{e'-q} \alpha_2 \} \in k^M_2(L)$:
{\allowdisplaybreaks
\begin{align*}
\{1+\pi^q \alpha_1,1+\pi^{e'-q} \alpha_2 \}
      = ~&
        \{1+ \pi^q \alpha_1+ \pi^{e'} \alpha_1 \alpha_2,
          1+\pi^{e'-q} \alpha_2\}\\
        & - \{(1+ \pi^q \alpha_1+ \pi^{e'} \alpha_1 \alpha_2)
            (1+ \pi^q \alpha_1)^{-1},
          1+\pi^{e'-q} \alpha_2\} \\
    \os{(1)}{=} ~&
        - \{1+ \pi^q \alpha_1+ \pi^{e'} \alpha_1 \alpha_2,
          -\pi^q \alpha_1 \} \\
    \os{(2)}{=} ~&
          - \{1+ \pi^{e'} \alpha_1 \alpha_2 (1+\pi^q \alpha_1)^{-1},
          -\pi^q \alpha_1 \} \\
    \os{(3)}{=} ~&
            - \{1+ \pi^{e'} \alpha_1 \alpha_2,
          \pi^q \alpha_1 \}.
\end{align*}
\par
}
\noindent
$(1)$ follows from the equality
$\{1+x_1x_2,x_1\}=-\{1+x_1x_2,-x_2\}$
        (applied to the first term)
   and the fact that the second term is contained in $U^{e'+1}k^M_2(L)=0$
   (\cite{bk}, (4.1), (5.1.i)).
$(2)$ follows from the equality $\{1+x,-x\}=0$,
and $(3)$ follows from loc.\ cit., (4.3).
The equality assertion in the lemma follows from this computation.
The last assertion follows from loc.\ cit., (4.3) and
       the fact that $\Omega_F^{N+1}=0$.
\end{pf}
To calculate the last two maps in \eqref{galois:cor7},
    we need the following lemma,
       which is a kind of explicit formula for $L$
      (see \S\ref{sect8.2} for $u$ and $\chi$):
\begin{lem}\label{lem:BK}
The following square commutes$:$
\addtocounter{equation}{2}
\begin{equation}\label{galois:cor6}
\begin{CD}
\Omega_F^N @>{\omega \mapsto u \otimes_k \omega}>>
         (\mu \otimes k) \otimes_k (\Omega_F^N/\BB_F^{N}) \\
@V{\rho^{e'}}VV   @V{\H^0(y,\chi)}VV \\
\H^{N+2}(L,\mu_p^{\otimes N+2}) @>{\simeq}>>
    \mu \otimes \H^1(y, \Omega_{y,\log}^N),
\end{CD}
\end{equation}
where the bottom arrow is the composite of the last two maps
       in $\eqref{galois:cor7}$ and
   $\rho^{e'}$ denotes the Bloch-Kato map
      sending $\alpha \cdot \dlog(\beta_1) \wedge \dotsb \wedge \dlog (\beta_N)$
         $(\alpha \in F$, $\beta_i \in F^{\times})$
         to $\{1+\pi^{e'}\wt{\alpha},\wt{\beta_1},\dotsc,\wt{\beta_N},\pi \}$
   $(\wt{\alpha}$ and $\wt{\beta_i}$'s
      are lifts of $\alpha$ and $\beta_i$'s,
        respectively$)$.
\end{lem}
\begin{pf}
This commutativity would be well-known to experts
    (cf.\ \cite{kk:ex} for the case $p > N+3$,
     see also \cite{fv}, VII.4).
However we include here a simple proof using a classical argument
      originally due to Hasse \cite{Has}
      to verify the above commutativity including signs.
By \cite{kk:cft2}, p.\ 612, Lemma 2,
   the bottom horizontal arrow of the diagram \eqref{galois:cor6} maps
$$
   \zeta \cup {\mathrm{Inf}}(x) \cup
    \{\wt{\beta_1},\wt{\beta_2},\dotsc,\wt{\beta_N},\pi\}
      \mapsto \zeta \otimes (-x) \cup
     (\dlog(\beta_1) \wedge \dotsb \wedge \dlog(\beta_N))
$$
for $\zeta \in \mu$, $x \in \H^1(F,\Z/p\Z)$
       and $\beta_i \in F^{\times}$.
Hence it is enough to show the following:
\par
\vspace{8pt}
\noindent
{\it Claim.
Fix a primitive $p$-th root of unity
     $\zeta_p \in \mu$,
     and consider the composite map
$$
\begin{CD}
F \lra \H^1(F, \Z/p\Z)
      \os{{\mathrm{Inf}}}{\lra} \H^1(L,\Z/p\Z)
      @>{1 \mapsto \zeta_p}>> \H^1(L,\mu_p)
      \os{\simeq}{\lra} L^{\times}/(L^{\times})^p,
\end{CD}
$$
where the first map is the boundary map
     of Artin-Schreier theory and
       the last isomorphism is the inverse of
          the boundary map of Kummer theory.
Then this composite map sends
    $- v(\zeta_p)^{-p} \alpha \in F$ to
         $1 + \pi^{e'} \wt{\alpha} \mod (L^{\times})^p$,
          where $\wt{\alpha}$ denotes a lift of $\alpha$ to $O_L$
       $($note that $U^{e'+1}L^{\times} \subset (L^{\times})^p)$.
See Definition $\ref{def:fqn}$ for the definition of $v(\zeta_p)$.}
\par
\vspace{8pt}
\noindent
{\it Proof of Claim.}
It suffices to show that
      $\alpha \in F$ maps to
       $1 - (1-\zeta_p)^p \, \wt{\alpha} \mod (L^{\times})^p$.
Consider the following equations in $T$ over $F$
      and $L$, respectively:
\begin{align}
         T^p - T &= \alpha,
\label{eq1} \\
        T^p &=1 - (1-\zeta_p)^p \cdot \wt{\alpha}.
\label{eq2}
\end{align}
We show that
           the Artin-Schreier character $G_F \ra \Z/p\Z$
              associated with \eqref{eq1}
    induces the Kummer character $G_L \ra \mu$
              associated with \eqref{eq2}
    by the composition $G_L \ra G_F \ra \Z/p\Z \ra \mu$.
Let $\beta \in \ol L$ be a solution to \eqref{eq2}.
By the congruity relation
$$
     (-1)^p \cdot p \equiv (1- \zeta_p)^{p-1} \mod \pi^{e'}O_L,
     \phantom{'}
$$
    one can easily show that $\beta$ is contained in $O_L^{\ur}$ and that
$$
      \gamma :=  (1-\beta)/(1-\zeta_p) \mod \pi O_{L}^{\ur} \in  \ol F
$$
is a solution to \eqref{eq1}.
Moreover, $\sigma \in G_L$
     satisfies $\sigma(\beta)/\beta=\zeta_p^m \in \mu$
         if and only if $\sigma(\gamma)-\gamma=m \in \Z/p\Z$,
           where $G_L$ acts on $\ol F$
              via the canonical projection $G_L \ra G_F$.
Thus we obtain the claim and Lemma \ref{lem:BK}.
\end{pf}
We now turn to the proof of \eqref{galois:cor4}.
Let $\alpha_1$, $\alpha_2 \in O_L \sm \{ 0 \}$,
   and $\beta_i \in O_L^{\times} \cup \{ \pi \}$
        ($1 \leq i \leq N$).
By Lemmas \ref{lem:prod} and \ref{lem:BK},
      the value of the symbol
$$
\begin{CD}
\{1+\pi^q\alpha_1,\beta_1,\dotsb,\beta_{n-1}\}
   \otimes \{ 1+\pi^{e'-q}\alpha_2,
      \beta_n,\dotsb,\beta_{N} \}
        \in U^qk^M_{n}(L) \otimes U^{e'-q}k^M_{n'}(L)
\end{CD}
$$
under \eqref{galois:cor5} agrees with
   the value of the following element of $\mu \otimes (\Omega_F^N/\BB_F^{N})$
      under $\H^0(y,\chi)$:
$$
\begin{CD}
\begin{cases}
   u \otimes_k (-1)^{n+N}
      [ q \cdot \ol {\alpha_1 \alpha_2} \cdot
     (\wedge_{1 \leq i \leq N}~\dlog\ol {\beta_i}) \mod \BB_F^{N}]
   \quad & (\hbox{if } \beta_i \in O_L^{\times}
           \hbox{ for all } i), \\
   u \otimes_k (-1)^{n+N+i'}
      [\ol {\alpha_2} \cdot  d \ol {\alpha_1}
   \wedge (\wedge_{1 \leq i \leq N,i \not= i'}~
           \dlog \ol {\beta_i}) \mod \BB_F^{N}]
            & (\hbox{if } \beta_{i}=\pi
           \hbox{ for exactly one } i=:i'), \\
   0       & (\hbox{otherwise}),
\end{cases}
\end{CD}
$$
where for $x \in O_L$ (resp.\ $x \in O_L^{\times}$),
     $\ol x$ denotes its residue class in $F$
       (resp.\ in $F^{\times}$).
Thus comparing this presentation of \eqref{galois:cor5}
   with the definition of $f^{q,n}$
           (cf.\ \S\ref{sect8.2}),
         we conclude that the equality \eqref{galois:cor4} holds.
This completes the proof of Theorem \ref{lem:trace}.
\qed
\medskip
\section{\bf Duality of $p$-adic vanishing cycles}\label{sect9}
\medskip
In this section we prove Theorem \ref{prop:sdual2} below,
    which will be used in \S\ref{sect10}.
\subsection{Statement of the result}\label{sect9.1}
Let the notation be as in \S\ref{sect8.0}.
   We prove the following:
\begin{thm}\label{prop:sdual2}
Let $n$ be $1 \leq n \leq N+1$ and
   put $n':=N+2-n$.
Assume that $X$ is proper over $\Spec(O_K)$.
Then for an integer $i$,
   the pairing induced by $\Theta^n$ and $\tr_Y$
     $($cf.\ Theorem $\ref{thm:milne})$
\stepcounter{equation}
\begin{equation}\label{dual:pair3}
\hspace{-30pt}
\begin{CD}
a^i : \H^i(Y,U^1M^n) \times
    \H^{N-i}(Y,U^1M^{n'})
   \;\os{\Theta^n}{\lra}\; \mu \otimes \H^{N+1}(Y,\nu_Y^N)
           @>{\id \otimes \tr_Y}>> \mu
\end{CD}
\hspace{-30pt}
\end{equation}
   is a non-degenerate
       pairing of finite $\Z/p\Z$-modules.
\end{thm}
\noindent
To prove this theorem,
     we first calculate the map $f^{q,n}$ defined
       in \S\ref{sect8.2} (cf.\ Lemma \ref{sublem3.5} below).
Let $U^{\bullet}M^n$ and $V^{\bullet}M^n$ be as in
    Definition \ref{def:vcyc}.
We further define the subsheaf
    $T^qM^n \subset U^qM^n$ $(q \geq 1)$
      as the part generated by $V^qM^n$ and symbols of the form
$$
\{1+\pi^q \alpha^p, \beta_1, \cdots, \beta_{n-1} \}
$$
with $\alpha \in \iota^*\O_X$ and each
     $\beta_i \in \iota^*j_*\O^{\times}_{X_K}$.
By definition we have
$$
\begin{CD}
U^{q+1}M^n \subset V^qM^n \subset T^qM^n \subset U^qM^n.
\end{CD}
$$
For $q \geq 1$, we put
$$
\begin{CD}
\gr_{U/T}^qM^n:=U^qM^n/T^qM^n, \qquad
\gr_{T/V}^qM^n:=T^qM^n/V^qM^n \\
    \hbox{ and }\quad
\gr_{V/U}^qM^n:=V^qM^n/U^{q+1}M^n.
\end{CD}
$$
Let us recall that $e'=pe/(p-1)$ is an integer divided by $p$
       (because $\zeta_p \in K$).
By Theorem \ref{thm:hyodo} (3), (4),
      the sheaf $U^{e'}M^n$ is zero, and
        for $q$ with $1 \leq q \leq e'-1$ we have isomorphisms
\begin{equation}\label{dual:isom:bkh}
\begin{CD}
\rho_1^{q,n}:
    \gr_{U/T}^qM^n @>{\simeq}>> \hspace{-53.5pt} \omega_Y^{n-1}/ \ZZ_Y^{n-1},\\
\rho_2^{q,n}: \gr_{T/V}^qM^n @>{\simeq}>>
       \begin{cases}
           \ZZ_Y^{n-1}/\BB_Y^{n-1}  \qquad &(p \hspace{-4pt}\not\vert q),\\
           0                        \qquad &(p \vert q),
       \end{cases}\\
\rho_3^{q,n}:
    \gr_{V/U}^qM^n @>{\simeq}>> \hspace{-56pt} \omega_Y^{n-2}/ \ZZ_Y^{n-2},
\end{CD}
\end{equation}
given by the following, respectively:
$$
\begin{CD}
\rho_1^{q,n}: ~@.
    \{1+\pi^q \alpha, \beta_1, \cdots, \beta_{n-1} \} \mod T^qM^n
    ~@. \mapsto~ @.\hspace{-47pt}
     \ol{\alpha} \cdot
          (\wedge_{i=1}^{n-1} ~\dlog \ol{\beta_i}) \mod \ZZ_Y^{n-1},\\
\rho_2^{q,n}: ~@.
\{1+\pi^q \alpha^p, \beta_1, \cdots, \beta_{n-1} \} \mod V^qM^n
   ~@. \mapsto~ @.  \ol{\alpha}{}^p \cdot
   (\wedge_{i=1}^{n-1} ~\dlog \ol{\beta_i}) \mod \BB_Y^{n-1},
            \quad (p \hspace{-4pt}\not \vert q),\\
\rho_3^{q,n}: ~@.
\{1+\pi^q \alpha, \beta_1 \cdots, \beta_{n-2}, \pi \} \mod U^{q+1}M^n
   ~@. \mapsto~ @.\hspace{-49pt}
     \ol{\alpha} \cdot
          (\wedge_{i=1}^{n-2} ~\dlog \ol{\beta_i}) \mod \ZZ_Y^{n-2}.
\end{CD}
$$
Here $\alpha$ (resp.\ each $\beta_i$)
   denotes a local section of $\iota^*\O_X$
       (resp.\ $\iota^*j_*\O^{\times}_{X_K}$), and
       $\ol {\alpha}$ (resp.\ $\ol{\beta_i}$)
        denotes its residue class in $\O_Y$
       (resp.\ in ${\cal L}_Y^{\gp}$ under \eqref{log:surj}).
The following lemma
    follows from straight-forward computations
      on symbols, whose proof is left to the reader
       (cf.\ Remark \ref{rem:log} (2), Definition \ref{def:fqn}):
\addtocounter{thm}{2}
\begin{lem}\label{sublem3.5}
Let $n$ and $n'$ be as in Theorem $\ref{prop:sdual2}$, and
assume $1 \leq q \leq e'-1$. Then$:$
\begin{enumerate}
\item[(1)]
$f^{q,n}$ annihilates
   the subsheaf of $U^{q}M^n \otimes U^{e'-q}M^{n'}$
     generated by
$U^{q+1}M^n \otimes U^{e'-q}M^{n'}$,
$U^qM^n \otimes U^{e'-q+1}M^{n'}$,
$V^{q}M^n \otimes T^{e'-q}M^{n'}$ and
$T^{q}M^n \otimes V^{e'-q}M^{n'}$.
\item[(2)]
The composite map
$$
\begin{CD}
\omega_Y^{n-1}/\ZZ_Y^{n-1}
    \otimes
\omega_Y^{n'-2}/\ZZ_Y^{n'-2}
      @>{(\rho_1^{q,n} \otimes
\rho_3^{e'-q,n'})^{-1}}>>
   \gr_{U/T}^qM^n \otimes \gr_{V/U}^{e'-q}M^{n'}
   \os{f^{q,n}}{\lra}  \mu' \otimes (\omega_Y^N/\BB_Y^N)
\end{CD}
$$
sends a local section $x \otimes y$ to
   $u \otimes_k (-1)^n \cdot (dx) \wedge y$.
Similarly,
   the composite map
$$
\begin{CD}
\omega_Y^{n-2}/\ZZ_Y^{n-2}
   \otimes \omega_Y^{n'-1}/\ZZ_Y^{n'-1}
@>{(\rho_3^{q,n} \otimes
\rho_1^{e'-q,n'})^{-1}}>>
    \gr_{V/U}^qM^n \otimes \gr_{U/T}^{e'-q}M^{n'}
   \os{f^{q,n}}{\lra}  \mu' \otimes (\omega_Y^N/\BB_Y^N)
\end{CD}
$$
sends a local section $x \otimes y$ to
   $u \otimes_k (-1)^{N} \cdot x \wedge dy$.
\item[(3)]
If $q$ is prime to $p$,
    then the composite map
$$
\begin{CD}
\ZZ_{Y}^{n-1}/\BB_{Y}^{n-1}
     \otimes \ZZ_{Y}^{n'-1}/\BB_{Y}^{n'-1}
   @>{(\rho_2^{q,n} \otimes \rho_2^{e'-q,n'})^{-1}}>>
    \gr_{T/V}^qM^n \otimes \gr_{T/V}^{e'-q}M^{n'}
   \os{f^{q,n}}{\lra}  \mu' \otimes (\omega_Y^N/\BB_Y^N)
\end{CD}
$$
sends a local section $x \otimes y$ to
   $u \otimes_k (-1)^{N+n} q \cdot x \wedge y$.
\end{enumerate}
\end{lem}
\medskip
\noindent
\subsection{Proof of Theorem \ref{prop:sdual2}}\label{sect9.2}
In this subsection, we reduce the theorem to Lemma \ref{sublem6} below.
The finiteness of the groups in the pairing
   \eqref{dual:pair3} follows from
         the finiteness of $k$, the properness of $Y$ and the fact that
             the sheaves $U^1M^n$ and $U^1M^{n'}$
               are finitely successive extensions of
   coherent $(\O_Y)^p$-modules
                    (cf.\ \eqref{dual:isom:bkh}).
To show the non-degeneracy of \eqref{dual:pair3},
   we introduce an auxiliary descending filtration
      $Z^r M^n$ ($r \geq 1$) on $U^1M^n$ defined as
$$
\begin{CD}
   Z^{r} M^n :=
    \begin{cases}
   U^{q}M^n   \qquad (\hbox{if }r\equiv 1 \mod 3
   \hbox{ and } q =(r+2)/3),\\
   T^{q}M^n   \qquad (\hbox{if }r\equiv 2 \mod 3
   \hbox{ and } q =(r+1)/3),\\
   V^{q}M^n   \qquad (\hbox{if }r\equiv 0 \mod 3
   \hbox{ and } q =r/3).
    \end{cases}
\end{CD}
$$
Note that $Z^1M^n=U^1M^n$
   and $Z^{r}M^n=0$ for $r \geq 3e'-2$.
We first show
\begin{lem}\label{sublem4}
Assume $1 \leq r \leq 3e'-3$.
Then$:$
\begin{enumerate}
\item[(1)]
The composite map
$$
\begin{CD}
   \H^N(Y,U^1M^n \otimes Z^{3e'-2-r}M^{n'})
       @>>> \H^N(Y,U^1M^n \otimes U^1M^{n'}) \\
       @>{\H^{2N+3}(Y,\Theta^n)}>> \mu \otimes \H^{N+1}(Y,\nu_Y^N)
       @>{\id \otimes \tr_Y}>{\simeq}> \mu
\end{CD}
$$
     induces a map
\stepcounter{equation}
\begin{equation}\label{dual:factor1}
\begin{CD}
   \H^N(Y,(U^1M^n/Z^{r+1}M^n) \otimes Z^{3e'-2-r}M^{n'})
       @>>> \mu.
\end{CD}
\end{equation}
\item[(2)]
The composite map
$$
\begin{CD}
   \H^N(Y,\gr_Z^{r}M^n \otimes Z^{3e'-2-r}M^{n'}) @>>>
   \H^N(Y,(U^1M^n/Z^{r+1}M^n) \otimes Z^{3e'-2-r}M^{n'})\\
   @>{\eqref{dual:factor1}}>>
       \hspace{-195pt} \mu
\end{CD}
$$
    induces a map
\begin{equation}\label{dual:factor2}
\begin{CD}
   \H^N(Y,\gr_Z^{r}M^n \otimes \gr_Z^{3e'-2-r}M^{n'})
       @>>> \mu.
\end{CD}
\end{equation}
\item[(3)]
   We put
$$
\begin{CD}
\FF^{r,n} ~@. :=  @.~ (U^1M^n/Z^{r+1}M^n) \otimes Z^{3e'-2-r}M^{n'}, \\
{\cal G}^{r,n} ~@. := @.~ \FF^{r,n} \left/
    \left( \gr_Z^{r}M^n \otimes Z^{3e'-1-r}M^{n'} \right) \right., \\
{\cal H}^{r,n} ~@. := @.~ \gr_Z^{r}M^n \otimes \gr_Z^{3e'-2-r}M^{n'}
    \phantom{AAAAAA}
\end{CD}
$$
$($note that $\gr_Z^{r}M^n \otimes Z^{3e'-1-r}M^{n'}$
      is a subsheaf of $\FF^{r,n}$, because
         a $\Z/p\Z$-sheaf is flat over $\Z/p\Z)$.
Then the map \eqref{dual:factor1}
      induces a map
\begin{equation}\label{dual:factor3}
\begin{CD}
   \H^N(Y,{\cal G}^{r,n})
       @>>> \mu.
\end{CD}
\end{equation}
If $r \geq 2$,
     then this map makes the following diagram commutative$:$
\begin{equation}\label{dual:factor4}
\begin{CD}
   \phantom{ia}
    \H^N(Y,{\cal H}^{r,n}) \oplus \H^N(Y, \FF^{r-1,n})
       @>>> \H^N(Y,{\cal G}^{r,n}) \\
    @V{\eqref{dual:factor2} \oplus}V
   {\eqref{dual:factor1}~{\mathrm{for}}~r-1}V  @VV{\eqref{dual:factor3}}V \\
    \mu \oplus \mu @>{{\mathrm{product}}}>> \mu,
\end{CD}
\end{equation}
where the top horizontal arrow is induced by
    a natural inclusion
      ${\cal H}^{r,n} \oplus \FF^{r-1,n}
   \subset {\cal G}^{r,n}$.
\end{enumerate}
\end{lem}
\begin{pf*}{\it Proof of Lemma \ref{sublem4}}
We prove only (1).
(2) and (3) are similar and left to the reader.
We use the notation we fixed in (3).
Let $q$ be the maximal integer with $3(q-1)<r$.
Noting that
$Z^{r}M^n \otimes Z^{3e'-2-r}M^{n'}
   \subset U^{q}M^n \otimes U^{e'-q}M^{n'}$,
     consider the composite map
$$
\begin{CD}
\H^N(Y,Z^{r}M^n \otimes Z^{3e'-2-r}M^{n'})
   @>>> \H^N(Y,U^1M^n \otimes Z^{3e'-2-r}M^{n'})
   @>{(*)}>> \mu,
\end{CD}
$$
where the arrow $(*)$ denotes the first composite map in (1).
By Theorem \ref{lem:trace} (cf.\ \eqref{trace:eq}),
     this composite map agrees with
       that induced by $\chi \circ f^{q,n}$.
By Lemma \ref{sublem3.5} (1),
     $f^{q,n}$ annihilates the subsheaf
$Z^{r+1}M^n \otimes Z^{3e'-2-r}M^{n'}$
    of $Z^{r}M^n \otimes Z^{3e'-2-r}M^{n'}$.
Hence the arrow $(*)$ induces a map of the form \eqref{dual:factor1} by
      the short exact sequence of sheaves
$$
\begin{CD}
0 \lra Z^{r+1}M^n \otimes Z^{3e'-2-r}M^{n'}
    \lra U^1M^n \otimes Z^{3e'-2-r}M^{n'}
    \lra \FF^{r,n} \lra 0
\end{CD}
$$
    and Sublemma \ref{sublem:vanish} (2) (cf.\ \eqref{dual:isom:bkh}).
Thus we obtain the lemma.
\end{pf*}
Now we turn to the proof of Theorem \ref{prop:sdual2}.
By the trace maps \eqref{dual:factor1} and \eqref{dual:factor2},
     there are induced pairings
\begin{equation}\label{dual:pair:4}
\begin{CD}
b^{i,r} @.~:~ @. \H^i(Y,U^1M^n/Z^{r+1}M^n) ~@. \times @.~
   \H^{N-i}(Y,Z^{3e'-2-r}M^{n'})
           @. @>>> @. \mu,\\
c^{i,r} @.~:~ @. \H^i(Y,\gr_Z^{r}M^n) ~ @. \times @.~
   \H^{N-i}(Y,\gr_Z^{3e'-2-r}M^{n'})
           @. @>>> @. \mu,
\end{CD}
\end{equation}
   for $i$ and $r$ with $1 \leq r \leq 3e'-3$.
Note that $b^{i,3e'-3}=a^i$ and $b^{i,1}=c^{i,1}$.
By the commutative diagram \eqref{dual:factor4},
    there is a commutative diagram with exact rows for $r \geq 2$
     (after changing the signs of $(\flat)$ suitably)
$$
{\scriptsize
\begin{CD}
   \H^{i-1}(Y,U^1M^{n}/Z^{r}M^n) @.~\ra~@.
   \H^{i}(Y,\gr_Z^{r}M^n) @.~\ra~@.
   \H^{i}(Y,U^1M^{n}/Z^{r+1}M^n) @.~\ra~@.
   \H^{i}(Y,U^1M^{n}/Z^{r}M^n)   @.~\ra~@.
    \H^{i+1}(Y,\gr_Z^{r}M^n)\\
   @V{b^{i-1,r-1}}VV @. @V{c^{i,r}}VV @. @V{b^{i,r}}VV @. @V{b^{i,r-1}}VV
   @. @V{c^{i+1,r}}VV\\
   \H^{\ell+1}(Y,Z^{t+1}M^{n'})^* @.~\os{(\flat)}\ra~@.
   \H^{\ell}(Y,\gr_Z^{t}M^{n'})^*         @.~\ra~@.
   \H^{\ell}(Y,Z^{t}M^{n'})^*             @.~\ra~@.
   \H^{\ell}(Y,Z^{t+1}M^{n'})^*           @.~\os{(\flat)}\ra~@.
   \H^{\ell-1}(Y,\gr_Z^{t}M^{n'})^*,
\end{CD}
}
$$
where we put $\ell:=N-i$, $t:=3e'-2-r$ and
    $E^*:=\Hom(E,\mu)$ for a $\Z/p\Z$-module $E$.
Hence Theorem \ref{prop:sdual2} is reduced to the following lemma
      by induction on $r \geq 1$ and the five lemma.
\addtocounter{thm}{5}
\begin{lem}\label{sublem6}
$c^{i,r}$ in \eqref{dual:pair:4} is non-degenerate
    for any $i$ and $r$ with $1 \leq r \leq 3e'-3$.
\end{lem}
\noindent
We prove this lemma in the next subsection.
\subsection{Proof of Lemma \ref{sublem6}}
We first give a brief review of linear Cartier operators.
Let $(s,{\cal L}_s)$ and ${\cal L}_{Y}$
    be as in \S\ref{sect8.2}, and
    let $(Y',{\cal L}_{Y'})$ be the log scheme
      defined by a cartesian diagram
\begin{equation}\label{def:Yp}
\begin{CD}
   (Y',{\cal L}_{Y'})  @>{\pr_2}>>  (Y,{\cal L}_Y) \\
   @V{\pr_1}V{\quad\quad\quad\;\square}V           @VVV \\
     (s, {\cal L}_{s}) @>{\F^{\abs}_{(s, {\cal L}_{s})}}>> (s, {\cal L}_{s}),
\end{CD}
\end{equation}
where $\F^{\abs}_{(s, {\cal L}_{s})}$ denotes the absolute Frobenius
       on $(s, {\cal L}_{s})$.
Let $\ul{\pr}{}_2:Y' \ra Y$ be the underlying morphism of
     schemes of $\pr_2$,
      and let
     $\F_{Y/s} : Y \ra Y'$ be the unique morphism of schemes
      such that $\ul{\pr}{}_2 \circ \F_{Y/s}$ agrees with the
   absolute Frobenius on $Y$.
Note that $\F_{Y/s}$ is a finite morphism of schemes.
We put
$\omega_{Y'}^N := \omega_{(Y',{\cal L}_{Y'})/(s,{\cal L}_s)}^N$
      for simplicity, where
        we regarded $(Y',{\cal L}_{Y'})$ as a smooth log scheme
           over $(s,{\cal L}_s)$ by $\pr_1$ in \eqref{def:Yp}.
By \cite{kf:log}, 5.3 and
    the same argument as for \cite{katz}, 7.2,
    there is a $\O_{Y'}$-linear isomorphism
\begin{equation}\label{def:cartier}
\begin{CD}
   C_{\lin}^{-1}: \omega_{Y'}^N @>{\simeq}>> \F_{Y/s*}(\omega_{Y}^N/\BB_Y^N).
\end{CD}
\end{equation}
(The following composite map gives the inverse Cartier operator
       $C^{-1}$ defined in \cite{hyodo}:
$$
\begin{CD}
\omega_{Y}^N  @>{{\mathrm{canonical}}}>>
   \ul{\pr}{}_{2*}\omega_{Y'}^N @>{\ul{\pr}{}_{2*}(C_{\lin}^{-1})}>>
   \ul{\pr}{}_{2*}\F_{Y/s*}(\omega_{Y}^N/\BB_Y^N)=\omega_{Y}^N/\BB_Y^N.\;)
\end{CD}
$$
Now we start the proof of Lemma \ref{sublem6}.
Let $C_{\lin}$ be the inverse of $C_{\lin}^{-1}$.
By \cite{hyodo2}, 3.2 and the same argument as for \cite{milne:duality}, 1.7,
there are $\O_{Y'}$-bilinear perfect pairings of locally free
   $\O_{Y'}$-modules of finite rank
\begin{equation}\notag
\begin{CD}
     \F_{Y/s*} (\omega_Y^{n-1}/\ZZ_Y^{n-1})
         @.~\times~@. \F_{Y/s*}
          (\omega_Y^{n'-2}/\ZZ_Y^{n'-2}) @.~ \lra ~@. \omega_{Y'}^N,
            @.\quad (x,y) \mapsto @. ~C_{\lin}((dx) \wedge y), \\
     \F_{Y/s*} (\ZZ_Y^{n-1}/\BB_Y^{n-1})
         @.~\times~@. \F_{Y/s*} (\ZZ_Y^{n'-1}/\BB_Y^{n'-1})
          @.~ \lra ~@. \omega_{Y'}^N,
            @.\quad (x,y) \mapsto @. C_{\lin}(x \wedge y), \\
     \F_{Y/s*} (\omega_Y^{n-2}/\ZZ_Y^{n-2})
         @.~\times~@. \F_{Y/s*} (\omega_Y^{n'-1}/\ZZ_Y^{n'-1})
          @.~ \lra ~@. \omega_{Y'}^N,
            @.\quad (x,y) \mapsto @. ~C_{\lin}(x \wedge dy).
\end{CD}
\end{equation}
By \cite{hyodo2}, Theorem 3.1
      and the Serre-Hartshorne duality \cite{H:RD},
      $\omega_{Y'}^N$ is a dualizing sheaf for $Y'$
           in the sense of \cite{H:AG}, p.\ 241, Definition.
Hence by \eqref{dual:isom:bkh} and
          Lemma \ref{sublem3.5},
          the pairing
$$
\begin{CD}
\H^i(Y,\gr_Z^{r}M^n) \times
               \H^{N-i}(Y,\gr_Z^{3e'-2-r}M^{n'})
            @>{f^{q,n}}>> \mu \otimes \H^N(Y,\omega_Y^N/\BB_Y^N)
            @>{\id \otimes \tr'_{Y/s}}>> \mu \otimes k
\end{CD}
$$
($q$ is the maximal integer with $3(q-1)<r$) is a non-degenerate pairing
      of finite-dimensional $k$-vector spaces.
Here $\tr'_{Y/s}$ denotes the $k$-linear trace map
$$
\begin{CD}
\H^N(Y,\omega_Y^N/\BB_Y^N)
     = \H^N(Y', \F_{Y/s*}(\omega_Y^N/\BB_Y^N))
      @>{C_{\lin}}>>
     \H^N(Y',\omega_{Y'}^N) @>>> k.
\end{CD}
$$
Finally, $c^{i,r}$ is non-degenerate
      by commutative squares
$$
{\small
\begin{CD}
    \H^N(Y,\gr_Z^{r}M^n \otimes \gr_Z^{3e'-2-r}M^{n'})
        @>{f^{q,n}}>>
     \mu \otimes \H^N(Y,\omega_Y^N/\BB_Y^N)
       @>{\id\otimes\tr'_{Y/s}}>> \mu \otimes k \\
     @V{\eqref{dual:factor2}}VV @VV{\chi}V
        @VV{\id \otimes \tr_{k/{\Bbb F}_p}}V \\
     \mu  @<{\id \otimes \tr_Y}<<
        \mu \otimes \H^{N+1}(Y,\nu_Y^N)
        @>{\id \otimes \tr_Y}>>
          \mu \otimes {\Bbb F}_p,
\end{CD}
}
$$
where the left square commutes by Theorem \ref{lem:trace}
   and the right square commutes
     by a similar argument as for \cite{sato:ss}, 3.4.1.
This completes the proof of Lemma \ref{sublem6} and Theorem \ref{prop:sdual2}.

\medskip

\section{\bf Duality of $p$-adic \'etale Tate twists}\label{sect10}
\medskip
In this section we prove Theorem \ref{thm0-4} using
    Theorem \ref{prop:sdual2}.
\subsection{Statement of the results}\label{sect10.1}
The setting is the same as in \S\ref{sect4.1}.
In this section,
   we assume that $X$ is {\it proper} over $B=\Spec(A)$, and
   that $A$ is either an algebraic integer ring (global case)
       or a henselian discrete valuation
         ring whose fraction field has characteristic $0$
          and whose residue field is finite of characteristic $p$
             (local case).
Let $d$ be the absolute dimension of $X$.
Throughout this section,
      $n$ and $r$ denote integers with $0 \leq n \leq d$ and $1 \leq r$.
The aim of this section is to prove the following duality results:
\begin{thm}\label{thm:main}
Assume that $A$ is local.
Then$:$
\begin{enumerate}
\item[(1)]
There is a canonical trace map
     $\tr_{(X,Y)} : \H^{2d+1}_Y(X,\T_r(d)_X) \ra \Z/p^r\Z$,
       which is bijective if $X$ is connected.
\item[(2)]
For $i \in \Z$, the natural pairing arising from \eqref{def:prod2} and
$\tr_{(X,Y)}$
\stepcounter{equation}
\begin{equation}\label{dual:main}
\hspace{-50pt}
\begin{CD}
\H^i_Y(X,\T_r(n)_X) \times \H^{2d+1-i}(X,\T_r(d-n)_X)
       @>>>  \Z/p^r\Z
\end{CD}
\hspace{-50pt}
\end{equation}
is a non-degenerate pairing of finite $\Z/p^r\Z$-modules.
\end{enumerate}
\end{thm}
\stepcounter{thm}
\begin{thm}[{{\bf \ref{thm0-4}}}]\label{thm:main2}
Assume that $A$ is global.
Then$:$
\begin{enumerate}
\item[(1)]
There is a canonical trace map
     $\tr_X : \H^{2d+1}_{c}(X,\T_r(d)_X) \ra \Z/p^r\Z$,
       where the subscript $c$ means the \'etale cohomology
         with compact support $($see $\S\ref{sect10.2}$ below$)$.
If $X$ is connected,
     then $\tr_X$ is bijective.
\item[(2)]
For $i \in \Z$, the natural pairing arising from \eqref{def:prod2} and $\tr_X$
\stepcounter{equation}
\begin{equation}\label{dual:main2}
\hspace{-50pt}
\begin{CD}
\H^i_c(X,\T_r(n)_X) \times \H^{2d+1-i}(X,\T_r(d-n)_X)
       @>>>  \Z/p^r\Z
\end{CD}
\hspace{-50pt}
\end{equation}
is a non-degenerate pairing of finite $\Z/p^r\Z$-modules.
\end{enumerate}
\end{thm}
\noindent
In \S\ref{sect10.2}, we will define the localized trace map $\tr_{(X,Y)}$
      and the global trace map $\tr_X$.
After showing a compatibility of these trace maps,
     we will reduce Theorem \ref{thm:main2} (2) to Theorem \ref{thm:main} (2).
We will prove Theorem \ref{thm:main} (2) in \S\S\ref{sect10.3}--\ref{sect10.5}.
\stepcounter{thm}
\begin{rem}\label{rem:PTdual}
If $A$ is local,
     there is a natural pairing of finite $\Z/p^r\Z$-modules
\stepcounter{equation}
\begin{equation}\label{PTdual}
\hspace{-50pt}
\begin{CD}
\H^i(V,\mu_{p^r}^{\otimes n}) \times
     \H^{2d-i}(V,\mu_{p^r}^{\otimes d-n})
       @>>> \H^{2d}(V,\mu_{p^r}^{\otimes d})
          \simeq \Z/p^r\Z,
\end{CD}
\hspace{-50pt}
\end{equation}
where $V$ denotes $X_K$ with $K:= \Frac(A)$.
As is well-known,
   this pairing is non-degenerate by
   the Tate duality for $K$
   and the Poincar\'e duality for $V_{\ol K}$.
Theorem $\ref{thm:main}$ $(2)$ does not follow from these facts,
    although Theorem $\ref{thm:main}$ implies the non-degeneracy of
\eqref{PTdual}.
We will deduce Theorem $\ref{thm:main}$ $(2)$
      from Theorems $\ref{thm:milne}$ and $\ref{prop:sdual2}$.
\end{rem}
\subsection{Trace maps}\label{sect10.2}
We first construct the localized trace map $\tr_{(X,Y)}$,
    assuming that $A$ is local.
Let $\iota:Y \hra X$ be the natural closed immersion.
By Lemma \ref{lem4'} and Theorem \ref{thm:milne},
   $\H^{i}_Y(X,\T_r(d)_X)$
     is zero for any $i \geq 2d+2$.
We define $\tr_{(X,Y)}$ as the composite
$$
\begin{CD}
\tr_{(X,Y)} : \H^{2d+1}_Y(X,\T_r(d)_X)
     @>{(\gys_{\iota}^{d})^{-1}}>{\simeq}>
      \H^{d}(Y,\nu_{Y,r}^{d-1}) @>{\tr_{Y}}>> \Z/p^r\Z,
\end{CD}
$$
which is bijective if $X$ is connected
      (i.e., $Y$ is connected).
We next define the global trace map $\tr_X$, assuming that $A$ is global.
For a scheme $Z$ which is separated of finite type over $B$
   and an object $\K \in D^+(Z_{\et},\Z/p^r\Z)$,
    we define $\H^*_c(Z,\K)$ as $\H^*_c(B,Rf_!\K)$,
  where $f$ denotes the structural morphism $Z \to B$ and
   $\H^*_c(B,\bullet)$ denote the \'etale cohomology groups with compact 
support of $B$
        (cf.\ \cite{milne:adual}, II.2).
By the Kummer sequence \eqref{DT:weight1}
      and the isomorphism $\H^3_c(B,\Gm) \simeq \qz$
       (cf.\ \cite{milne:adual}, II.2.6),
         there is a trace map $\H^3_c(B,\T_r(1)_B) \to \Z/p^r\Z$
(cf.\ \cite{jss}, Corollary 4.3 (a)).
We normalize this map so that for a closed point $i_s : s \hra B$,
   the composite map
$$
\begin{CD}
\H^1(s,\Z/p^r\Z) @>{\gys_{i_s}^1}>> \H^3_c(B,\T_r(1)_B) @>>> \Z/p^r\Z
\end{CD}
$$
coincides with the trace map of $s$ (defined in \ref{thm:milne} (1)).
We define the trace map $\tr_X$ as the composite
$$
\begin{CD}
\tr_X: \H^{2d+1}_c(X,\T_r(d)_X) @>>>
          \H^3_c(B,\T_r(1)_B) @>>> \Z/p^r\Z,
\end{CD}
$$
where the first arrow arises from the trace morphism in Theorem 
\ref{prop:trace}.
The bijectivity assertion for $\tr_X$ in
  Theorem \ref{thm:main2} (1) will follow from \ref{thm:main2} (2).
We show here the following:
\begin{lem}\label{lem:com}
Assume that $A$ is global.
Then there is a commutative diagram
$$
\begin{CD}
\H^{2d+1}_Y(X,\T_r(d)_X) @>{\tr_{(X,Y)}}>> \Z/p^r\Z\\
     @V{\iota_*}VV  @|\\
\H^{2d+1}_{c}(X,\T_r(d)_X) @>{\tr_{X}}>> \Z/p^r\Z,
\end{CD}
$$
where the arrow $\iota_*$ denotes the canonical adjunction map and
$\tr_{(X,Y)}$ denotes the sum of the localized trace maps
for the connected components of $Y$.
\end{lem}
\begin{pf}
Let $\{Y_i\}_{i \in I}$ be the connected components of $Y$.
Let $x$ be a closed point on $Y$
        with $i_x : x \hra X$.
Noting that $\H^{2d+1}_Y(X,\T_r(d)_X) \simeq
    \bigoplus_{i \in I}~\Z/p^r\Z$,
      consider a diagram
$$
\begin{CD}
\H^1(x,\Z/p^r\Z) @>{\gys_{i_x}^d}>>
     \H^{2d+1}_Y(X,\T_r(d)_X) @>{\tr_{(X,Y)}}>> \Z/p^r\Z\\
@|  @V{\iota_*}VV  @|\\
\H^1(x,\Z/p^r\Z) @>{\gys_{i_x}^d}>>
\H^{2d+1}_{c}(X,\T_r(d)_X) @>{\tr_{X}}>> \Z/p^r\Z.
\end{CD}
$$
Since the left square commutes,
    it suffices to show that the composite of the upper row is bijective
        and that the outer rectangle is commutative.
The composite of the upper row agrees with
    the trace map for $x$ by Theorem \ref{thm:milne} (1).
In particular it is bijective.
The composite of the lower row
    agrees with the trace map for $x$ by Theorem \ref{thm:relative} (2).
We are done.
\end{pf}
We reduce Theorem \ref{thm:main2} (2) to Theorem \ref{thm:main} (2).
Assume that $A$ is global.
We use the notation in \S\ref{sect4.1}.
Put $X_{\Sigma} := \coprod_{s \in \Sigma}~X \times_B B_s$.
Since $j^*\T_r(n)_X \simeq \mu_{p^r}^{\otimes n}$,
     there is a distinguished triangle
$$
\begin{CD}
\T_r(n)_X @>{\iota^*}>>
    R\iota_*\iota^*\T_r(n)_X @>>>
    j_!\mu_{p^r}^{\otimes n} [1]
           @>{j_!}>> \T_r(n)_X [1],
\end{CD}
$$
where the arrow $\iota^*$ (resp.\ $j_!$)
    denotes the canonical adjunction morphism $\id \ra R\iota_*\iota^*$
     (resp.\ $Rj_!j^* \ra \id$).
By Lemma \ref{lem:com} and
     the proper base-change theorem:
      $\H^*(Y,\iota^*\T_r(n)_X) \simeq
       \H^*(X_{\Sigma},\T_r(n)_{X_{\Sigma}})$,
     we obtain a commutative diagram with exact rows
    (after changing the signs of $(\flat)$ suitably)
$$
{\small
\begin{CD}
\H^{i-1}(X_{\Sigma},\T_r(n)_{X_{\Sigma}}) @.~ \ra ~@.
     \H^i_{c}(V,\mu_{p^r}^{\otimes n}) @.~ \ra ~@.
     \H^i_{c}(X,\T_r(n)_X) @.~ \ra ~@.
     \H^i(X_{\Sigma},\T_r(n)_{X_{\Sigma}}) @.~ \ra ~@.
      \H^{i+1}_{c}(V,\mu_{p^r}^{\otimes n})\\
    @V{\eqref{dual:main}}VV @. @V{a}VV @. @V{\eqref{dual:main2}}VV @.
    @V{\eqref{dual:main}}VV @. @V{a}VV\\
\H^{\ell+1}_Y(X,\T_r(m)_{X})^* @.~ \os{(\flat)}{\ra} ~@.
     \H^{\ell}(V,\mu_{p^r}^{\otimes m})^* @.~ \ra ~@.
     \H^{\ell}(X,\T_r(m)_X)^* @.~ \ra ~@.
     \H^{\ell}_Y(X,\T_r(m)_{X})^* @.~ \os{(\flat)}{\ra} ~@.
     \H^{\ell-1}(V,\mu_{p^r}^{\otimes m})^*. \\
\end{CD}
}
$$
Here the superscript $*$ means the Pontryagin dual
    and we put $\ell:= 2d+1-i$ and $m:= d-n$.
The lower row is the dual of the localization long exact sequence
      and the vertical arrows arise from duality pairings.
The arrows $a$ are isomorphisms of finite groups
     by the Artin-Verdier duality and the absolute purity
      \cite{Th}, \cite{fujiwara}.
Thus Theorem \ref{thm:main2} (2) is reduced to Theorem \ref{thm:main} (2)
     by the five lemma.
\subsection{Reduction to the case $r=1$}\label{sect10.3}
We start the proof of Theorem \ref{thm:main} (2),
     which will be completed in \S\ref{sect10.5}.
By the distinguished triangle in Proposition \ref{prop:bock} (3),
         the problem is reduced to the case $r=1$.
Furthermore we may assume that
     $K=\Frac(A)$ contains a primitive $p$-th root of unity $\zeta_p$.
Indeed, otherwise the scalar extension
         $X_{A'}:=X \otimes_{A} A'$,
            where $A'$ denotes the integer ring of $K(\zeta_p)$,
     again satisfies the condition \ref{cond2} over $\Spec(A)$.
Hence once we show Theorem \ref{thm:main} (2) for $X_{A'}$,
    we will obtain Theorem \ref{thm:main} (2) for $X$
       by a standard norm argument and Corollary \ref{cor:proj}.
\subsection{Descending induction on $n$}\label{sect10.4}
Assume that $\zeta_p \in K$ and $r=1$.
We prove this case of Theorem \ref{thm:main} (2)
     by descending induction on $n \leq d$.
Let $N$ be the relative dimension $\dim(X/B)$.
If $n=N+1(=d)$,
       \eqref{dual:main} is isomorphic to the pairing
$$
\begin{CD}
\H^{i-N-2}(Y,\nu_{Y,1}^N) \times
     \H^{2N+3-i}(Y,\Z/p\Z) \lra
      \H^{N+1}(Y,\nu_{Y,1}^N) \simeq \Z/p\Z
\end{CD}
$$
by the proper base-change theorem and Lemma \ref{lem4'}.
This pairing is a non-degenerate pairing of finite $\Z/p\Z$-modules
   by Theorem \ref{thm:milne}.
To proceed the descending induction on $n$,
       we study the inductive structure
          of $\{\T_1(n)_X \}_{n \geq 0}$ on $n$.
We fix some notation.
Let $\iota : Y \hra X$ and $j : V(=X_K) \hra X$ be as before.
Let $\nu_Y^n$, $\mu'$ and $\mu$
     be as in \eqref{def:simple}.
See also the remark after \eqref{def:simple}.
Put $\lam_Y^n:=\lam_{Y,1}^n$ and $\T(n)_X := \T_1(n)_X$.
Now for $n$
     with $1 \leq n \leq N+1$, we define the morphism
$$
\begin{CD}
\ind_n:(j_*\mu_p) \otimes \T(n-1)_X
      \left( :=(j_*\mu_p) \otimes^{\L} \T(n-1)_X \right)
      @>>> \T(n)_X
\end{CD}
$$
by restricting the product structure
       $\T(1)_X \otimes^{\L} \T(n-1)_X
      \ra \T(n)_X$ to the $0$-th cohomology sheaf $j_*\mu_p$
         of $\T(1)_X$.
\begin{lem}\label{sublem1}
Let
\stepcounter{equation}
\begin{equation}\label{dual:DT1}
\hspace{-30pt}
\begin{CD}
\KK(n)[-1] @>{b_n}>> \mu' \otimes \iota^*\T(n-1)_X
      @>{\iota^*(\ind_n)}>> \iota^*\T(n)_X
      @>{a_n}>> \KK(n)
\end{CD}
\hspace{-30pt}
\end{equation}
    be a distinguished triangle in $D^b(Y_{\et},\Z/p\Z)$.
Then$:$
\begin{enumerate}
\item[(1)]
The triple $(\KK(n),a_n,b_n)$
          is unique up to a unique isomorphism
       in $D^b(Y_{\et},\Z/p\Z)$, and
           $b_n$ is determined by
             the pair $(\KK(n),a_n)$.
\item[(2)] $\KK(n)$ is concentrated in
      $[n-1,n]$ and $a_n$ induces isomorphisms
$$
{\cal H}^q(\KK(n)) \simeq
     \begin{cases}
       \mu' \otimes \nu_{Y}^{n-2}
           \quad &(q = n-1), \\
        FM^n      \quad &(q=n),
     \end{cases}
$$
where
     $M^n$ denotes the \'etale sheaf
      $\iota^*R^nj_*\mu_p^{\otimes n}$ on $Y$, and
      $FM^n$ denotes the \'etale subsheaf of $M^n$
       defined in $\S\ref{sect3.4}$.
\item[(3)]
     There is a distinguished triangle
     in $D^b(Y_{\et},\Z/p\Z)$
\begin{equation}\label{dual:DT2}
\hspace{-30pt}
\begin{CD}
\KK(n)[-1] @>{c_n}>> \mu' \otimes R \iota^!\T(n-1)_X
       @>{R \iota^!(\ind_n)}>> R \iota^! \T(n)_X @>{d_n}>> \KK(n),
\end{CD}
\hspace{-30pt}
\end{equation}
where $c_n$ and $d_n$ are morphisms
     determined by the pair $(\KK(n),a_n)$.
\item[(4)]
There is an anti-commutative diagram
\begin{equation}\label{ancn:anti}
\begin{CD}
     \iota^*\T(n)_X  @>{\mathrm{canonical}}>>
      \mu' \otimes \iota^*Rj_*\mu_p^{\otimes n-1}\\
     @V{a_n}VV
       @VV{\id_{\mu'}\otimes \iota^*(\delta_{V,Y}^{\loc}(\T(n-1)_X))}V\\
     \KK(n)     @>{c_n[1]}>>  \mu'\otimes R\iota^!\T(n-1)_X[1].
\end{CD}
\end{equation}
\end{enumerate}
\end{lem}
\begin{pf}
(2) follows from
     the long exact sequence of cohomology sheaves
      associated with \eqref{dual:DT1}
      and the isomorphism of sheaves
     $\mu' \otimes \iota^*R^qj_*\mu_p^{\otimes n-1}
        \simeq \iota^*R^qj_*\mu_p^{\otimes n}$
          (cf.\ \eqref{defdef2}).
The details are straight-forward and left to the reader.
By (2) and Lemma \ref{lem:CD},
     there is no non-zero morphism from
            $\mu' \otimes \iota^*\T(n-1)_X$ to $\KK(n)[-1]$.
Lemma \ref{sublem1} (1) follows from this fact and
        Lemma \ref{lem:CD2} (3).
We next prove (3).
Let
$$
\begin{CD}
    \FF[-1] @>{b'}>> (j_*\mu_p) \otimes \T(n-1)_X
       @>{\ind_n}>> \T(n)_X @>{a'}>> \FF
\end{CD}
$$
be a distinguished triangle in $D^b(X_{\et},\Z/p\Z)$.
By a similar argument as for the claim (2),
      the cohomology sheaves of $\FF$
         are supported on $Y$.
This implies that $\FF=R\iota_*\iota^*\FF$.
By the uniqueness assertion of the claim (1),
     the triple $(\iota^*\FF,\iota^*(a'),\iota^*(b'))$
       is isomorphic to $(\KK(n),a_n,b_n)$
         by a unique isomorphism.
Under this identification, we have $a' \simeq R\iota_*(a_n)$.
Moreover, $b'$ is determined by the pair
        $(\FF,a')=(R\iota_*\KK(n),R\iota_*(a_n))$
          by a similar argument as for the claim (1).
Hence applying $R\iota^!$ to the above triangle,
      we obtain the distinguished triangle
        \eqref{dual:DT2} with
          $c_n=R\iota^!(b')$ and $d_n=R\iota^!R\iota_*(a_n)$.
Finally, (4) follows from an elementary computation on connecting morphisms,
whose details are straight-forward and left to the reader.
\end{pf}
\medbreak
\noindent
In what follows,
     we fix a pair $(\KK(n), a_n)$
       fitting into \eqref{dual:DT1} for each $n$ with $1 \leq n \leq N+1$.
By Lemma \ref{sublem1},
    the morphisms $b_n$, $c_n$ and $d_n$ fitting into
     \eqref{dual:DT1} and \eqref{dual:DT2} are determined by $(\KK(n), a_n)$.
Next we construct a pairing on
     $\{\KK(n) \}_{1 \leq n \leq N+1}$
        using $\{a_n \}_{1 \leq n \leq N+1}$.
Let us note that for objects $\K_1, \K_2 \in D^-(Y_{\et},\Z/p\Z)$,
       and $\K_3 \in D^+(Y_{\et},\Z/p\Z)$,
       we have
$$
\begin{CD}
\Hom_{D(Y_{\et},\Z/p\Z)}(\K_1 \otimes^{\L} \K_2,\K_3)
      \simeq
       \Hom_{D(Y_{\et},\Z/p\Z)}\left(\K_1,
          R\sHom_{Y,\Z/p\Z}(\K_2,\K_3) \right).
\end{CD}
$$
For $\K \in D^-(Y_{\et},\Z/p\Z)$,
     we define
\begin{equation}\label{def:D}
\begin{CD}
{\Bbb D}(\K):=R\sHom_{Y,\Z/p\Z}(\K,
        \mu' \otimes \nu_{Y}^N[-N-2])
             \in  D^+(Y_{\et},\Z/p\Z).
\end{CD}
\end{equation}
\addtocounter{thm}{4}
\begin{lem}\label{sublem2}
Let $n$ be as before and
      put $n':=N+2-n$.
Then there is a unique morphism
\stepcounter{equation}
\begin{equation}\label{dual:pair1}
\hspace{-50pt}
\begin{CD}
    (\KK(n)[-1]) ~{\otimes}^{\L}~ \KK(n') @>>>
                 \mu' \otimes \nu_{Y}^N[-N-2]
                  \quad \hbox{ in } \; D^-(Y_{\et},\Z/p\Z)
\end{CD}
\hspace{-50pt}
\end{equation}
      whose adjoint morphism
       $\KK(n)[-1] \ra \DD(\KK(n'))$ fits
         into a commutative diagram
            with distinguished rows
         $($cf.\ $\eqref{dual:DT1}$, $\eqref{dual:DT2})$
$$
{\small
\begin{CD}
    \KK(n)[-1]
      @>{c_{n}}>> \mu' \otimes R\iota^!\T(n-1)_X
      @>{R\iota^!(\ind_{n})}>> R\iota^!\T(n)_X
      @>{d_{n}}>> \KK(n) \\
@VVV @V{(\sharp)}VV @V{(\sharp)}VV @VVV \\
    \DD(\KK(n'))
      @>{\DD(a_{n'})}>> \DD(\iota^*\T(n')_X)
      @>{\DD(\ind_{n'})}>> \DD(\mu' \otimes \iota^*\T(n'-1)_X)
      @>{-\DD(b_{n'})}>> \DD(\KK(n')[-1]).
\end{CD}
}
$$
Here the vertical arrows $(\sharp)$ come from
           the product structure of $\{\T(n)_X\}_{n \geq 0}$,
            the identity map of $\mu'$ and
    the Gysin isomorphism $\gys_{\iota}^{N+1}$ in Lemma $\ref{lem4'}$
        $($the commutativity of the central square is
           easy and left to the reader$)$.
\end{lem}
\begin{pf}
The assertion follows from Lemma \ref{lem:CD2} (1) and the fact that
\begin{align*}
    & \Hom_{D^+(Y_{\et},\Z/p\Z)}
     \big(\KK(n),
    \DD(\mu' \otimes \iota^*\T(n'-1)_X) \big) \\
      \simeq \, & \Hom_{D^-(Y_{\et},\Z/p\Z)}
      \big( \KK(n) \otimes^{\L} (\mu' \otimes \iota^*\T(n'-1)_X),
         \mu' \otimes \nu_Y^N [-N-2] \big) = 0,
\end{align*}
where the last equality follows from
     Lemma \ref{sublem1} (2) and Lemma \ref{lem:CD}.
\end{pf}
Now we turn to the proof of Theorem \ref{thm:main} (2)
     and claim the following:
\stepcounter{thm}
\begin{prop}\label{prop:sdual}
Let $n$ and $n'$ be as in Lemma $\ref{sublem2}$.
Then for $i \in \Z$, the pairing
\stepcounter{equation}
\begin{equation}\label{dual:pair2}
\hspace{-50pt}
\begin{CD}
\H^i(Y,\KK(n))
      \times \H^{2N+2-i}(Y,\KK(n')) \lra
                  \mu \otimes
                   \H^{N+1}(Y, \nu_{Y}^N)
                    @>{\id \otimes \tr_Y}>> \mu,
\end{CD}
\hspace{-50pt}
\end{equation}
induced by $\eqref{dual:pair1}$,
is a non-degenerate pairing of finite $\Z/p\Z$-modules.
\end{prop}
\noindent
We will prove this proposition
       in the next subsection.
We first finish the proof of Theorem \ref{thm:main} (2)
      by descending induction on $n \leq N+1$,
   admitting Proposition \ref{prop:sdual}.
See the beginning of this subsection for the case $n=N+1$.
Indeed, we obtain Theorem \ref{thm:main} (2)
      from Proposition \ref{prop:sdual},
           applying the following general lemma to
        the commutative diagram in Lemma \ref{sublem2}:
\stepcounter{thm}
\begin{lem}\label{subcor1}
Let ${\cal K}_1 \ra {\cal K}_2 \ra {\cal K}_3 \ra {\cal K}_1[1]$
     and ${\cal L}_3 \ra {\cal L}_2 \ra  {\cal L}_1 \ra {\cal L}_3[1]$
     be distinguished triangles in $D^b(Y_{\et},\Z/p\Z)$,
      and suppose that we are given a commutative diagram
$$
\begin{CD}
      {\cal K}_1 @>>> {\cal K}_2 @>>> {\cal K}_3 @>>> {\cal K}_1[1]\\
        @V{\alpha_1}VV @V{\alpha_2}VV @V{\alpha_3}VV @V{\alpha_1[1]}VV\\
      \DD({\cal L}_1) @>>> \DD({\cal L}_2) @>>>
        \DD({\cal L}_3) @>>> \DD({\cal L}_1) [1]\\
\end{CD}
$$
      $($with distinguished rows$)$ in $D^+(Y_{\et},\Z/p\Z)$.
For $m \in \{ 1,2,3 \}$ and $i \in \Z$, let
\begin{equation}\notag
\hspace{-50pt}
\begin{CD}
\beta_{m}^i : \H^i(Y,{\cal K}_{m})
      \times \H^{2N+3-i}(Y,{\cal L}_{m}) \lra
                  \mu \otimes
                   \H^{N+1}(Y, \nu_{Y}^N)
                    @>{\id \otimes \tr_Y}>> \mu
\end{CD}
\hspace{-50pt}
\end{equation}
     be the pairing induced by
       the adjoint morphism
        ${\cal K}_{m} \otimes^{\L} {\cal L}_{m} \ra
          \mu' \otimes\nu_{Y}^N [-N-2]$
            of $\alpha_{m}$.
Put $\ell:=2N+3-i$.
Then there is a commutative diagram
       with exact rows
$$
{
\begin{CD}
      \H^{i-1}(Y,{\cal K}_3) @.~\lra~@.
      \H^{i}(Y,{\cal K}_1) @.~\lra~@.
      \H^{i}(Y,{\cal K}_2) @.~\lra~@.
      \H^{i}(Y,{\cal K}_3) @.~\lra~@.
      \H^{i+1}(Y,{\cal K}_1) \\
      @V{\gamma_3^{i-1}}VV @. @V{\gamma_1^{i}}VV @.
      @V{\gamma_2^{i}}VV @. @V{\gamma_3^{i}}VV @. @V{\gamma_1^{i+1}}VV\\
      \H^{\ell+1}(Y,{\cal L}_3)^* @.~\os{(\flat)}{\lra}~@.
      \H^{\ell}(Y,{\cal L}_1)^* @.~\lra~@.
      \H^{\ell}(Y,{\cal L}_2)^* @.~\lra~@.
      \H^{\ell}(Y,{\cal L}_3)^* @.~\os{(\flat)}{\lra}~@.
      \H^{\ell-1}(Y,{\cal L}_1)^*
\end{CD}
}
$$
       after changing the signs of $(\flat)$ suitably.
Here for a $\Z/p\Z$-module $E$, $E^*$ denotes $\Hom_{\Z/p\Z}(E,\mu)$,
   and $\gamma_{m}^i$ denotes the natural map induced by $\beta_{m}^i$.
Furthermore, if $\gamma_1^i$ and $\gamma_3^i$ are
      bijective for any $i$, then $\gamma_2^i$ is
        bijective for any $i$.
\end{lem}
\begin{pf}
For each $m$ and $i$,
     $\gamma_{m}^i$ factors as follows:
$$
\begin{CD}
\H^i(Y,{\cal K}_{m}) @>{\H^i(Y,\alpha_{m})}>>
     \Ext^{i-N-2}_{Y,\Z/p\Z}({\cal L}_{m}, \mu' \otimes \nu_{Y}^N)
        @>>> \Hom_{\Z/p\Z}(\H^{2N+3-i}(Y,{\cal L}_{m}),\mu),
\end{CD}
$$
where the last map arises from a Yoneda pairing
     and the trace isomorphism
       $\H^{N+1}(Y,\mu' \otimes \nu_Y^N) \simeq \mu$.
The commutativity of the diagram of cohomology groups
     in the lemma follows from the functoriality of this decomposition.
The last assertion follows from the five lemma.
\end{pf}
\subsection{Proof of Proposition \ref{prop:sdual}}\label{sect10.5}
Let us recall that the canonical pairings
\begin{equation}\label{pair:sdual}
\begin{CD}
\H^i(Y,\nu_Y^q) @. \times @.
      \H^{N+1-i}(Y,\lam_Y^{N-q}) @>{\eqref{logHW:pair}}>> \Z/p\Z,\\
\H^i(Y,U^1M^n) @.~ \times ~@.
      \H^{N-i}(Y,U^1M^{n'}) @>{\eqref{dual:pair3}}>> \mu
\end{CD}
\end{equation}
($q=n'-2$ or $n-2$) are non-degenerate pairings of finite groups
     for any $i$ by
       Theorems \ref{thm:milne} and \ref{prop:sdual2}, respectively.
We deduce Proposition \ref{prop:sdual} from these results.
Let ${\Bbb U}(n)$ be an object of $D^b(Y_{\et},\Z/p\Z)$
      fitting into a distinguished triangle
$$
\begin{CD}
\lam_Y^{n}[-n-1]
       @>>> {\Bbb U}(n) @>>> \KK(n) @>>> \lam_Y^{n}[-n],
\end{CD}
$$
where the last morphism is defined as
      the composite
$\KK(n) \ra {\cal H}^{n}(\KK(n))[-n]
       \simeq FM^{n}[-n]
          \ra \lam_Y^{n}[-n]$
(cf.\ Lemma \ref{sublem1} (2),
      Theorem \ref{thm:vcyc}, Corollary \ref{cor:vcyc}).
By Lemma \ref{sublem1} (2) and Lemma \ref{lem:CD2} (3),
${\Bbb U}(n)$ is concentrated in $[n-1,n]$
     and unique up to a unique isomorphism.
We have
\begin{equation}\label{pair:sdual2}
{\cal H}^q({\Bbb U}(n)) \simeq
     \begin{cases}
       \mu' \otimes \nu_{Y}^{n-2}
           \qquad &(q = n-1), \\
       U^1M^n      \qquad &(q=n).
     \end{cases}
\end{equation}
For $\K \in D^b(Y_{\et},\Z/p\Z)$,
     let
     ${\Bbb D}(\K)$ be as in \eqref{def:D}.
In view of Lemma \ref{subcor1} and the non-degeneracy of the pairings
        in \eqref{pair:sdual},
        we have only to show the following:
\addtocounter{thm}{2}
\begin{lem}\label{sublem3}
\begin{enumerate}
\item[(1)]
There is a unique morphism
$$
\begin{CD}
f : {\Bbb U}(n)[-1] @>>> {\Bbb D}(FM^{n'}[-n'])
       \quad \hbox{ in } \; D^+(Y_{\et},\Z/p\Z)
\end{CD}
$$
fitting into a commutative diagram with distinguished rows
\stepcounter{equation}
\begin{equation}\label{dual:commute2}
\hspace{-80pt}
{\small
\begin{CD}
    {\Bbb U}(n)[-1]
            @.~ \lra ~@. \KK(n)[-1] @.~ \lra ~@. \lam_Y^{n}[-n-1]
            @.~ \lra ~@. {\Bbb U}(n) \\
    @V{f}VV @. @V{\eqref{dual:pair1}}VV @.
     @V{(-1)^n\cdot f_1}VV @. @V{f[1]}VV \\
    {\Bbb D}(FM^{n'}[-n'])
      @.~ \lra ~@. {\Bbb D}(\KK(n'))
      @.~ \lra ~@. {\Bbb D}(\mu' \otimes\nu_{Y}^{n'-2}[-n'+1])
      @.~ \lra ~@. {\Bbb D}(FM^{n'}[-n'])[1].
\end{CD}
}
\hspace{-80pt}
\end{equation}
Here the lower row arises from
    a distinguished triangle obtained by truncation
$$
\begin{CD}
     \mu' \otimes \nu_{Y}^{n'-2}[-n'+1]
     @>>> \KK(n') @>>> FM^{n'}[-n']
     @>>> \mu' \otimes \nu_{Y}^{n'-2}[-n'+2]
\end{CD}
$$
$($cf.\ Lemma $\ref{sublem1}$ $(2))$,
and we have chosen the signs of the connecting morphisms
       $(=$the last arrows$)$ of the both rows suitably.
The arrow $f_1$ is defined as the adjoint morphism of the map
     $\lam_Y^{n}[-n-1] \otimes^{\L} (\mu' \otimes\nu_{Y}^{n'-2})[-n'+1]
       \ra \mu' \otimes\nu_{Y}^{N}[-N-2]$
       induced by the identity map of $\mu'$
    and the pairing $\eqref{def:prod}$.
\item[(2)]
There is a commutative diagram
       with distinguished rows in $D^+(Y_{\et},\Z/p\Z)$
\begin{equation}\label{dual:commute3}
\hspace{-80pt}
{\small
\begin{CD}
    \mu' \otimes  \nu_Y^{n-2}[-n]
         @.~ \ra ~@. {\Bbb U}(n)[-1] @.~ \ra ~@. U^1M^{n}[-n-1]
         @.~ \ra ~@. \mu' \otimes  \nu_Y^{n-2}[-n+1]\\
    @V{f_2}VV @. @V{f}VV @. @V{f_3}VV @. @V{f_2[1]}VV\\
    \DD(\lam_Y^{n'}[-n'])
     @.~ \ra ~@. \DD(FM^{n'}[-n']) @.~ \ra ~@. \DD(U^1M^{n'}[-n'])
     @.~ \ra ~@. \DD(\lam_Y^{n'}[-n'])[1].\\
\end{CD}
}
\hspace{-80pt}
\end{equation}
Here the upper row is the distinguished triangle obtained by truncation
      $($cf.\ $\eqref{pair:sdual2})$,
the lower row arises from the short exact sequence
$0 \ra U^1M^n \ra FM^n \ra \lam_Y^n \ra 0$,
     and we have chosen the signs of the connecting morphisms
       $(=$the last arrows$)$ of the both rows suitably.
The arrow $f_2$ is defined in a similar way as for $f_1$ in $(1)$,
      and $f_3$ denotes the morphism induced by $\Theta^n[-1]$.
See $\S\ref{sect8.1}$ for $\Theta^n$.
\end{enumerate}
\end{lem}
\noindent
To prove Lemma \ref{sublem3},
     we first show Lemma \ref{sublem3'} below.
Note that for $\K \in D^b(X_{\et},\Z/p\Z)$,
            $Rj_*j^*\K$ and $R\iota^!\K$ are both bounded
              (cf.\ Lemma \ref{lem4''}).
For $\K, {\cal L} \in D^b(X_{\et},\Z/p\Z)$,
                $\K \otimes^{\L} {\cal L}$ is bounded,
          because a $\Z/p\Z$-sheaf is flat over $\Z/p\Z$.
\addtocounter{thm}{2}
\begin{lem}\label{sublem3'}
For $\K,{\cal L} \in D^b(X_{\et},\Z/p\Z)$,
    there is a commutative diagram
\stepcounter{equation}
\begin{equation}\label{dual:commute4}
\begin{CD}
     (Rj_*j^*\K) \otimes ^{\L} {\cal L} @>>> Rj_*j^*(\K \otimes ^{\L} {\cal 
L})\\
      @V{\delta^{\loc}_{V,Y}(\K)\otimes^{\L}\id}VV
         @VV{\delta^{\loc}_{V,Y}(\K \otimes ^{\L} {\cal L})}V\\
    (\iota_*R\iota^!\K[1]) \otimes ^{\L} {\cal L} @>>>
          \iota_*R\iota^!(\K \otimes ^{\L} {\cal L})[1],
\end{CD}
\end{equation}
where the horizontal arrows are natural product morphisms.
\end{lem}
\noindent
The commutativity of the induced diagram of cohomology sheaves of
       \eqref{dual:commute4} would be well-known.
However, we include a proof of the lemma,
     because we need the commutativity in the derived category
   to show especially Lemma \ref{sublem3} (2).
\begin{pf*}{Proof of Lemma \ref{sublem3'}}
For two complexes $M^{\bullet}$ and $N^{\bullet}$,
     let $(M^{\bullet}\otimes N^{\bullet})^{\tot}$
      be as in Proof of Proposition \ref{prop:proj}.
For a map $h^{\bullet}:M^{\bullet} \ra N^{\bullet}$ of complexes,
     let $\cone(h)^{\bullet}$ be as in Proof of Proposition \ref{prop:bock}
        and let $u_h:N^{\bullet} \ra \cone(h)^{\bullet}$ be
          the canonical map.
Let $C^b(X_{\et},\Z/p\Z)$ be the category of
      bounded complexes of $\Z/p\Z$-sheaves on $X_{\et}$.
Take an $\iota^!$-acyclic resolution $K^{\bullet} \in C^b(X_{\et},\Z/p\Z)$
      of $\K$ (see the remark before Lemma \ref{sublem3'})
        and a bounded complex
           $L^{\bullet}\in C^b(X_{\et},\Z/p\Z)$ which represents
          ${\cal L}$.
Note that $K^{\bullet}$
            is a $j_*$-acyclic resolution of $\K$ as well.
We further take an injective resolution $J^{\bullet} \in C^+(X_{\et},\Z/p\Z)$
         of $(K^{\bullet} \otimes L^{\bullet})^{\tot}$.
Let $f:K^{\bullet} \ra j_*j^*K^{\bullet}$ and
      $g:J^{\bullet} \ra j_*j^*J^{\bullet}$ be the canonical maps,
     and let $f':(K^{\bullet} \otimes L^{\bullet})^{\tot}
             \ra ((j_*j^*K^{\bullet}) \otimes L^{\bullet})^{\tot}$
     be the map induced by $f$.
Then in $D^b(X_{\et},\Z/p\Z)$,
    the diagram \eqref{dual:commute4} decomposes as follows:
$$
{\scriptsize
\hspace{-15pt}
\begin{CD}
     ((j_*j^*K^{\bullet}) \otimes L^{\bullet})^{\tot} @=
        ((j_*j^*K^{\bullet}) \otimes L^{\bullet})^{\tot}
       @= ((j_*j^*K^{\bullet}) \otimes L^{\bullet})^{\tot}
          @>{\varphi_1}>> j_*j^*J^{\bullet}
            @= j_*j^*J^{\bullet} \\
      @V{\delta^{\loc}_{V,Y}(\K)\otimes^{\L}\id}V{\qquad \qquad \quad (1)}V
       @V{(u_f \otimes \id)^{\tot} }V{\qquad \qquad \quad (2)}V
           @V{u_{f'}}V{\qquad \quad \quad\;\;\;(3)}V
             @V{u_g}V{\qquad \quad (4)}V
         @VV{\delta^{\loc}_{V,Y}(\K \otimes ^{\L} {\cal L})}V\\
    ((\iota_*\iota^!K^{\bullet}[1]) \otimes L^{\bullet})^{\tot}
      @>{\varphi_2}>{\simeq}>
        (\cone(f)^{\bullet} \otimes L^{\bullet})^{\tot} @>{\varphi_3}>{=}>
           \cone(f')^{\bullet}
            @>{\varphi_4}>> \cone(g)^{\bullet} @<{\varphi_5}<{\simeq}<
          \iota_*\iota^!J^{\bullet}[1],
\end{CD}
}
$$
where $\varphi_1$, $\varphi_2$, $\varphi_4$ and $\varphi_5$
       are canonical maps of complexes
          and $\varphi_2$ and $\varphi_5$ are isomorphisms in
             $D^b(X_{\et},\Z/p\Z)$.
The arrow $\varphi_3$ is defined as the natural
              identification of complexes, and the composite of
             the lower row agrees with
	     the bottom arrow in \eqref{dual:commute4}.
The squares (2) and (3) commute
          in the category of complexes,
     and the squares
       (1) and (4) commute in $D^b(X_{\et},\Z/p\Z)$
          by the definition of connecting morphisms.
Thus the diagram \eqref{dual:commute4} commutes.
\end{pf*}
\begin{pf*}{Proof of Lemma \ref{sublem3}}
There is a commutative diagram in $D^b(X_{\et},\Z/p\Z)$
$$
{\small
\begin{CD}
    (Rj_*\mu_p^{\otimes n-1}[-1]) \otimes ^{\L}
       \T(n')_X @>>> Rj_*\mu_p^{\otimes N+1}[-1]
          @>{\eqref{sigma(n)}}>> \iota_*\nu_Y^N[-N-2]\\
      @V{(\delta_1[-1])\otimes^{\L}\id}VV
         @V{\delta_2[-1]}VV   @| \\
    \iota_*R\iota^!\T(n-1)_X \otimes ^{\L}
         \T(n')_X @>>>
          \iota_*R\iota^!\T(N+1)_X
          @<{-\gys_{\iota}^{N+1}}<{\simeq}<  \iota_*\nu_Y^N[-N-2],
\end{CD}
}
$$
where the left horizontal arrows are product morphisms
    and we wrote $\delta_1$ and $\delta_2$ for
      $\delta^{\loc}_{V,Y}(\T(n-1)_X)$ and $\delta^{\loc}_{V,Y}(\T(N+1)_X)$,
        respectively.
The left square commutes by
          Lemma \ref{sublem3'}, and
    the right square commutes by \eqref{DT:loc}.
By this commutative diagram and
      the anti-commutativity of \eqref{ancn:anti},
       the square
\begin{equation}\label{dual:commute}
\begin{CD}
    (\iota^*\T(n)_X[-1]) \otimes^{\L} \iota^*\T(n')_X
         @>{\mathrm{product}}>>
           \iota^*Rj_*\mu_{p}^{\otimes N+2}[-1]\\
      @V{(a_n[-1])\otimes^{\L}a_{n'}}VV  @VV{\eqref{tau2}}V\\
    (\KK(n)[-1])  \otimes ^{\L} \KK(n')
            @>{\eqref{dual:pair1}}>>
                   \mu'\otimes \nu_Y^N[-N-2]
\end{CD}
\end{equation}
commutes in $D^b(Y_{\et},\Z/p\Z)$
          (cf.\ the diagram in Lemma \ref{sublem2}).
Now we prove Lemma \ref{sublem3},
     using a similar argument as for Lemma \ref{sublem2}.
We first show (1).
Because there is no non-zero morphism
        from ${\Bbb U}(n)$ to
          ${\Bbb D}(\mu' \otimes\nu_{Y}^{n'-2}[-n'+1])$,
        it suffices to show the commutativity of
           the central square in \eqref{dual:commute2}.
Our task is to show that the composite morphism
$$
\begin{CD}
     (\KK(n)[-1])~ \otimes^{\L}
         (\mu' \otimes \nu_{Y}^{n'-2})[-n'+1] \; \lra \;
          (\KK(n)[-1]) \otimes^{\L} \KK(n')
           @>{\eqref{dual:pair1}}>> \mu' \otimes \nu_{Y}^N[-N-2]
\end{CD}
$$
induces $(-1)^n \cdot f_1$,
       which follows from the commutativity of \eqref{dual:commute} and
     Lemmas \ref{sublem1} (2) and \ref{lem:CD}.
The details are straight-forward and left to the reader.
We next show (2).
There are no non-zero morphisms
    from $\mu' \otimes  \nu_Y^{n-2}[-n+1]$ to $\DD(U^1M^{n'}[-n'])$,
     and the left square in \eqref{dual:commute3}
       commutes by a similar argument as for (1).
Hence there is a unique morphism
         $f_4: U^1M^n[-n] \ra \DD(U^1M^{n'}[-n'-1])$
            fitting into \eqref{dual:commute3}
               (cf.\ Lemma \ref{lem:CD2} (2)),
      which necessarily agrees with
         $f_3$ by the commutativity of \eqref{dual:commute}
        and the construction of these maps.
Thus we obtain the lemma.
\end{pf*}
This completes the proof of Proposition \ref{prop:sdual} and
     Theorems \ref{thm:main} and \ref{thm:main2}.
\qed
\subsection{Consequences in the local case}\label{sect10.6}
We state some consequences of Theorem \ref{thm:main}.
Let the notation be as in Theorem \ref{thm:main}
        and Remark \ref{rem:PTdual}.
Let $\H^i_{\ur}(V,\mu_{p^r}^{\otimes n})$ be
    the image of the canonical map
    $\H^i(X,\T_r(n)_X)\ra \H^i(V,\mu_{p^r}^{\otimes n})$.
\begin{cor}\label{cor:dual}
$\H^i_{\ur}(V,\mu_{p^r}^{\otimes n})$
         and $\H^{2d-i}_{\ur}(V,\mu_{p^r}^{\otimes d-n})$
are exact annihilators of each other under
      the non-degenerate pairing $\eqref{PTdual}$.
\end{cor}
\begin{pf}
By Theorem \ref{thm:comp} and a similar argument as for
     Lemma \ref{lem:com},
     one can easily check that the composite map
$$
\begin{CD}
\H^{2d}(V,\mu_{p^r}^{\otimes d}) @>{\delta_{V,Y}^{\loc}(\T_r(d)_X)}>>
    \H^{2d+1}_Y(X,\T_r(d)_X) @>{\tr_{(X,Y)}}>> \Z/p^r\Z
\end{CD}
$$
agrees with the trace map $\tr_V$.
Hence the diagram with exact rows
$$
\begin{CD}
\H^i(X,\T_r(n)_X) @>>>
        \H^i(V,\mu_{p^r}^{\otimes n})
          @>{\delta^{\loc}}>>
           \H^{i+1}_Y(X,\T_r(n)_X)\\
      @V{\eqref{dual:main}}V{\simeq}V  @V{\eqref{PTdual}}V{\simeq}V
         @V{\eqref{dual:main}}V{\simeq}V\\
\H^{2d+1-i}_Y(X,\T_r(d-n)_X)^*
     @>{(\delta^{\loc})^*}>>
        \H^{2d-i}(V,\mu_{p^r}^{\otimes d-n})^* @>>>
           \H^{2d-i}(X,\T_r(d-n)_X)^*
\end{CD}
$$
commutes up to signs.
Here the superscript $*$ means the Pontryagin dual,
     and the bijectivity of the left and the right vertical arrows
        follows from Theorem \ref{thm:main}.
Now the assertion follows from a simple diagram chase on this diagram.
\end{pf}
Corollary \ref{cor:dual} includes some non-trivial duality theorems
        in the local class field theory.
More precisely, we have the following consequence,
       where $K:=\Frac(A)$ and $\Br(C):=\H^2(C,\Gm)$:
\begin{cor}
Let $C$ be a proper smooth curve over $K$
     with semistable reduction.
Then there is a
    non-degenerate pairing of finite $\Z/p^r\Z$-modules
$$
\begin{CD}
\Pic(C)/p^r \times {}_{p^r}\Br(C) \lra \Z/p^r\Z.
\end{CD}
$$
\end{cor}
\noindent
This corollary recovers the $p$-adic part of
      the Lichtenbaum duality \cite{lich} for $C$
            and the Tate duality \cite{tate:duality}
               for the Jacobian variety of $C$
             (cf.\ \cite{saito}, p.\ 413).
However our proof is not new,
   because we use Artin's proper base-change theorem
       for Brauer groups.
\begin{pf}
Take a proper flat regular model $X$ over $B$ of $C$
         with semistable reduction.
Let $Y$ be the closed fiber of $X/B$, and
     define $\Br(X):=\H^2(X,\Gm)$.
There is a commutative diagram with exact rows
$$
\begin{CD}
0 @>>> \Pic(X)/p^r @>>> \H^2(X,\T_r(1)_X)
         @>>> {}_{p^r}\Br(X) @>>> 0 \\
    @.  @VVV  @VVV  @VVV  @.\\
0 @>>> \Pic(C)/p^r @>>> \H^2(C,\mu_{p^r})
         @>>> {}_{p^r}\Br(C) @>>> 0
\end{CD}
$$
See \eqref{DT:weight1} for the upper row.
In view of Corollary \ref{cor:dual},
      our task is to show $\Pic(C)/p^r= \H^2_{\ur}(C,\mu_{p^r})$.
Because the left vertical arrow is surjective,
   it is enough to show $\Br(X)=0$.
Now by Artin's proper base-change theorem:
          $\Br(X) \simeq \H^2(Y,\Gm)$
     (cf.\ \cite{g:brauer}, III.3.1),
    we are reduced to showing $\H^2(Y,\Gm)=0$,
   which follows from the classical Hasse principle for
          the function fields of $Y$
            (cf.\ \cite{saito}, \S3, p.\ 388).
Thus we obtain the corollary.
\end{pf}
\medskip
\noindent
{\bf Acknowledgements.}
\hspace{1pt}
The research for this article was supported by
      JSPS Postdoctoral Fellowship for Research Abroad and EPSRC grant.
The author expresses his gratitude to University of Southern California
     and The University of Nottingham for their great hospitality, and
     to Professors Shuji Saito, Kazuhiro Fujiwara, Wayne Raskind,
     Thomas Geisser and Ivan Fesenko for valuable comments and discussions.
Thanks are also due to Doctor Kei Hagihara,
     who allowed the author to include
     his master's thesis at Tokyo University as Appendix A.

\medskip
\appendix
\theoremstyle{plain}
      \newtheorem{athm}{Theorem}[subsection]
      \newtheorem{alem}[athm]{Lemma}
      \newtheorem{asublem}[athm]{Sublemma}
      \newtheorem{aprop}[athm]{Proposition}
      \newtheorem{acor}[athm]{Corollary}
      \newtheorem{aconj}[athm]{Conjecture}
      \newtheorem{arem}[athm]{Remark}

\section{\bf An application of $p$-adic Hodge theory
    to the coniveau filtration}
\begin{center}
by {\sc kei hagihara}
\footnote{Appendix A is based on his master's thesis at Tokyo University on
1999. He expresses his gratitude to Professors Kazuya Kato and Shuji Saito
for helpful conversations and much encouragement. He is supported by the
21st century COE program at Graduate School of Mathematical Sciences,
University of Tokyo.}
\end{center}
\begin{center}
{\small
Department of Mathematical Sciences,
University of Tokyo
}
\end{center}
\par
\vspace{-3mm}
\begin{center}
{\scriptsize
Komaba, Meguro-ku, Tokyo 153-8914 Japan
}
\end{center}
\vspace{-3mm}
\begin{center}
{\scriptsize
e-mail: kei-hagi@@318uo.ms.u-tokyo.ac.jp
}
\end{center}
\medskip
\stepcounter{subsection}
\noindent
A.1.
Every `suitable' cohomology theory $\H^{*}$ for schemes, for example
\'etale cohomology, is naturally accompanied with an important filtration
called {\it coniveau filtration}, which is defined as follows:
\begin{align*}
N^{r}\H^{i}(X)
:= & \Image\big(\varinjlim_{Z \in X^{\ge r}}\, \H^{i}_{Z}(X)
\longrightarrow \H^{i}(X) \big) \\
    = & \ker\big(\H^{i}(X) \longrightarrow \varinjlim_{Z \in X^{\ge r}}\,
\H^{i}(X-Z)\big),
\end{align*}
where $\H^{*}_{Z}(X)$ denote cohomology groups with support in $Z$ and we put
\[
X^{\ge r}=\{ Z \subset X \;| \;
     \hbox{closed in $X$ and $\codim_{X}(Z) \ge r$} \}
\]
for non-negative integer $r$.
This filtration, built into any cohomology group, is intimately related to
algebraic cycles and often enables us to control their behavior by various
cohomological tools, although the filtration \textit{per se} has not been
well-understood yet. The aim of this appendix is to analyze this
interesting filtration on \'etale cohomology groups by means of $p$-adic
Hodge theory. More precisely, we give an upper bound of it, assuming that
$X$ is a variety over a `$p$-adic' field.

\subsubsection{}\label{hypo}
To state our results more precisely, we fix the setting. Let $A$ be a
henselian discrete valuation ring $A$ whose fraction field $K$ has
characteristic 0 and whose residue field $k$ is perfect of characteristic
$p>0$. Consider the following diagram of schemes:
\begin{equation}\notag
\begin{CD}
     Y      @>{i}>>   {\frak X}  @<{j}<<    X     \\
@VV{\quad\quad\;\;\;\;\square}V  @VV{\qquad\quad\square}V  @VVV    \\
\Spec (k)  @>>>      \Spec (A)       @<<<     \Spec (K), \\
\end{CD}
\end{equation}
where the vertical arrows are proper and flat, and both squares are 
cartesian. We assume that ${\frak X}$ is a regular semistable family over 
$A$, i.e., ${\frak X}$ is regular, $X$ is smooth over $K$ and $Y$ is a 
reduced divisor on ${\frak X}$ with normal crossings. Fix an algebraic 
closure $\ol K$ of $K$, let $\ol{A}$ be the integral closure of $A$ in $\ol 
K$ and let $\ol{k}$ be its residue field. We denote $Y \otimes_k \ol k$, 
${\frak X} \otimes_A \ol A$ and $X \otimes_K \ol K$ by $\ol{Y}$, 
$\ol{{\frak X}}$ and $\ol{X}$, respectively, and write $\ol i$ and $\ol j$ 
for the canonical maps $\ol{Y} \rightarrow \ol{{\frak X}}$ and $\ol{X} 
\rightarrow \ol{{\frak X}}$, respectively. For simplicity we always suppose 
that $\ol{X}$ and $\ol{{\frak X}}$ are connected. Throughout this appendix, 
we use the general notation fixed in \S\S\ref{sect1}.6--\ref{sect1}.7 of 
the main body.

\subsubsection{}
By standard theorems in \'etale cohomology theory, we have spectral sequences
\[
\begin{CD}
E^{a,b}_{2}=\H^{a}(Y,i^{*}R^{b}j_{*}{\Bbb Z}/p^{n}(m)) \Longrightarrow
\H^{a+b}(X,{\Bbb Z}/p^{n}(m)),\\
E^{a,b}_{2}=\H^{a}(\ol{Y},\ol{i}^{*}R^{b}\ol{j}_{*}{\Bbb
Z}/p^{n})\Longrightarrow \H^{a+b}(\ol{X},{\Bbb Z}/p^{n}),
\end{CD}
\]
where ${\Bbb Z}/p^{n}(m)$ denotes the sheaf $\mu_{p^n}^{\otimes m}$ on
$X_{\et}$. We define the filtration $F^{\bullet} \subset \H^{q}(X,{\Bbb
Z}/p^{n}(m))$
as that induced by the former spectral sequence. Alternatively, one can define
\[
F^{r}\H^{q}(X,{\Bbb Z}/p^{n}(m)):=\Image\,(\H^{q}({\frak X}, \tau_{\le
q-r}Rj_{*}{\Bbb Z}/p^{n}(m)) \to \H^{q}(X, {\Bbb Z}/p^{n}(m)) ).
\]
Now we have two filtrations $N^{\bullet}$ and $F^{\bullet}$ on
$\H^{q}(X,{\Bbb Z}/p^{n}(m))$. One defines the filtrations $N^{\bullet}$
and $F^{\bullet}$ on $\H^{i}(\ol{X},{\Bbb Z}/p^{n})$ as well in the same
way.
\subsubsection{}
Our results are stated as follows.
\addtocounter{athm}{3}
\begin{athm}\label{thm:hmain1}
Let $r, s$ and $n$ be non-negative integers with $0 \le r \le s/2$. Then
\[
N^{r}\H^{s}(X, {\Bbb Z}/p^{n}(s-r)) \subset F^{r}\H^{s}(X,{\Bbb Z}/p^{n}(s-r)).
\]
\end{athm}

\begin{athm}\label{thm:hmain2}
Let $r, s$ and $n$ be non-negative integers with $0 \le r \le s/2$. Then
\[
N^{r}\H^{s}(\ol{X}, {\Bbb Z}/p^{n}) \subset F^{r}\H^{s}(\ol{X}, {\Bbb
Z}/p^{n}).
\]
\end{athm}

\begin{arem}
If $r$ is outside of this interval, these assertions are straight-forward
by coniveau spectral sequences $($cf.\ \cite{BO}$)$.
\end{arem}

\addtocounter{subsubsection}{3}
\subsubsection{}
If $r=1$, $s=3$ and ${\frak X}$ is smooth, then Theorem \ref{thm:hmain1} is 
originally due to Langer and Saito (\cite{LS}, 5.4). Their proof is 
$K$-theoretic and reduces the problem to semi-purity of cohomology groups 
with coefficients in ${\cal K}_{2}$-sheaves. On the other hand, our proof 
is $p$-adic Hodge theoretic, i.e., we will reduce the problem to 
semi-purity of cohomology groups with coefficients in \'etale sheaves of 
$p$-adic vanishing cycles.

\subsubsection{}
The filtration $F^{\bullet}$ is highly non-trivial in the $p$-adic 
coefficients case, in contrast with the `$\ell$-adic coefficients' case, 
where for instance in the good reduction case, the corresponding filtration 
is trivial. In fact, as an application of Theorem \ref{thm:hmain2} we will 
prove the following:
\addtocounter{athm}{2}
\begin{acor}\label{coro}
Let $s, r$ and $n$ be non-negative integers with $0 \le r \le s/2$. Assume 
that ${\frak X}$ is ordinary, i.e., $\H^{a}(\ol{Y},{\cal 
B}^{b}_{\ol{Y}})=0$ for all $a$ and $b$ $($see Theorem $\ref{thm:hyodo}$ 
for ${\cal B}^{b}_{\ol{Y}})$. Then we have
\[ {\mathrm{length}}_{{\Bbb Z}/p^{n}}N^{r}\H^{s}(\ol{X}, {\Bbb Z}/p^{n})
\le \sum_{r \le a \le s}{\mathrm{length}}_{W_{n}(\ol{k})}\H^{a}(\ol{Y},
W_{n}\omega^{s-a}_{\ol{Y}}),  \]
where $W_{n}\omega^{\bullet}_{Y}$ denotes the de Rham-Witt complex defined 
in \cite{hyodo}.
\end{acor}

\addtocounter{subsubsection}{1}
\subsubsection{}
Bloch and Esnault \cite{BE} proved that $\Gamma (Y, \Omega^{m}_{Y}) \neq 0
\Rightarrow N^{1}\H^{m}(\ol{X}, {\Bbb Z}/p) \neq \H^{m}(\ol{X}, {\Bbb
Z}/p)$, assuming that $X$ has ordinary good reduction and that the spectral
sequence
\[E^{a,b}_{2}=\H^{a}(\ol{Y},\ol{i}^{*}R^{b}\ol{j}_{*}{\Bbb Z}/p)
\Longrightarrow \H^{a+b}(\ol{X},{\Bbb Z}/p)\]
degenerates at $E_{2}$-terms. Corollary \ref{coro} recovers and generalizes
this fact.

\subsubsection{}
We will prove Theorems \ref{thm:hmain1}, \ref{thm:hmain2} and Corollary
\ref{coro} in \S A.2, \S A.3 and \S A.4 below, respectively.
\par
\bigskip
\stepcounter{subsection}
\noindent
A.2. {\bf Proof of Theorem \ref{thm:hmain1}.}
We first reduce Theorem \ref{thm:hmain1} to Lemma \ref{lem:key} below. For 
$Z \in X^{\ge r}$, let ${\frak Z}$ be the closure of $Z$ in ${\frak X}$. 
There is a commutative diagram
\begin{equation}\label{keyCD}
\hspace{-50pt}
\begin{CD}
\H^{s}_{Z}(X, {\Bbb Z}/p^{n}(s-r)) @>>> \H^{s}(X, {\Bbb Z}/p^{n}(s-r))\\
@AAA  @AAA   \\
\H^{s}_{{\frak Z}}({\frak X}, \tau_{\le s-r}Rj_{*}{\Bbb Z}/p^{n}(s-r))
@>>> \H^{s}({\frak X}, \tau_{\le s-r}Rj_{*}{\Bbb Z}/p^{n}(s-r)).
\end{CD}
\hspace{-50pt}
\end{equation}
Put ${\cal A}:=Rj_{*}{\Bbb Z}/p^{n}(s-r)$. Since $\H^{i}_{{\frak Z}}({\frak 
X}, {\cal A}) \simeq \H^{i}_{Z}(X, {\Bbb Z}/p^{n}(s-r))$, there is a long 
exact sequence
\[
\cdots \to \H^{s}_{{\frak Z}}({\frak X}, \tau_{\le s-r}{\cal A}) \to
\H^{s}_{Z}(X, {\Bbb Z}/p^{n}(s-r)) \to \H^{s}_{{\frak Z}}({\frak X},
\tau_{\ge s-r+1}{\cal A})\to \H^{s+1}_{{\frak Z}}({\frak X}, \tau_{\le
s-r}{\cal A}) \to \cdots .
\]
Now Theorem \ref{thm:hmain1} is reduced to
\addtocounter{athm}{1}
\begin{alem}\label{lem:key}
$\H^{s}_{{\frak Z}}({\frak X},\tau_{\ge s-r+1}{\cal A})=0$ for any $Z \in 
X^{\ge r}$ and any $r,s \in \Z$ as in the theorem.
\end{alem}
\noindent
Indeed, by this lemma the left vertical arrow in \eqref{keyCD} is 
surjective, and Theorem \ref{thm:hmain1} follows from a diagram chase on 
\eqref{keyCD}.

\addtocounter{subsubsection}{2}
\subsubsection{}
The rest of this subsection is devoted to Lemma \ref{lem:key}. The 
following sublemma follows from a simple argument on flatness, whose proof 
is left to the reader:
\addtocounter{athm}{1}

\begin{alem}\label{lem:dim}
For $Z \in X^{\ge r}$, put $Z_{p}:={\frak Z} \otimes_{A} k$ with ${\frak 
Z}$ the closure of $Z \subset {\frak X}$. Then $Z_{p} \in Y^{\ge r}$.
\end{alem}

\addtocounter{subsubsection}{1}
\subsubsection{}
Let $Z$, ${\frak Z}$ and $Z_p$ be as in Lemma \ref{lem:dim}. For $q \in 
\Z$, we put ${\cal C}_{n}(q):=i^{*}Rj_{*}{\Bbb Z}/p^{n}(q)$. Since 
$R^{m}j_{*}{\Bbb Z}/p^{n}(q) \simeq i_*i^*R^{m}j_{*}{\Bbb Z}/p^{n}(q)$ for 
$m > 0$, we have
\[\H^{s}_{{\frak Z}}({\frak X}, \tau_{\ge s-r+1}Rj_{*}{\Bbb Z}/p^{n}(s-r))
\simeq \H^{s}_{Z_{p}}(Y, \tau_{\ge s-r+1}\,{\cal C}_{n}(s-r)).\]
We prove that the right hand side is zero. There are distinguished triangles
\[
\begin{CD}
\tau_{\ge s-r+1}\,{\cal C}_{n-1}(s-r) \longrightarrow \tau_{\ge 
s-r+1}\,{\cal C}_{n}(s-r) \longrightarrow \tau_{\ge s-r+1}\,{\cal 
C}_{1}(s-r),\\
i^{*}R^{m}j_{*}{\Bbb Z}/p\,(s-r)[-m] \longrightarrow \tau_{\ge m}\,{\cal 
C}_{1}(s-r) \longrightarrow \tau_{\ge m+1}\,{\cal C}_{1}(s-r)
\end{CD}
\]
in $D^{+}(Y_{\et})$, where the former triangle is obtained from the short 
exact sequence
\[0 \longrightarrow {\Bbb Z}/p^{n-1}(s-r) \longrightarrow {\Bbb
Z}/p^{n}(s-r) \longrightarrow {\Bbb Z}/p\,(s-r) \longrightarrow 0\]
and, in fact, distinguished because the map $i^{*}R^{s-r}j_{*}{\Bbb
Z}/p^{n}(s-r) \rightarrow i^{*}R^{s-r}j_{*}{\Bbb Z}/p\,(s-r)$ is surjective 
(cf.\ Theorem \ref{thm:hyodo} (1) in the main body). By these distinguished 
triangles and Lemma \ref{lem:dim}, Lemma \ref{lem:key} is reduced to the 
following semi-purity result:
\addtocounter{athm}{1}
\begin{athm}[Semi-purity]\label{lem:fin}
For any $Z_p \in Y^{\geq r}$ and any $a,m,q$ with $a \le r-1$, we have
\[\H^{a}_{Z_{p}}(Y, i^{*}R^{m}j_{*}{\Bbb Z}/p\,(q))=0.\]
\end{athm}
\addtocounter{subsubsection}{1}
\subsubsection{}\label{hred}
By a standard norm argument, Theorem \ref{lem:fin} is reduced to the case 
where $K$ contains primitive $p$-th roots of unity. In this case, we have 
${\Bbb Z}/p\,(q) \simeq {\Bbb Z}/p\,(m)$ on $X$, and it suffices to 
consider the case $m=q$. Hence Theorem \ref{lem:fin} is reduced to the 
vanishing
\[ \H^{a}_{Z_{p}}(Y, \omega^{b}_{Y,\log}) = \H^{a}_{Z_{p}}(Y,
\omega^{b}_{Y}/{\cal B}^{b}_{Y}) = \H^{a}_{Z_{p}}(Y, {\cal B}^{b}_{Y}) = 0
    \hbox{ for any } a,b \hbox{ with } a \le r-1\]
by the Bloch-Kato-Hyodo theorem (cf.\ Theorem \ref{thm:hyodo}). By a 
similar argument as for \cite{milne:duality}, 1.7, the sheaves 
$\omega^{b}_{Y}/{\cal B}^{b}_{Y}$ and ${\cal B}^{b}_{Y}$ are locally free 
$(\O_Y)^p$-modules of finite rank. By \cite{hyodo}, (1.5.1), there is an 
exact sequence
\[\begin{CD}
0 @>>> \omega^b_{Y,\log} @>>> \omega^{b}_{Y} @>{1-C^{-1}}>>
\omega^{b}_{Y}/{\cal B}_Y^b @>>> 0. \end{CD}\]
Therefore we are further reduced to the following lemma:
\addtocounter{athm}{1}
\begin{alem}\label{hclaim1}
Let $\FF$ be a locally free $(\O_Y)^p$-module of finite rank.
Then $\H^{a}_{Z_{p}}(Y, \FF)$ is zero for any $a \le r-1$.
\end{alem}
\addtocounter{subsubsection}{1}
\subsubsection{}
Since the absolute Frobenius morphism $\F_Y:Y \to Y$ is finite, 
$\H^*_{Z_p}(Y,\FF)$ is isomorphic to $\H^*_{Z_p}(Y,\F_{Y*}(\FF))$. Hence we 
are reduced to the case that $\FF$ is a locally free $\O_Y$-module of 
finite rank. Take an \'etale covering $\{ U_{i}\}_{i \in I}$ of $Y$ which 
trivializes $\FF$. By a local-global spectral sequence (\cite{sga4}, V.6.4
(3)), it suffices to prove that
\[{\cal H}^{a}_{Z_{p} \times_{Y} U_{i}}(U_{i}, {\cal O}_{U_{i}})=0 \text{
for any } a \le r-1 \text{ and any } i \in I,\]
where ${\cal H}^*_Z(X,\bullet)$ denotes the sheaf of cohomology groups with 
support (\cite{sga4}, V.6). One can easily check this triviality by the 
comparison theorem on Zariski and \'etale cohomology groups for coherent 
sheaves (\cite{sga4}, VII.4.3) and standard facts on depth (see, e.g., 
\cite{G:LC}, (3.8)), noting that $Y$ and $U_i$ ($i \in I$) are 
Cohen-Macauley (\cite{AK}, VII.4.8). This completes the proof of Claim, 
Lemma \ref{lem:key} and Theorem \ref{thm:hmain1}.
\hfill \qed
\addtocounter{athm}{2}
\begin{acor}
Let $q, r$ and $s$ be integers with $0 \le r \le s/2$. Then
    \[N^{r}\H^{s}(X, {\Bbb Z}/p\,(q)) \subset F^{r}\H^{s}(X, {\Bbb 
Z}/p\,(q)). \]
\end{acor}
\begin{pf}
Indeed the restriction on Tate twists is unnecessary in this $n=1$ case by 
a standard norm argument.
\end{pf}
\par
\bigskip
\stepcounter{subsection}
\noindent
A.3. {\bf Proof of Theorem \ref{thm:hmain2}.}
Because we do not need to care about Tate twists on $\ol X$, the proof 
becomes much simpler. As in the proof of Theorem \ref{thm:hmain1}, it is 
enough to show that
\[
\H^{s}_{Z_{p}}(\ol{Y}, \tau_{\ge s-r+1}\ol{i}^{*}R\ol{j}_{*}{\Bbb Z}/p^{n})=0\]
for arbitrary $Z_p \in \ol{Y}^{\ge r}$. Since ${\Bbb Z}/p^{n}(1) \simeq
{\Bbb Z}/p^{n}$ on $\ol X$ for any $n$, we have an exact sequence
\[0 \longrightarrow \ol{i}^{*}R^{q}\ol{j}_{*}{\Bbb Z}/p^{n-1}(q)
\longrightarrow \ol{i}^{*}R^{q}\ol{j}_{*}{\Bbb Z}/p^{n}(q) \longrightarrow
\ol{i}^{*}R^{q}\ol{j}_{*}{\Bbb Z}/p\,(q) \longrightarrow 0\]
(cf.\,\cite{bk}, p.\ 142, line 9 and \cite{hyodo}, (1.11.1)). Hence it
suffices to show that
\[
\H^{a}_{Z_{p}}(\ol{Y}, \ol{i}^{*}R^{q}\ol{j}_{*}{\Bbb Z}/p\,(q))=0
\]
for any $a, q$ with $0 \le a \le r-1$. Take a finite extension $k_0/k$ over
which $Z_p \subset \ol Y$ is defined, and take a closed subset $Z_{p,0}$ of
$Y \otimes_{k} k_0$ such that $Z_{p,0} \otimes_{k_0} \ol{k} \simeq Z_{p}$
under the isomorphism $(Y \otimes_{k} k_0) \otimes_{k_0} \ol{k} \simeq
\ol{Y}$.
Now let $K'$ be a finite extension of $K$ whose residue field $k'$ contains
$k_0$, and let $A'$ be the integer ring of $K'$. By \cite{sga4}, VII.5.8,
our task is to show that
\[ \H^{a}_{Z_{p,0} \otimes_{k_0} k'}(Y\otimes_{k}k',
{i'}^{*}R^{q}j'_{*}{\Bbb Z}/p\,(q))=0 \]
for any $a, q$ with $0 \le a \le r-1$, where $i'$ (resp.\ $j'$) denotes the 
morphism $Y \otimes_k k' \rightarrow {\frak X} \otimes_A A'$ (resp.\ $X 
\otimes_K K' \rightarrow {\frak X} \otimes_A A'$). This assertion follows 
from the same argument as in Theorem \ref{lem:fin}. Thus we obtain Theorem 
\ref{thm:hmain2}.
\hfill\qed
\par
\bigskip
\stepcounter{subsection}
\noindent
A.4. {\bf Proof of Corollary \ref{coro}.}
Let $\mlogwitt {\ol{Y}} n *$ be the modified logarithmic Hodge-Witt sheaves
(cf.\ \S\ref{sect3.3} of the main body). The ordinarity assumption implies
that
\[\H^{a}(\ol{Y}, \ol{i}^{*}R^{b}\ol{j}_{*}{\Bbb Z}/p^{n}) \simeq
\H^{a}(\ol{Y}, \mlogwitt {\ol{Y}} n b)\]
(\cite{bk}, (9.2), \cite{hyodo}, (1.10))
and that
\[\H^{a}(\ol{Y},\mlogwitt {\ol{Y}} n b)\otimes_{{\Bbb
Z}/p^{n}}W_{n}(\ol{k}) \simeq \H^{a}(\ol{Y},\mwitt {\ol{Y}} n b) \]
(\cite{bk}, (7.3), \cite{I2}, (2.3)). Hence Corollary \ref{coro} follows
from Theorem \ref{thm:hmain2}.
\hfill\qed
\medskip
\begin{arem}
The theorem of Bloch-Esnault in \cite{BE}, $(1.2)$ is a direct consequence
of Corollary $\ref{coro}$ with $r=n=1$ and ${\frak X}/A$ smooth. They
derived some interesting results on algebraic cycles from this case.
Therefore Corollary $\ref{coro}$ would provide us with much more
information.
\end{arem}

\end{document}